\documentclass[11pt]{article}
\usepackage[utf8]{inputenc}
\usepackage{fullpage,amsmath,amssymb,bm,url,mathrsfs}
\usepackage{algorithmic}
\usepackage{algorithm}
\usepackage{amsthm}
\usepackage{hyperref}
\usepackage{makecell}
\usepackage{tablefootnote}
\usepackage{comment}
\usepackage{soul}
\usepackage{graphics,graphicx}
\usepackage{subcaption,caption,cleveref}
\hypersetup{
    colorlinks=true,
    linkcolor=blue,
    citecolor=blue,
    filecolor=magenta,      
    urlcolor=cyan,
}

\usepackage[round,colon,authoryear]{natbib}
\usepackage{xcolor}
\usepackage{dsfont}
\usepackage{bbm}
\usepackage{enumitem}

\def\calM{{\mathcal M}}
\def\calN{{\mathcal N}}

\def\calP{{\mathcal P}}

\def\bcalX{{\boldsymbol{\mathcal X}}}

\def\EE{{\mathbb E}}

\def\PP{{\mathbb P}}
\def\QQ{{\mathbb Q}}
\def\RR{{\mathbb R}}

\def\v{{\boldsymbol v}}

\def\bmu{{\boldsymbol \mu}}

\DeclareMathOperator*{\argmin}{arg\,min}

\def\bs{\boldsymbol{s}}

\def\calM{{\cal  M}} 
\def\calN{{\cal  N}} 
 
\def\calP{{\cal  P}}

\newcommand{\bfm}[1]{\ensuremath{\mathbf{#1}}}

   \def\bE{\bfm E}  \def\EE{\mathbb{E}}

   \def\bI{\bfm I}

   \def\bM{\bfm M}  
     
   \def\bO{\bfm O}  
     \def\PP{\mathbb{P}}
     \def\QQ{\mathbb{Q}}
     \def\RR{\mathbb{R}}
\def\bs{\bfm s}     
     
\def\bu{\bfm u}   \def\bU{\bfm U}  
\def\bv{\bfm v}   \def\bV{\bfm V}  
     
\def\bx{\bfm x}   \def\bX{\bfm X}  
     
\def\bz{\bfm z}   \def\bZ{\bfm Z}  
                   
  \def\bTheta{\bfm \Theta}

\def\bSigma{\bfm \Sigma}
\def\bLambda{\bfm \Lambda}

\def\MLE{\textsf{\tiny MLE}}

\def\hat{\widehat}

\newcommand\fro[1]{\| #1 \|_{\rm F}}
\newcommand\op[1]{\|#1\|}

\newtheorem{theorem}{Theorem}

\newtheorem{lemma}{Lemma}

\newtheorem{conjecture}{Conjecture}
\theoremstyle{remark}

\newcommand{\E}{\mathbb{E}}
\newcommand{\Prob}{\mathbb{P}}

\newcommand{\eps}{\varepsilon}

\begin{document}

\title{Optimal Estimation and Computational Limit of Low-rank Gaussian Mixtures}

\author{Zhongyuan Lyu and Dong Xia\footnote{Dong Xia's research was partially supported by Hong Kong RGC Grant ECS 26302019, GRF 16303320 and GRF 16300121.}\\
%Department of Mathematics, 
{ Hong Kong University of Science and Technology}}

\date{(\today)}

\maketitle

\begin{abstract}

Structural matrix-variate observations routinely arise in diverse fields such as multi-layer network analysis and brain image clustering. While data of this type have been extensively investigated with fruitful outcomes being delivered, the fundamental questions like its statistical optimality and computational limit are largely under-explored. In this paper, we propose a low-rank Gaussian mixture model (LrMM) assuming each matrix-valued observation has a planted low-rank structure. Minimax lower bounds for estimating the underlying low-rank matrix are established allowing a whole range of sample sizes and signal strength. Under a minimal condition on signal strength, referred to as the {\it information-theoretical limit} or {\it statistical limit}, we prove the minimax optimality of a maximum likelihood estimator which, in general, is computationally infeasible. If the signal is stronger than a certain threshold, called the {\it computational limit}, we design a computationally fast estimator based on spectral aggregation and demonstrate its minimax optimality. Moreover, when the signal strength is smaller than the computational limit, we provide evidences based on the low-degree likelihood ratio framework to claim that no polynomial-time algorithm can consistently recover the underlying low-rank matrix. Our results reveal multiple phase transitions in the minimax error rates and the statistical-to-computational gap.  Numerical experiments confirm our theoretical findings.   We further showcase the merit of our spectral aggregation method on the worldwide food trading dataset.  

\end{abstract}

\section{Introduction}\label{sec:intro}
The recent decade has witnessed a burgeoning demand in processing and analyzing large-scale matrix-variate data which routinely arise in diverse fields.  In gene expression analysis, e.g., the BHL (brain, heart and lung) dataset \citep{BHL-data,mai2021doubly}, the measurement of gene expression on different types of  tissues is often repeated for multiple times. The resultant observation for each tissue becomes a matrix and thus the cluster analysis is operated on matrix-valued observations.  A multi-layer network \citep{le2018estimating,jing2021community,lyu2021latent,paul2020spectral} usually consists of multiple networks on the same set of vertices. Since each observed layer is equivalently represented as an adjacent matrix, problems such as community detection \citep{paul2020spectral}, layer clustering \citep{jing2021community} and common probability matrix estimation \citep{le2018estimating} are generally attacked by statistical analysis on a collection of adjacency matrices. Other notable examples include brain image clustering \citep{sun2019dynamic,wang2017generalized}, EEG data analysis \citep{hu2020matrix,gao2021regularized}, etc. Oftentimes, the dimensions of observed matrices are ultra-large or the number of matrix-valued observations is relatively small, which has motivated the exploration of hidden low-dimensional structures, e.g. sparsity and low-rankness, in matrix-valued observations. All the aforementioned works assumed, among others, certain types of low-rank structures for the underlying parameters of interest and have delivered fruitful outcomes in real-world applications.

Inspired by those foregoing works, throughout this paper, we assume that {\it each matrix-valued observation has a low-rank expectation which might vary for different observations}. Towards that end,  several specific low-rank statistical models, tailored for concrete applications, and respective estimating procedures have been proposed. For instance, a {\it mixture} multi-layer stochastic block model (SBM) was introduced in \cite{jing2021community} for uncovering the global and local communities in multi-layer networks. At the core of this model is the assumption that every layer has a low-rank expected adjacency matrix. Their estimator was based on the (regularized) low-rank tensor decomposition. A special multi-layer SBM was proposed by \cite{paul2020spectral} and estimated by a spectral method. In order to analyze the brain fMRI data, \cite{sun2019dynamic} proposed a tensor Gaussian mixture model and designed an estimator via (fusedly-)truncated low-rank tensor decomposition. Despite these prior efforts, usually motivated by particular applications,  on the low-rank estimates from a mixture of matrix-valued observations, many fundamental questions remain unanswered. What is the role and benefit of low-rankness? How do the sample size and signal strength (see the definition after eq.(\ref{eq:calMr})) characterize the intrinsic difficulty, i.e., are there any phase transitions? What is the statistically optimal rate, which estimator can achieve the rate and is this estimator computationally feasible? What is the fastest error rate achievable by estimators requiring only polynomial-time complexity? This paper aims to answer all these questions and provides a complete picture for the statistical and computational limits in the low-rank estimation from a {\it mixture} of matrix-valued observations. 

We now introduce the {\it low-rank Gaussian mixture model} (LrMM) to formalize the questions. For simplicity, we focus on the mixture of two components and will briefly discuss the case of multiple components in Section~\ref{sec:discuss}. The $d_1\times d_2$ matrix $\bX$ is said to follow an isotropic matrix normal \citep{gupta2018matrix} distribution $\calN(\bM, \bI_{d_1}\otimes \bI_{d_2})$ if ${\it vec}(\bX)\sim \calN\big({\it vec}(\bM),\bI_{d_1d_2}\big)$, where $\bI_d$ represents the $d\times d$ identity matrix and $\bM$ is a deterministic matrix. Clearly, this implies that $\bX$ is equal to $\bM+\bZ$ {\it in distribution} where $\bZ$ has i.i.d. standard normal entries. Denote\footnote{With a slight abuse of notation, we also denote $p_{\bM}$ the associated probability density function.}
\begin{equation}\label{eq:LrMM}
p_{\bM}=\frac{1}{2}\calN(\bM,\bI_{d_1}\otimes \bI_{d_2})+\frac{1}{2} \calN(-\bM, \bI_{d_1}\otimes \bI_{d_2})
\end{equation}
the symmetric mixture of two-component Gaussian mixture model. Then $\bX\sim p_{\bM}$ means that  $\bX$ is sampled from $\calN(\bM,\bI_{d_1}\otimes \bI_{d_2})$ and $\calN(-\bM,\bI_{d_1}\otimes \bI_{d_2})$ with probability both $1/2$, respectively. Put it differently, $\bX$ equals $s\bM+\bZ$ in distribution with $s$ being a Rademacher random variable, called the {\it label} of $\bX$, satisfying $\PP(s=\pm 1)=1/2$. Throughout the paper, we assume that $\bM$ has a small rank, i.e., $r={\rm rank}(\bM)\ll \min\{d_1,d_2\}$. Note that, under model (\ref{eq:LrMM}), the marginal expectation of $\bX$ is actually zero. The former claim of {\it low-rank expectation} in the last paragraph actually refers to the conditional expectation $\EE(\bX|s)=s\bM$ which is low-rank. 
 We remark that the condition of equal prior probabilities is not essential and can be slightly relaxed. The assumption of symmetry of the two components is only for ease of exposition. If the two components have distinct mean matrices, say $\bM_1$ and $\bM_2$, respectively, one can first estimate the average $(\bM_1+\bM_2)/2$, subtract it from all observations and reduce the problem to the symmetric case. Similarly, the assumption of isotropic noise is relaxable as long as the covariance tensor is known. The case of unknown covariance is much more challenging  \citep{davis2021clustering,bakshi2020robustly,belkin2010toward,cai2019chime,ge2015learning,moitra2010settling} even in the vector case and is beyond the scope of the current paper. 

Given i.i.d. observations $\bX, \bX_1,\cdots,\bX_n$ sampled from the mixture distribution $p_{\bM}$ in (\ref{eq:LrMM}), our goals are to estimate the latent low-rank matrix $\bM$, establish the minimax error rates and design computationally efficient estimators. We assume  $d_1\asymp d_2\asymp d$ meaning that there exist absolute constants $c_0, C_0>0$ satisfying $c_0d\leq \min\{d_1,d_2\}\leq \max\{d_1,d_2\}\leq C_0d$. The parameter space of interest is, for any $\lambda>0$, 
\begin{equation}\label{eq:calMr}
\calM_{d_1,d_2}(r,\lambda):=\left\{\bM\in\RR^{d_1\times d_2}: {\rm rank}(\bM)= r,\ \lambda\asymp \sigma_r(\bM)\leq \cdots\leq \sigma_1(\bM)\asymp \lambda\right\}
\end{equation}
where $\sigma_k(\cdot)$ denotes the $k$-th largest singular value of a matrix. For notational brevity, we shall write $\calM(r,\lambda)$ for short. The {\it signal strength} of low-rank models is usually determined by the smallest non-zero singular value \citep{koltchinskii2016perturbation,zhang2018tensor,xia2021normal,cheng2021tackling,gavish2014optimal}. The set $\calM_{d_1,d_2}(r,\lambda)$ is the collection of all $d_1\times d_2$ rank-$r$ matrices whose signal strength is of order $\lambda$. For simplicity, we only focus on the well-conditioned matrices, i.e., with a bounded condition number. The minimax error rate of estimating $\bM$ is defined by 
$
\inf_{\hat\bM}\sup_{\bM\in\calM(r,\lambda)} \EE\ell(\hat\bM, \bM)
$
, where the infimum is taken over all possible estimator $\hat\bM$ constructed from the i.i.d. observations $\bX_1,\cdots,\bX_n$ and the loss function is $\ell(\hat\bM,\bM):=\min_{\eta=\pm1}\|\hat\bM-\eta\bM\|_{\rm F}$ with $\|\cdot\|_{\rm F}$ standing for the Frobenius norm. Note that, due to the symmetry of model (\ref{eq:LrMM}), $\bM$ is estimable up to a sign flip. 

If $\bM$ has a full rank with $r=\min\{d_1,d_2\}$, model (\ref{eq:LrMM}) reduces to the canonical two component isotropic Gaussian mixture model (GMM) in the dimension $d_1d_2\asymp d^2$, which has been extensively investigated in the literature. See \cite{balakrishnan2017statistical, chen1995optimal,ho2016convergence,xu2016global,wu2020optimal} and references therein. For instance, \cite{wu2019randomly} proved that the minimax rate \footnote{Note that there is an additional term $d^2(\theta n)^{-1}$ derived in \cite{wu2019randomly} which is actually negligible if inspecting all other terms carefully.} is
\begin{align}\label{eq:GMM_bd}
\inf_{\hat\bM}\sup_{\|\bM\|_{\rm F}=\theta} \EE \ell(\hat \bM ,\bM)\asymp \min\left\{\frac{1}{\theta}\frac{d}{n^{1/2}}+\frac{d}{n^{1/2}},\ \theta\right\}
\end{align}
, and showed that a simple spectral method, together with a trivial estimate for the case of small $\theta$, is minimax optimal. This rate implies intriguing phenomenons of phase transitions concerning the sample size $n$ and signal strength $\theta$. For instance, if the sample size $n\geq d^2$, their result reveals three different minimax rates: $\theta$ for $\theta\leq d^{1/2}n^{-1/4}$, $\theta^{-1}dn^{-1/2}$ for $d^{1/2}n^{-1/4}\leq \theta\leq 1$ and $dn^{-1/2}$ for $\theta\geq 1$.  Interestingly,  it also implies that non-trivial estimate is impossible, i.e., information-theoretically impossible, if the signal strength is smaller than $d^{1/2}n^{-1/4}$.  
Undoubtedly,  if $\bM$ is low-rank with $r\ll d$,  one can naturally foresee the existence of multiple phase transitions for the minimax error rates.  Establishing these rates becomes more challenging for several reasons.  On the methodological front,  a naive spectral method cannot attain the minimax optimal rate and thus additional procedures are necessary.  On the theoretical front,  the low-rank structure dictates a smaller intrinsic dimension and brings about new behaviors to the phase transitions of the minimax error rates.  See,  e.g. \cite{koltchinskii2015optimal,ma2015volume} and references therein.  On the computational front,  it is well recognized that the low-rankness sometimes bears a so-called {\it statistical-to-computational gap} \citep{barak2016noisy, zhang2018tensor} in the sense that there exist regimes where statistically optimal estimators can be computationally infeasible,  e.g.,  requiring an exponential-time complexity. 

\begin{figure}
	\center
	\includegraphics[width=0.98\linewidth]{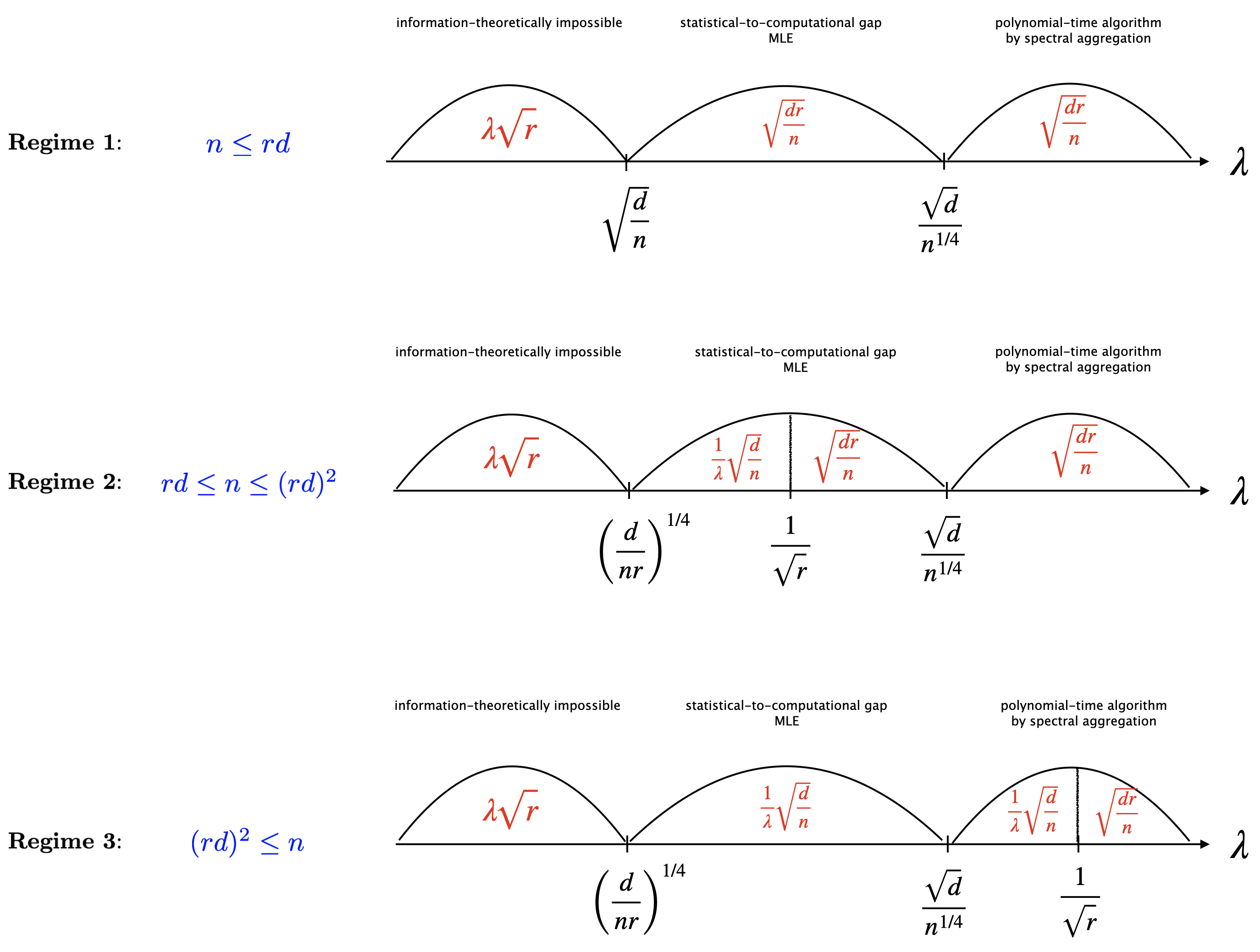}
	\caption{The minimax rates, phase transitions and statistical-to-computational gaps of LrMM,  model (\ref{eq:LrMM}). Here $r$ is the rank, the matrix dimension $d_1\asymp d_2\asymp d$, $n$ is the sample size and $\lambda$ denotes the smallest non-zero singular value. There exist three regimes concerning the sample size which are colored in {\color{blue}blue}. The minimax error rates (up to logarithmic factor) of estimating $\bM\in \calM(r,\lambda)$  in different regimes are colored in {\color{red}red}.  Here {\it information-theoretically impossible} means that non-trivial estimates are impossible because of weak signal strength.  Within the low-degree likelihood ratio framework \citep{kunisky2019notes},  we provide evidence showing that no polynomial-time algorithms can consistently estimate $\bM$ if $\lambda$ is smaller than $d^{1/2}n^{-1/4}$.}
	\label{fig:minimax-rate}
\end{figure}

The summary of our contributions is as follows. We establish the minimax rate of estimating the rank-$r$ matrix $\bM$ for the LrMM model that reads as
\begin{align}\label{eq:minimax-rate}
\inf_{\hat\bM}\ \sup_{\bM\in \calM_{d_1,d_2}(r,\lambda)} \EE \ell(\hat\bM, \bM)\ \asymp \ \min\left\{\frac{1}{\lambda}\Big(\frac{d}{n}\Big)^{1/2}+\Big(\frac{dr}{n}\Big)^{1/2},\ \lambda r^{1/2} \right\}
\end{align}
where the infimum is taken over all possible estimators, regardless of their computational feasibility. This rate implies that, when the sample size $n\geq dr$, it is information-theoretically impossible to estimate $\bM$ if the signal strength $\lambda$ is smaller than $d^{1/4}(rn)^{-1/4}+(d/n)^{1/2}$. Under minimal conditions, we prove that the maximum likelihood estimator (MLE) can achieve the rate (\ref{eq:minimax-rate}) up to a logarithmic factor.  
Unfortunately, there are no known polynomial-time algorithms with guaranteed performance to solve MLE. Earlier works \citep{tosh2017maximum,sanjeev2001learning} show that solving MLE is generally NP-hard. We then propose a computationally fast estimator based on spectral aggregation. This approach can be viewed as a modified method of second moment \citep{pearson1894contributions, wu2020optimal} adapted with a spectral projection to leverage the low-rank structure. We prove that this computationally efficient estimator can achieve the minimax rate (\ref{eq:minimax-rate}) as long as the signal strength $\lambda$ is larger than $d^{1/2}n^{-1/4}$, which is much stronger than the information-theoretical requirement for the minimal signal strength. This difference unveils the statistical-to-computational gap in LrMM. Lastly, we adopt the low-degree likelihood ratio framework \citep{kunisky2019notes} to {\it conjecture} that no polynomial-time estimator is consistent if $\lambda$ is smaller than $d^{1/2}n^{-1/4}$. 
The minimax rates, phase transitions and statistical-to-computational gaps are illustrated in Figure~\ref{fig:minimax-rate}. 

Our results are closely related yet crucially different from several existing works. In \cite{chen2021learning}, a low-rank mixture model was proposed for linear regression which is generally more challenging than our model (\ref{eq:LrMM}). They designed a computationally efficient estimator but provided no results respecting the statistical optimality or computational limits.  A multi-graph network model was introduced by \cite{wang2019common} which allows heterogeneous structure on each matrix-valued observation.  However,  their model has no mixture nature and there is no guarantee on minimax optimality.   More recently,  \cite{jing2021community} proposes a mixture multi-layer SBM and establishes the minimax error rate of spectral estimate only for the special regime when the sample size $n$ is smaller than $d$ and the signal strength,  reflected by the network sparsity,  is strong enough.  In addition,  our LrMM is directly related to low-rank tensor literature.  By stacking the matrix observations slice by slice,  we end up with a tensor of size $n\times d_1\times d_2$ whose expectation,  under model (\ref{eq:LrMM}),  has a low Tucker rank $(1,r,r)$.  See,  e.g.,  \cite{zhang2018tensor, jing2021community} and references therein.  Minimax rates for low-rank tensor denoising and noisy tensor completion have been investigated by \cite{zhang2018tensor} and \cite{xia2021statistically}, respectively.  However, they both require the sample size $n$ to be of the same order of $d$, which becomes unrealistic in the low-rank mixture model. Finally, it worths to remark that our bound (\ref{eq:minimax-rate}) reduces to the minimax bound of GMM (\ref{eq:GMM_bd}) if we let $\bM$ be full-rank. To see this, one can just replace $\lambda$ and $r$ in our bound (\ref{eq:minimax-rate}) by $\theta d^{-1/2}$ and $d$, respectively. 

The rest of paper is organized as follows. We establish the minimax lower bound in Section~\ref{sec:MLE} and prove that the maximum likelihood estimator, albeit computationally infeasible in general, achieves the minimax optimal rates. A computationally fast estimator based on spectral aggregation is proposed in Section~\ref{sec:spectral} which attains minimax optimal rates as long as the signal strength is strong. Section~\ref{sec:stat-to-comp-gap} justifies the statistical-to-computational gap by showing that there exists some regime where the MLE can attains minimax optimal rates but no-polynomial time algorithms can consistently recover the underlying low-rank matrix. We then showcase results of numerical simulations in Section~\ref{sec:numerical}, present a real-world data experiment in Section~\ref{sec:realdata},  and discuss open questions and potential directions in Section~\ref{sec:discuss}. 

\section{Maximum likelihood estimator and minimax optimality}\label{sec:MLE}
We slightly abuse the notation and denote $p_{\bM}(\cdot)$ the probability density function of $\bX\in\RR^{d_1\times d_2}$ under the LrMM model (\ref{eq:LrMM}). The family of density functions parameterized by  $\calM_{d_1,d_2}(r,\lambda)$ is written as (note that we assume $d_1\asymp d_2\asymp d$)
$$
\calP_{d_1,d_2}(r,\lambda):=\Big\{p_{\bM}: \bM\in\calM_{d_1,d_2}(r,\lambda)\Big\}
$$
which is indexed by rank-$r$ matrices with the signal strength $\lambda$.  Given $i.i.d.$ observations $\bX_1, \cdots, \bX_n$ sampled from $p_{\bM}$, the maximum likelihood estimator (not necessarily unique) is defined by 
\begin{align}\label{eq:MLE}
p_{\hat\bM_{\MLE}}:=\underset{p_{\bM}\in\calP_{d_1,d_2}(r,\lambda)}{\arg\max}\ \sum_{i=1}^n \log\big(p_{\bM}(\bX_i)\big)
\end{align}
While the MLE estimator (\ref{eq:MLE}) is generally NP-hard to compute, it often serves as a benchmark for understanding the information-theoretical limit of a statistical model. 

We begin with the regime $n=\tilde \Omega (dr)$\footnote{Here, $\tilde{\Omega}$ stands for the standard big-$\Omega$ notation up to a logarithmic factor.}, which falls into the typical low-dimensional setting\footnote{The low-dimensional setting here refers to the case that dimension $d$ is allowed to grow with sample size  $n$, while the order of $n$ still dominates.}.  The convergence rate of MLE in this regime has been thoroughly investigated for Gaussian mixture model.  See, for instance, \cite{leroux1992consistent,van1993hellinger,chen1995optimal,genovese2000rates,ghosal2001entropies}. The standard tool, e.g. \cite{van1993hellinger} and \cite[Theorem 7.4]{geer2000empirical}, establishes the convergence rate of MLE in the Hellinger distance defined by $d_{H}(p_{\bM_1}, p_{\bM_2}):=1-\int p_{\bM_1}^{1/2}(\bX)p_{\bM_2}^{1/2}(\bX)d\bX$ for two density functions $p_{\bM_1}(\cdot)$ and $p_{\bM_2}(\cdot)$. According to this tool, it suffices to bound the bracketing entropy number of a class of square root density functions around the truth $p_{\bM}^{1/2}$.  While existing literature \citep{ho2016convergence,ho2016strong,maugis2011non} have developed respective bracketing entropy bounds for Gaussian mixture model, they only focus on the fixed dimension $d$ and their method is inapplicable to matrix-variate observations with a planted low-rank structure. By a covering argument and the construction of bracket functions, we establish such a bracketing entropy bounds for LrMM and derive the upper bound in Hellinger distance for $d_{\textsf{H}}(\hat p_{\bM_{\MLE}}, p_{\bM})$. 

To bridge the density estimation and parameter estimation, we resort to a sharp characterization for the total variation distance (similarly, the Hellinger distance) between Gaussian mixture densities established recently by \cite{davies2021lower}. 
\begin{lemma}\label{lem:hellinger} (Lower bound of Hellinger distance)
Let $\bM_1$ and $\bM$ be two matrices, and denote $p_{\bM_1}$ and $p_{\bM}$ the two density functions defined by (\ref{eq:LrMM}). There exists absolute constants $c_0,c_1,c_2>0$ such that, if $\fro{\bM}+\fro{\bM_1}\le c_0$ then
$$
d_{H}(p_{\bM_1}, p_{\bM})\geq c_1\big(\|\bM\|_{\rm F}+ \|\bM_1\|_{\rm F}\big)\ell(\bM_1, \bM)
$$
Otherwise
$$
d_{H}(p_{\bM_1}, p_{\bM})\geq c_2\min\big\{1, \ell(\bM_1, \bM)\big\}
$$
where $\ell(\bM_1,\bM):=\min\{\|\bM_1-\bM\|_{\rm F},\ \|\bM_1+\bM\|_{\rm F}\}$.
\end{lemma}

Together with the upper bound of Hellinger distance $d_{\textsf{H}}(\hat p_{\bM_{\MLE}}, p_{\bM})$ and Lemma~\ref{lem:hellinger}, we obtain the error rate of the maximum likelihood estimator when $n=\tilde \Omega (dr)$,  namely the first part of Theorem \ref{thm:mle}. 

However, the above argument fails when it comes to the regime $n=\tilde O(dr)$\footnote{Again, $\tilde{O}$ stands for the standard big-$O$ notation up to a logarithmic factor.}, corresponding to an ultra high-dimensional setting.  The reason is that the minimax lower bound,  as we will see later in Theorem \ref{thm:lower_bound},  suggests that the optimal error rate should be of order  $(dr/n)^{1/2}$, which can be larger than 1.  Consequently,  the Hellinger distance is no longer an appropriate metric\footnote{The error rate of other bounded metric, say, the Wasserstein distance considered in \cite{doss2020optimal}, also becomes trivial when $d>n$.}, for instance, the lower bound in Lemma \ref{lem:hellinger} becomes trivial.  To this end,  we turn to Kullback-Leibler (KL) divergence defined by $D_{\textsf{KL}}(p_{\bM_1}\|p_{\bM_2}):=\int p_{\bM_1}(\bX)\log(p_{\bM_1}(\bX)/p_{\bM_2}(\bX))d\bX$.  Though KL divergence is not a metric itself, in many cases its convergence also implies consistency of parameter estimate in some metric of interest \citep{geer2000empirical}. Moreover, the KL divergence in its form is closely related to MLE and its unboundedness property is beneficial for our purpose since $(dr/n)^{1/'2}$ possibly diverges. By carefully characterizing the distribution of $\log(p_{\bM_1}(\bX)/p_{\bM_2}(\bX))$ and exploiting the concentration inequality of suprema of an empirical process \cite[Theorem 4]{adamczak2008tail}, we are able to derive an upper bound for the KL divergence $D_{\textsf{KL}}(p_{\bM}\|p_{\hat\bM_{\MLE}})$. We also establish the following lower bound relating KL divergence to the distance in the parameter space.  Combining Lemma~\ref{lem:KL} with the upper bound of $D_{\textsf{KL}}(p_{\bM}\| p_{\hat\bM_{\MLE}})$ leads to the desired error rate in the regime $n=\tilde O(dr)$,  i.e.,  the second part of Theorem \ref{thm:mle}.
\begin{lemma}\label{lem:KL}(Lower bound of KL divergence)
Let $\bM_1$ and $\bM$ be two matrices, and denote $p_{\bM_1}$ and $p_{\bM}$ the two density functions defined by (\ref{eq:LrMM}). There exists absolute constants $C_0,C_1>1,c_0>0$ such that if $\fro{\bM}\ge C_0$ and $\fro{\bM-\bM_1}\ge C_1$,  then
	$$D_{\textsf{KL}}(p_{\bM}\| p_{\bM_1})\ge c_0\cdot \ell^2(\bM,\bM_1)$$
\end{lemma}
Collecting two pieces, the error rate of the maximum likelihood estimator is summarized in the following theorem.  

\begin{theorem}\label{thm:mle}
Suppose $\bM\in\calM(r,\lambda)$ and let $\hat\bM_{\MLE}$ denote the maximum likelihood estimator by (\ref{eq:MLE}).  
%Under either of the following scenarios:
%\begin{itemize}
%	\item $n\le d$, $\lambda\ge C_1\sqrt{\frac{d}{n}}$ for some absolute constant $C_1$
%	\item $n\geq d$, $\lambda\ge C_2\left(\frac{d}{n}\right)^{\frac{1}{4}}$ for some absolute constant $C_2$
%\end{itemize}
\begin{enumerate}
	\item[(1)] If $dr\log nd<n$, then there exist absolute constants $c_1,c_2, C_1,C_2,C_3>0$ such that the following bound holds with probability at least $1-\exp(-c_1d\log^2(nd))$, 
\begin{align}\label{eq:MLE_bd1}
\ell (\hat \bM_{\MLE}, \bM)\le C_1\left(\sqrt{\frac{dr\log (nd)}{n}} + \frac{1}{\lambda}\sqrt{\frac{d\log (nd)}{n}}\right)
\end{align}
If further assume $\lambda\le C_2\exp(c_2d\log^2(nd))$, then
$$
\E \ell (\hat \bM_{\MLE}, \bM)\le C_3\left(\sqrt{\frac{dr\log (nd)}{n}} + \frac{1}{\lambda}\sqrt{\frac{d\log (nd)}{n}}\right)
$$
\item[(2)] If $dr\log nd\ge n$, then there exist absolute constants $C_4,C_5,C_6,C_7>0$ such that if \\$C_4r^{-1/2}\le \lambda\le C_5d^{1/2}$, then the following bound holds with probability at least $1-(nd)^{-4}$, 
\begin{align}\label{eq:MLE_bd2}
\ell (\hat \bM_{\MLE}, \bM)\le C_6\sqrt{\frac{dr\log (nd)}{n}}
\end{align}
And the following bound in expectation holds,
$$
\E \ell (\hat \bM_{\MLE}, \bM)\le C_7\sqrt{\frac{dr\log (nd)}{n}}
$$
\end{enumerate}
\end{theorem}

We note that the logarithmic factor in (\ref{eq:MLE_bd1}) emerges from the bracketing entropy bound and that in (\ref{eq:MLE_bd2}) arises from the tail inequality for suprema of empirical processes of unbounded functions. The high probability bound in the first part of Theorem~\ref{thm:mle} is proved without conditions on the sample size $n$, the rank $r$ or on the signal strength $\lambda$. It suggests intriguing phase transitions in the regime $n=\tilde \Omega (dr)$. When $\lambda>r^{-1/2}$, the MLE attains the rate $\tilde{O}\big((rd/n)^{1/2}\big)$, growing with respect to the rank $r$, which is the best achievable rate even if the labels of observations are all known, namely in the {\it oracle} scenario.  On the other hand, if $\lambda< r^{-1/2}$, the MLE attains the rate $\tilde{O}\big(\lambda^{-1}(d/n)^{1/2}\big)$ that is free of the underlying rank $r$. Moreover, a trivial estimate by $\hat\bM={\bf 0}$ attains the error rate $r^{1/2}\lambda$. Therefore, the MLE becomes pointless if $\lambda$ is smaller than $d^{1/4}(rn)^{-1/4}+(d/n)^{1/2}$, which is referred to as the information-theoretically impossible regime. In the second statement of Theorem \ref{thm:mle}, a more stringent condition is imposed on signal strength ($\lambda=O(d^{1/2})$) for technical difficulty, though we believe that MLE could attain the optimal rate $\tilde O((rd/n)^{1/2})$ in a wider range of $\lambda$ via more sophisticated analysis. On the other hand, as long as $\lambda=\Omega(d^{1/2}n^{-1/4})$, a computationally efficient estimator (see Section~\ref{sec:spectral}) is already able to attain the optimal rate. As we intend to reveal the optimal estimation rate under different signal strength, we only appeal to MLE when the signal strength is not strong enough. Therefore, the technical condition of $\lambda$ for MLE is not essential. 

The next theorem demonstrates the minimax optimality of the MLE by establishing a matching minimax lower bound up to the logarithmic factor. We note that the minimax lower bound (\ref{eq:lower_bd}) is a statistical lower bound because it takes no considerations of the computational feasibility. In Section~\ref{sec:spectral}, we introduce a computationally fast estimator that achieves these lower bounds but requires much more stringent conditions. 

\begin{theorem}\label{thm:lower_bound}
There exists an absolute constant $c_1>0$ such that
\begin{equation}\label{eq:lower_bd}
\inf_{\hat \bM}\ \sup_{\bM\in \mathcal{M}(r, \lambda)}\E\ell(\hat \bM, \bM)
%\geq c_1\left(\sqrt{\frac{d/n}{\lambda^2\wedge 1}}\wedge \lambda\right)
\geq c_1\left(\sqrt{\frac{dr}{n}}+\frac{1}{\lambda}\sqrt{\frac{d}{n}}\right)\bigwedge \lambda\sqrt{r},
\end{equation}
where the infimum is taken over all possible estimators and $a\wedge b=\min\{a, b\}$. 
\end{theorem}

\section{Computationally efficient estimator by spectral aggregation}\label{sec:spectral}
Since the MLE (\ref{eq:MLE}) is generally computationally infeasible, it is of crucial importance to design an estimator which is polynomial-time computable. While existing works have demonstrated the optimality of spectral method for both estimation \citep{wu2019randomly} and clustering \citep{loffler2019optimality} under the GMM, it turns out that a naive spectral estimate is statistically sub-optimal for our LrMM and additional subsequent treatments are necessary. 

For technical simplicity, we adopt the sample splitting in our estimating procedure. It will inevitably affect the constant factor in the error rate, e.g., the $C_1, C_3$ as in Theorem~\ref{thm:mle}. Since our main interest concerns only the convergence rate in terms of the model parameters, we spare no efforts to improve the constant factor. 

Without loss of generality, assume the sample size $n=4n_0$. We randomly split the original sample $\bX_1,\cdots,\bX_n$ into four disjoint subsets of equal size, denoted by $\{\bX_i^{(k)}\}_{i=1}^{n_0}$ for $k=1,2,3,4$. Our estimating procedure consists of three major steps: 
\begin{itemize}[leftmargin=9mm]
\item[-] {\it Step 1 (Spectral initialization)}. Stack the observations column by column into a $d_1\times (n_0d_2)$ matrix $[\bX_1^{(1)},\cdots,\bX_{n_0}^{(1)}]$, extract its leading left singular vector, denoted by $\hat \bu_1$. Then, construct the $d_2\times n_0$ matrix $[\bX_1^{(2)\top}\hat\bu_1,\cdots, \bX_{n_0}^{(2)\top}\hat\bu_1]$ and extract its left singular vector, denoted by $\hat \bv_1$. 

\item[-] {\it Step 2 (Spectral refinement)}. Extract the top-$r$ left and right singular vectors of 
\begin{align}\label{eq:refinement}
\tilde\bU, \tilde \bV\ \ \stackrel{{\rm SVD}_r}{\parbox{1.5cm}{\leftarrowfill}}\ \  \frac{1}{n_0}\sum_{i=1}^{n_0} (\hat \bu_1^{\top}\bX_i^{(3)}\hat\bv_1)\bX_i^{(3)}-\hat \bu_1\hat\bv_1^{\top}
\end{align}
%denoted by $\tilde \bU$ and $\tilde \bV$, respectively. 

\item[-] {\it Step 3 (Aggregation)}. Denote $\check\bM$ the best rank-$r$ approximation of
\begin{align}\label{eq:aggregation}
\check{\bM}\ \ \stackrel{{\rm rank}-r\ {\rm approx.}}{\parbox{2cm}{\leftarrowfill}}\ \ \frac{1}{n_0}\sum_{i=1}^{n_0} {\rm Tr}(\tilde\bU^{\top}\bX_i^{(4)}\tilde\bV)\bX_i^{(4)}-\tilde\bU\tilde\bV^{\top}
\end{align}
Compute the scaling factor by
$$
\hat\Lambda\ \ \longleftarrow \ \ \left[\max\left\{\frac{1}{n_0}\sum_{i=1}^{n_0}{\rm Tr}^2(\tilde\bU^{\top}\bX_i^{(4)}\tilde\bV)-r,\ \frac{dr^2}{\sqrt{n}}\right\}\right]^{1/2}
$$
The final estimator is defined by $\hat\bM=\hat\Lambda^{-1}\check\bM$. 
\end{itemize}
Due to eq. (\ref{eq:aggregation}), we refer to this procedure as the spectral aggregation. Note that (\ref{eq:aggregation}) is, in spirit, similar to the method of second moment as in Gaussian mixture model \citep{wu2020optimal}. The additional projection onto $\tilde\bU$ and $\tilde\bV$ serves the purpose of denoising to leverage the low-rank structure. In this regard, our estimating procedure can also be viewed as a method of projected moments. 
The spectral initialization in Step 1 is very similar to the tensor literature. See, for instance, \cite{montanari2014statistical,zhang2018tensor,xia2019sup} and references therein. A crucial difference here is that the estimate $\hat \bv_1$ relies on the estimate $\hat \bu_1$ to ensure that they are properly correlated in the sense that $\hat \bu_1^{\top}\bM\hat\v_1$ is bounded away from zero, which is a critical requirement for the refinement step (\ref{eq:refinement}). 

Note that the expectation of the RHS of eq. (\ref{eq:aggregation}), with respect to the randomness of $\{\bX_i^{(4)}\}_{i=1}^{n_0}$, is ${\rm Tr}(\tilde \bU^{\top}\bM\tilde\bV)\bM$. Therefore, its best rank-$r$ approximation needs to scaled to serve as a valid estimator for $\bM$. The quantity $\hat\Lambda$ is an estimate of this scaling factor. The performance of the final estimator $\hat\bM$ is guaranteed by the following theorem where we assume $\bM\in\calM(r,d)$ defined in (\ref{eq:calMr}) and $d_1\asymp d_2\asymp d$.  

\begin{theorem}\label{thm:spectral}
There exist absolute constants $c_0, C_0, C_1, C_2, C_3, C_4>0$ such that if the signal strength $\lambda\ge C_0d^{1/2}n^{-1/4}$ and $\min\{d, n\} \ge C_1r\log r$, then with probability at least $1-\exp(-c_0(n\wedge r^{-1}d))$, 
$$
\ell(\hat \bM, \bM)\leq C_2\left(\sqrt{\frac{dr}{n}}+ \frac{1}{\lambda}\sqrt{\frac{d}{n}}\right)
$$
, if we further assume $\lambda\le C_3\exp(c_0(n\wedge r^{-1}d)-\log n)$, then
$$
\E \ell (\hat \bM, \bM)\le C_4\left(\sqrt{\frac{dr}{n}}+ \frac{1}{\lambda}\sqrt{\frac{d}{n}}\right)
$$ 
\end{theorem}

By Theorem~\ref{thm:lower_bound} and Theorem \ref{thm:spectral},  we conclude that the estimator $\hat\bM$ can attain the minimax optimal error rate as long as the signal strength is larger than $d^{1/2}n^{-1/4}$, which we refer to as the {\it strong} signal phase.  This is much more stringent than the information-theoretical limit $d^{1/4}(rn)^{-1/4}+(d/n)^{1/2}$ suggested by the maximum likelihood estimator and minimax lower bound in Section~\ref{sec:MLE}.

%\section{Computationally efficient estimator by projected method of moments}

\section{Statistical and computational tradeoffs}\label{sec:stat-to-comp-gap}
Section~\ref{sec:MLE} and Section~\ref{sec:spectral} indicate the existence of a gap in the signal strengths required by the, {\it in general}, computationally infeasible maximum likelihood estimator and the computationally fast spectral-based estimator. Gap of this type is usually called the statistical-to-computational gap. 
In this section, we provide evidences claiming that no polynomial-time algorithm can consistently estimate $\bM$ if the signal strength is smaller than $d^{1/2}n^{-1/4}$. Our evidence is built on the low-degree likelihood ratio framework for hypothesis testing \citep{kunisky2019notes,loffler2020computationally, hopkins2018statistical}. This framework delivered convincing evidences justifying the statistical-to-computational gap for sparse Gaussian mixture model \citep{loffler2020computationally}  and tensor PCA model, and demonstrated the sharp phase transitions for the spiked Wigner matrix model \citep{kunisky2019notes}. 

The low-degree likelihood ratio framework aims to test two sequences of hypothesis. For our purpose, consider the following hypothesis testing:
\begin{align}\label{eq:test}
H_0^{(n)}: \bM={\bf 0}\quad {\rm versus}\quad H_1^{(n)}: \bM\in \calM(1,\lambda)
\end{align}
where $n$ denotes the sample size. By observing i.i.d. matrices $\bX_1,\cdots,\bX_n$ sampled from the mixture model (\ref{eq:LrMM}), the interest is to test whether the data is pure noise or there is a planted low-rank matrix. Without loss of generality, it suffices to focus on the rank-one case since the ``information" strength $\|\bM\|_{\rm F}$ increases if the rank is larger and, as a result, the hypothesis testing becomes easier for larger ranks. 

Classical textbook results, say, by Neyman-Pearson Lemma, dictate that the likelihood ratio test has preferable power and is uniformly most powerful under some scenarios. Direct computation of the likelihood ratio for testing (\ref{eq:test}) is rather involved due to the composite hypothesis in $H_1^{(n)}$. For simplicity, under the alternative hypothesis, we impose a {\it prior distribution} on $\bM$ assuming that $\bM=\lambda \bu\bv^{\top}$ with a fixed $\lambda$ and the entries of $\bu$ and $\bv$ independently taking the values $\pm d_1^{-1/2}$ and $\pm d_2^{-1/2}$, respectively, with probability $1/2$. Denote $\PP_n$, treated as the alternative hypothesis, the distribution of $(\bX_1,\cdots,\bX_n)$ under LrMM (\ref{eq:LrMM}) with $\bM$ sampled from the aforementioned prior distribution. Note that, for brevity, we suppress the dependence of $\PP_n$ on $\lambda$. 
 Let $\QQ_n$ be the distribution of $(\bX_1,\cdots,\bX_n)$ under the null hypothesis, i.e., each $\bX_i$ is sampled from LrMM (\ref{eq:LrMM}) with $\bM={\bf 0}$. Instead of  (\ref{eq:test}), we consider the following hypothesis testing
\begin{align}\label{eq:testH}
H_0^{(n)}: \bX_1,\cdots,\bX_n\ \stackrel{{\rm i.i.d.}}{\sim}\ \QQ_n\quad {\rm versus}\quad 
H_1^{(n)}: \bX_1,\cdots,\bX_n\ \stackrel{{\rm i.i.d.}}{\sim}\ \PP_n
\end{align}
Denote $L_n(\bcalX):=d\PP_n/d\QQ_n(\bX_1,\cdots,\bX_n)$ the likelihood ratio, where $\bcalX\in\mathbb{R}^{d_1\times d_2\times n}$ is constructed by stacking $n$ data matrices. A well-recognized fact is that the two distributions $\PP_n$ and $\QQ_n$ are {\it statistically indistinguishable} if $\|L_n\|^2:=\EE_{\QQ_n}[L_n(\bcalX)^2]$ remains bounded as $n\to\infty$. Here statistically indistinguishable means that no test can have both type I and type II error probabilities vanishing asymptotically. 

Let $L_n^{\leq D}(\bcalX)$ denote the orthogonal projection of $L_n(\bcalX)$ onto the linear subspace of polynomials $\RR^{d_1\times d_2\times n}\mapsto \RR$ of degree at most $D$. Similarly, define $\|L_n^{\leq D}\|^2:=\EE_{\QQ_n}[L_n^{\leq D}(\bcalX)^2]$. At the core of low-degree likelihood ratio framework is the following conjecture\footnote{We note that a recent work \cite{zadik2021lattice} introduces a very special counter-example to Conjecture~\ref{conjecture}.  However, our LrMM is more closely related to the spiked matrix and tensor model where Conjecture~\ref{conjecture} has contributed convincing evidences to the computational hardness. Therefore, we still postulate the correctness of Conjecture~\ref{conjecture} for our LrMM.},  adapted to the matrix-variate case for our purpose. Here, a test $\phi_n(\cdot)$ taking value $1$ means rejecting the null hypothesis and takes value $0$ if the null hypothesis is not rejected.  
\begin{conjecture}\label{conjecture}
Consider $\mathbb{P}_n$ and $\mathbb{Q}_n$ defined in (\ref{eq:testH}). If there exists $\eps>0$ and $D=D_n\geq (\log nd)^{1+\eps}$ for which $\|L_n^{\le D}\|=1+o(1)$, then there is no polynomial-time test $\phi_n:\mathbb{R}^{d_1\times d_2\times n}\mapsto  \{0,1\}$ such that the sum of type-I error and type-II error probabilities
$$
\E_{\mathbb{Q}_n}[\phi_n(\bcalX)]+\E_{\mathbb{P}_n}[1-\phi_n(\bcalX)]\to 0\quad {\rm as}\quad n\to \infty.
$$
\end{conjecture}

Basically, Conjecture~\ref{conjecture} means that the two distributions $\PP_n$ and $\QQ_n$ are indistinguishable by polynomial-time algorithms if $\|L_n^{\leq D}\|=1+o(1)$. Under the low-degree framework, we now state the computational lower bound of our signal strength for testing \eqref{eq:testH}. 
%This also implies that no-polynomial time algorithms can consistently estimate $\bM$ under LrMM as stated in the following theorem.

\begin{theorem}\label{thm:comp}
Consider $\mathbb{P}_n$ and $\mathbb{Q}_n$ defined in (\ref{eq:testH}). If $\lambda=o(d^{1/2}n^{-1/4})$, then $\|L_n^{\le D}\|^2=1+o(1)$.
%and there exists an absolute constant $c_0>0$ such that, for any polynomial-time estimate $\hat \bM$, the following bound holds
%$$
%\sup_{\bM\in\mathcal{M}(1, \lambda)} \lambda ^{-1}\cdot \EE \ell(\hat\bM, \bM) \ge c_0.
%$$
\end{theorem}

By Theorem~\ref{thm:comp}, conditioned on Conjecture~\ref{conjecture}, detecting the signal matrix in LrMM as in \eqref{eq:testH} becomes computationally hard as long as the signal strength $\lambda$ is at a smaller order of $d^{1/2}n^{-1/4}$. In principle, the estimation of signal matrix is at least as hard (computationally) as detection as in \eqref{eq:testH}, as the latter one only concerns the mere existence thereof, and hence we would expect, at least, the same lower bound also holds for estimation problem in LrMM. %However, it's unclear how to reduce the estimation problem directly to our testing problem as the null hypothesis is inevitably excluded from the parameter set, whose intrinsic difficulty might come from the low-degree testing framework itself. Notably,} 
Notably, if $n=1$, LrMM reduces to the typical matrix perturbation model \citep{cai2018rate, xia2021normal} where there exists no statistical-to-computational gap and the signal strength requirement $O(d^{1/2})$ is both the statistical and computational limit. Interestingly, if $n$ is at the order of $d$, the computational hardness occurs at the signal strength $O(d^{1/4})$ which coincides with the prior literature on spiked tensor model. See \cite{zhang2018tensor, kunisky2019notes} and references therein. %Our Theorem~\ref{thm:comp} significantly complements the aforementioned literature and delivers a 

\section{Numerical simulations}\label{sec:numerical}
\begin{figure}
	\centering
	\includegraphics[width=0.49\linewidth]{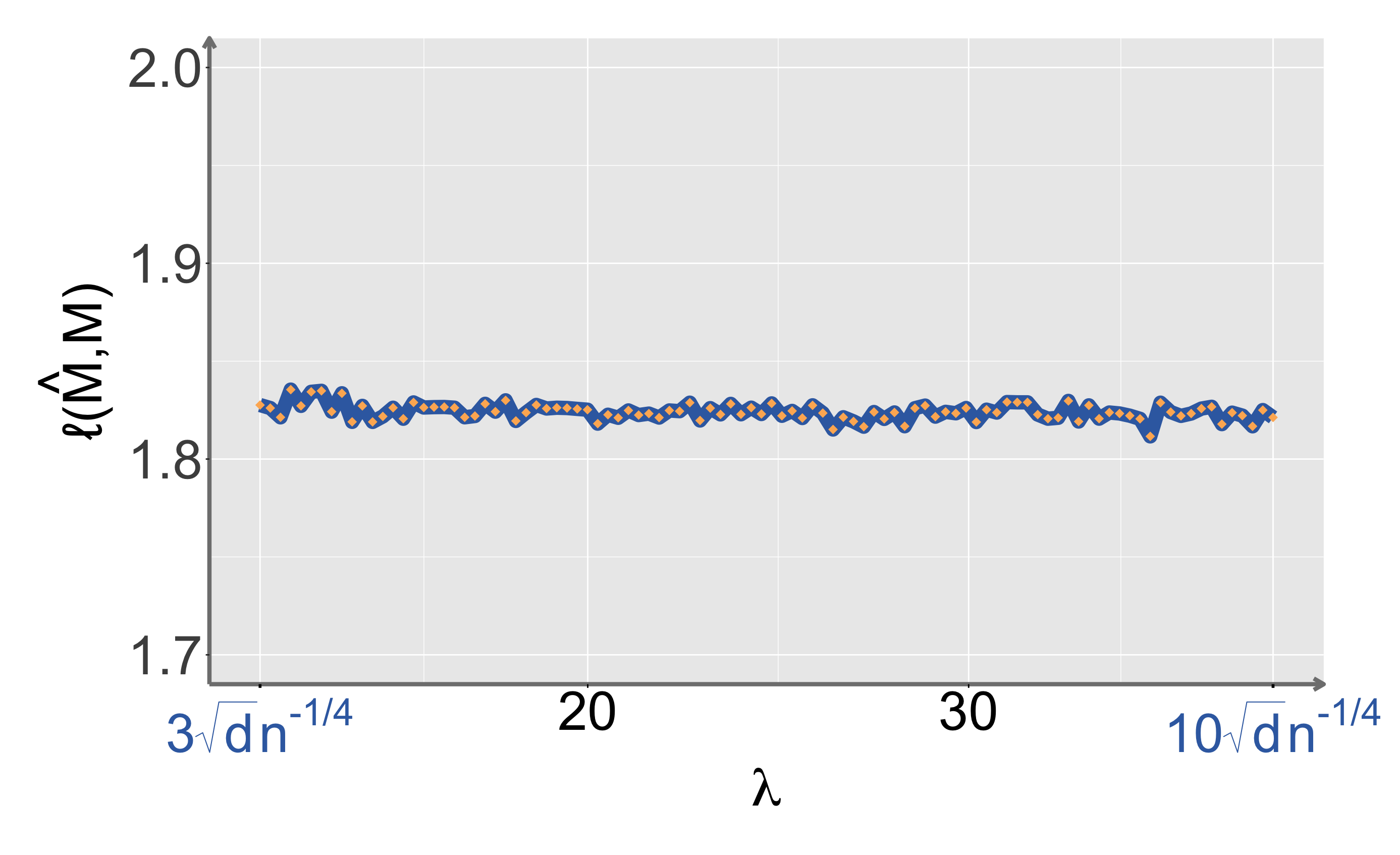}
	\includegraphics[width=0.49\linewidth]{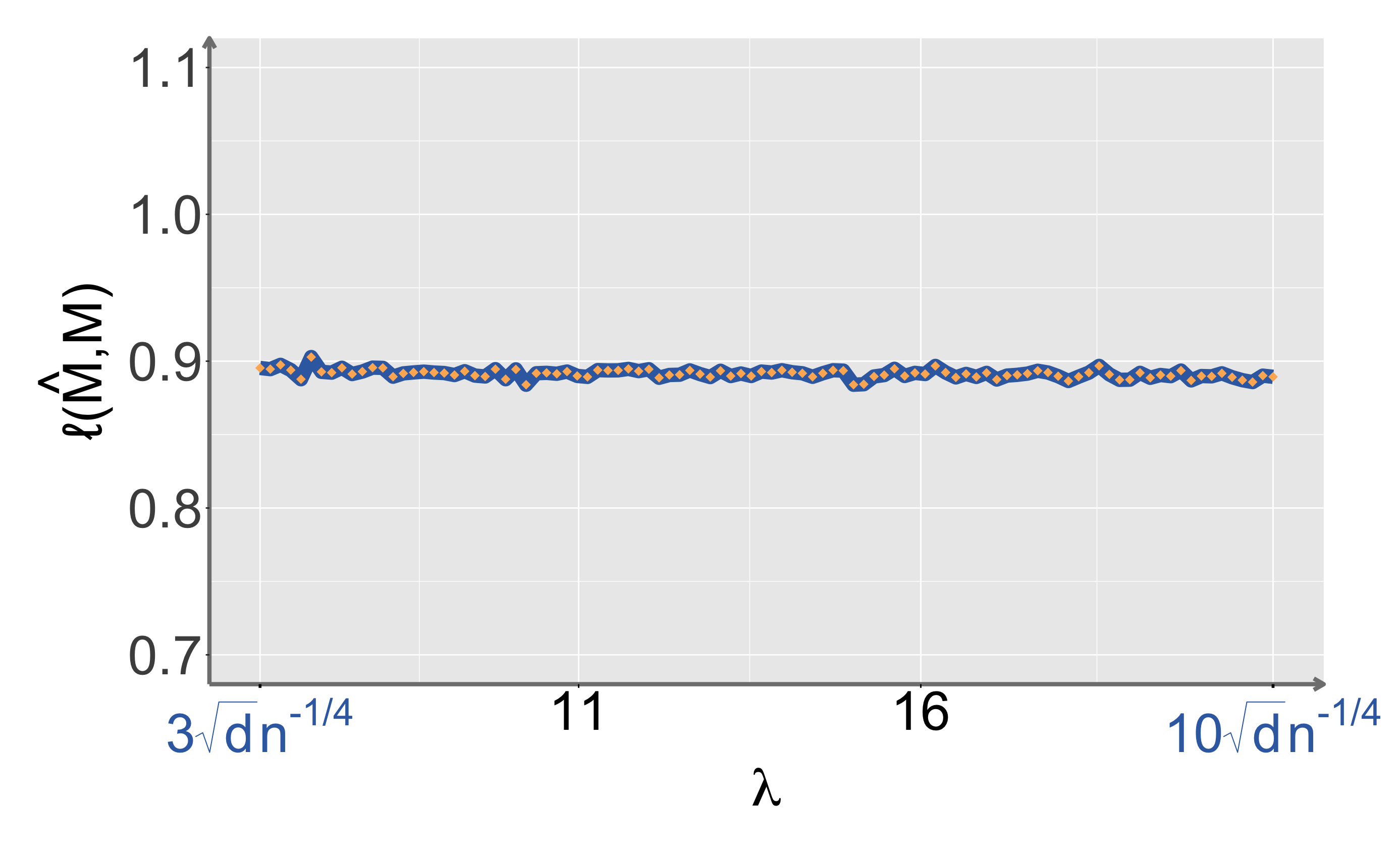}
	\includegraphics[width=0.49\linewidth]{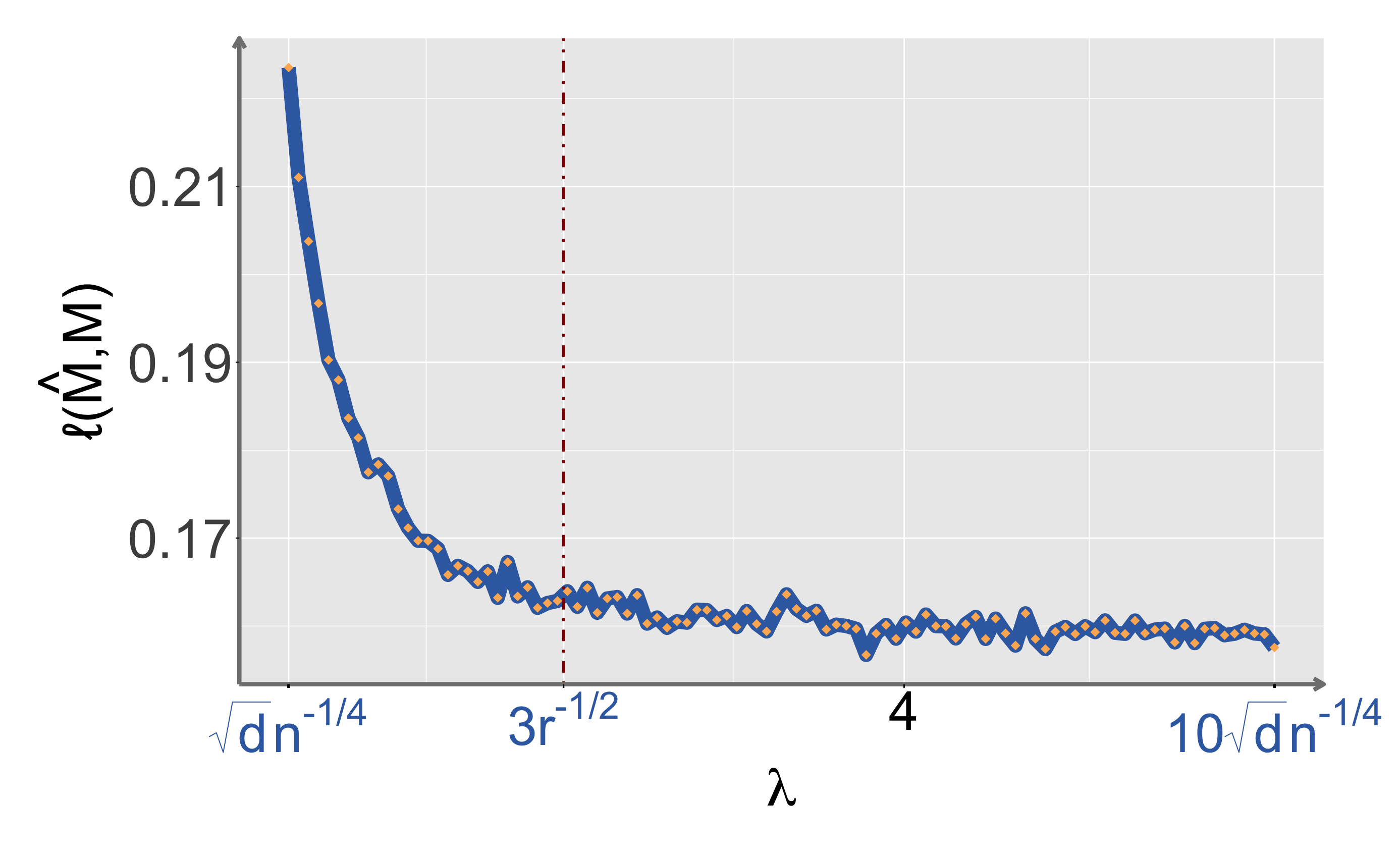}
	\caption{Experiments with $\lambda$ varying. Top-left panel: \textbf{Regime 1}; Top-right panel: \textbf{Regime 2}; Bottom panel: \textbf{Regime 3}.}
	\label{fig:lambda_varying}
\end{figure}

\begin{figure}
	\centering
	\includegraphics[width=0.49\linewidth]{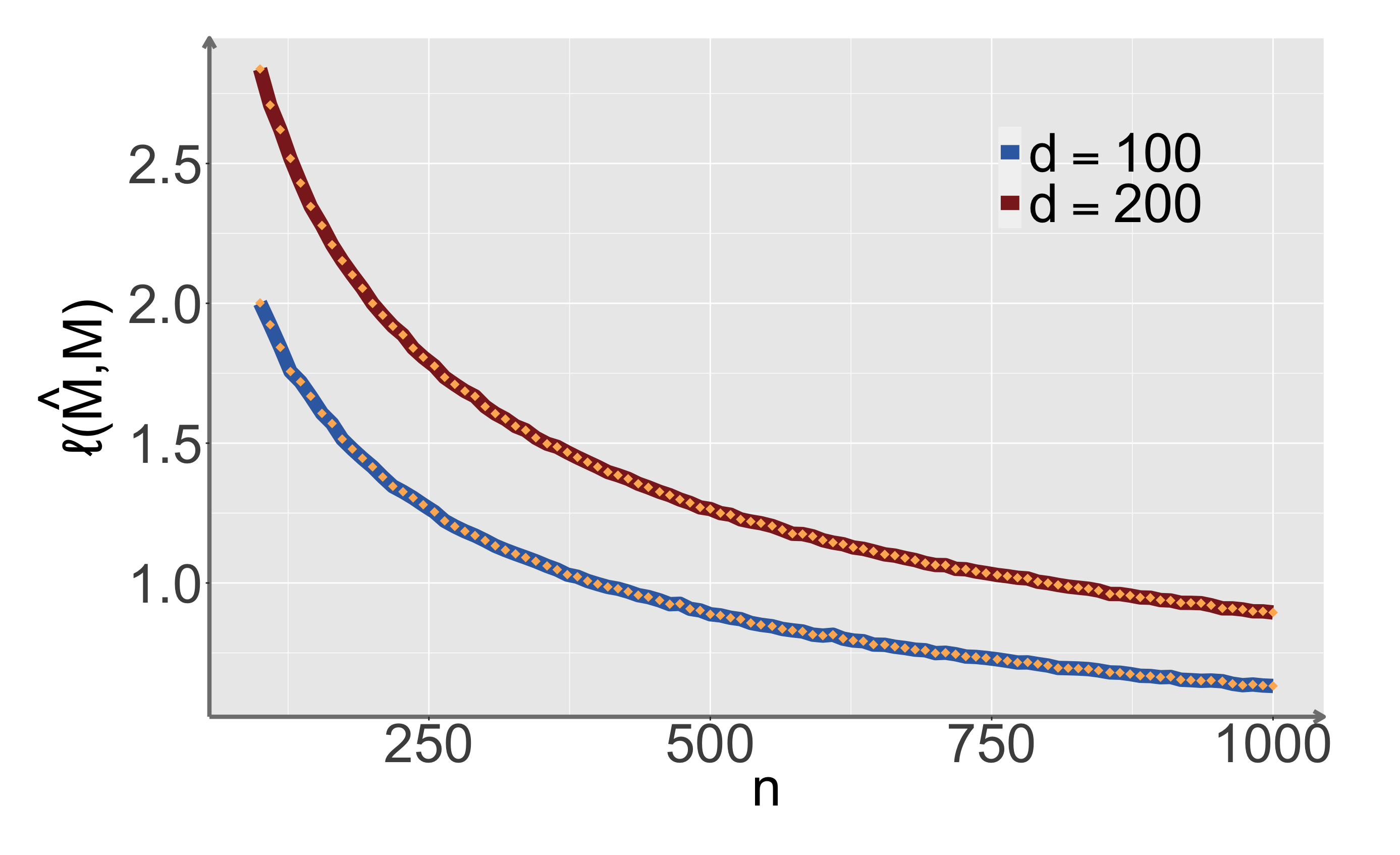}
	\includegraphics[width=0.49\linewidth]{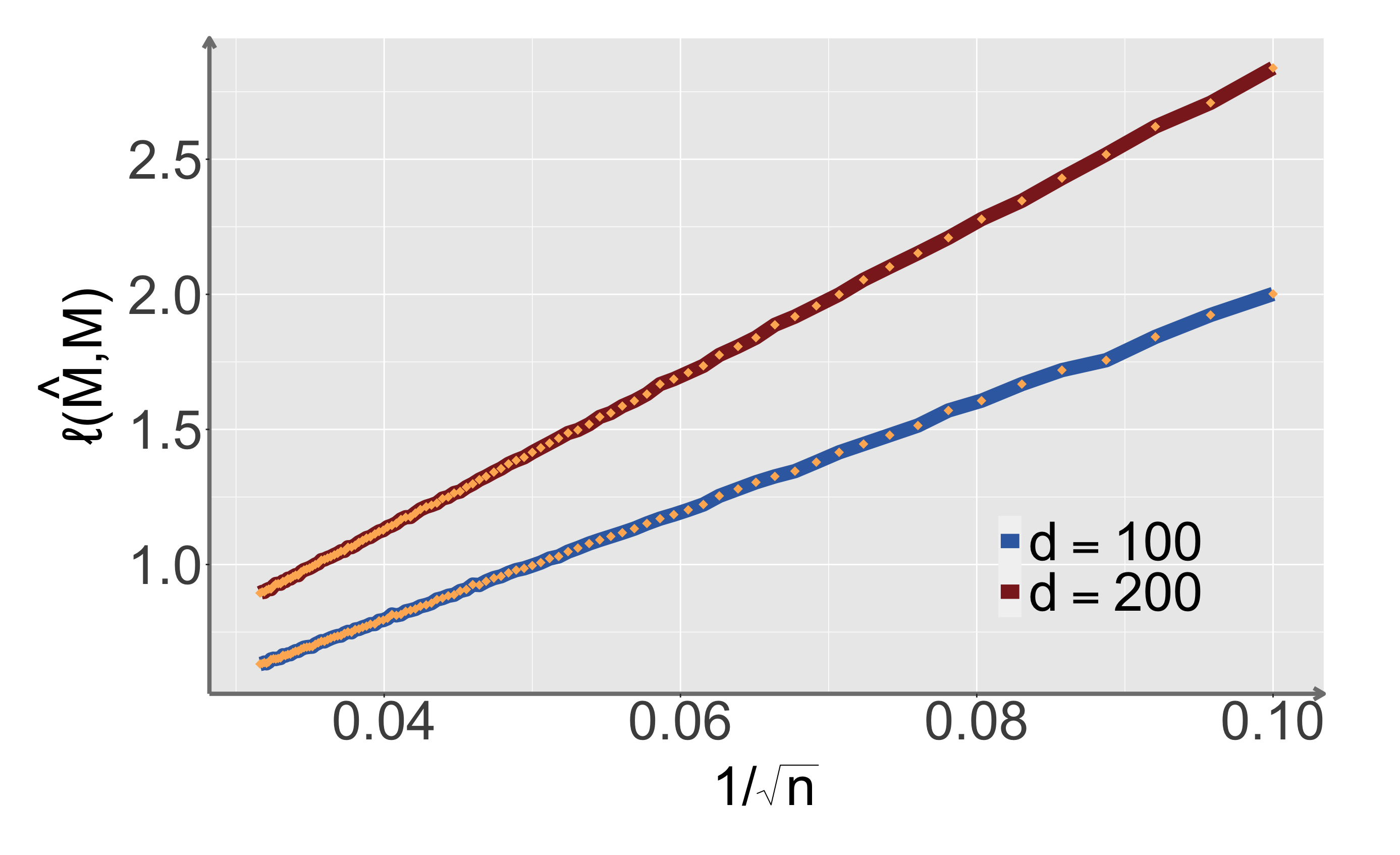}
	\caption{Experiments with $n$ varying. Left panel: $\ell(\hat \bM,\bM)$ against $n$; Right panel: $\ell(\hat \bM,\bM)$ against $n^{-1/2}$. }
	\label{fig:n_varying}
\end{figure}

\begin{figure}
	\centering
	\includegraphics[width=0.5\linewidth]{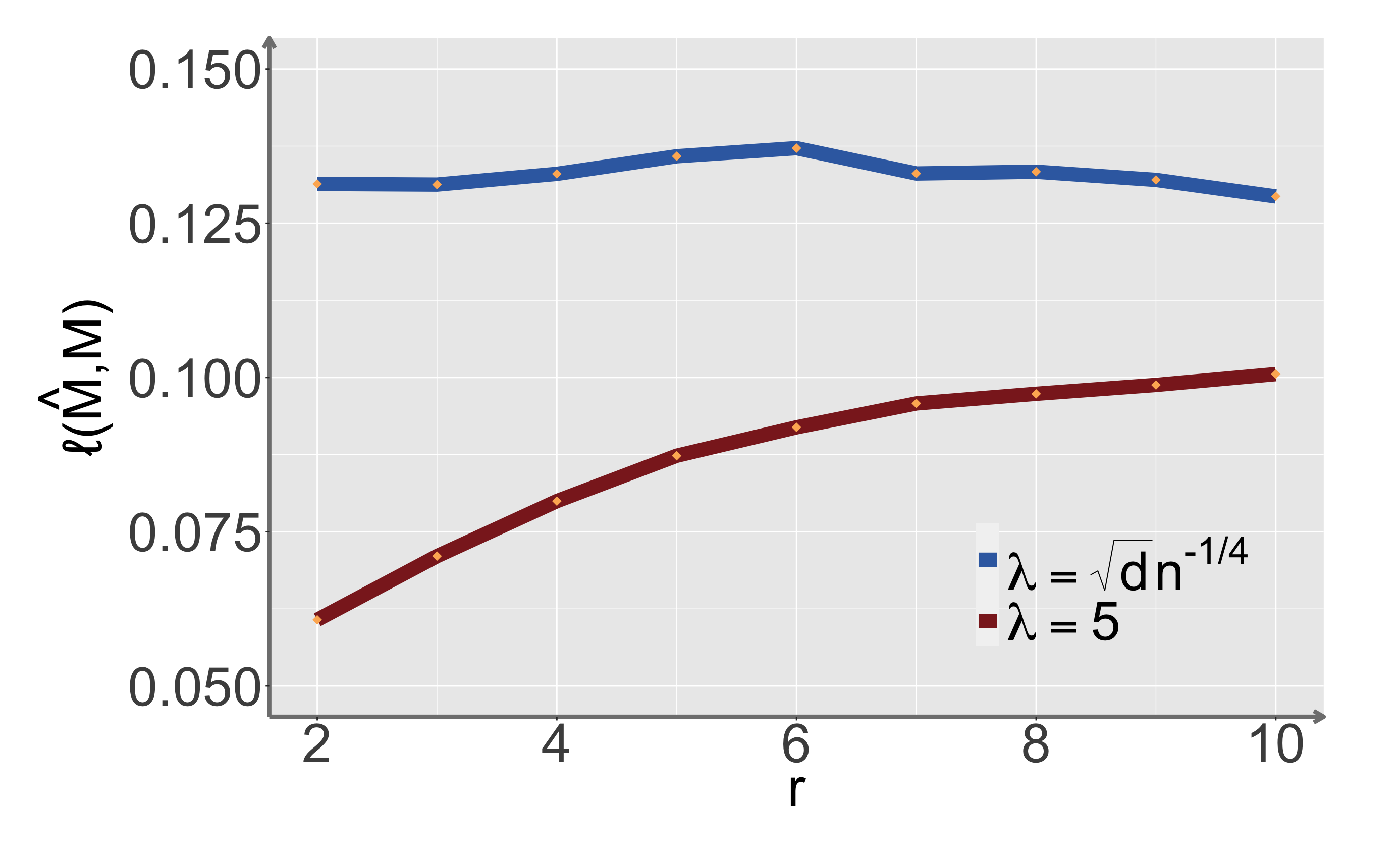}
	\caption{Experiments with $r$ varying. Blue curve $\lambda=d^{1/2}n^{-1/4}$ corresponds to the case where error rate is of order $\lambda^{-1}(d/n)^{1/2}$; Red curve $\lambda=5 (\ge \max_rr^{-1/2})$ corresponds to case where the error rate is of order $(dr/n)^{1/2}$.}
	\label{fig:r_varying}
\end{figure}

In this section, we present numerical experiments to confirm our theoretical findings in the strong signal phase and showcase the performance of our algorithm. Particularly, we apply the spectral aggregate algorithm on $n$ independent data matrices generated from LrMM model in \eqref{eq:LrMM}, with a signal matrix $\bM\in\mathbb{R}^{d\times d}$ of rank $r$ constructed as follows. We first generate two uniformly random $d\times r$ orthonormal matrices $\bU$ and $\bV$, say, by computing the column span (i.e. the image) of a random $d\times r$ Gaussian random matrix with i.i.d. $\mathcal{N}(0, 1)$ entries. Then we fix the smallest and the largest singular value to be $\lambda_r=\lambda$ and $\lambda_1=1.5\lambda$, respectively, and form a diagonal matrix $\bLambda=\text{diag}(\lambda_1,\cdots,\lambda_r)$ where the values of the diagonal terms are equally spaced in a decreasing order. Finally we get our signal matrix $\bM=\bU\bLambda\bV^{\top}$.  We study the effect of parameters $(n,d,r,\lambda)$ on the error $\ell(\hat \bM,\bM)$ via varying one/two parameters while fixing the rest of them.  In each experiment (for a given parameter group $(n,d,r,\lambda)$), the value of the error is the average based on 100 independent simulations with the same signal matrix $\bM$. As the aim of sample splitting step is to facilitate the theoretical analysis, we apply the spectral aggregation algorithm on all samples without sample splitting in all numerical experiments. For brevity, \textbf{Regime 1} is referred to the case when $\boldsymbol{n\leq dr}$, \textbf{Regime 2} is referred to the case when $\boldsymbol{dr\leq n\leq (dr)^2}$ and \textbf{Regime 3} is referred to the case when $\boldsymbol{n\geq (dr)^2}$.  The information are summarized as follows:
\begin{itemize}
	\item \textbf{Experiment 1}: $n=300$, $d=250$, $r=2$ (\textbf{Regime 1}). $\lambda$ is varying from $3\sqrt{d}n^{-1/4}$ to $10\sqrt{d}n^{-1/4}$.
	\item \textbf{Experiment 2}: $n=500$, $d=100$, $r=2$ (\textbf{Regime 2}). $\lambda$ is varying from $3\sqrt{d}n^{-1/4}$ to $10\sqrt{d}n^{-1/4}$.
	\item \textbf{Experiment 3}: $n=3000$, $d=20$, $r=2$ (\textbf{Regime 3}). $\lambda$ is varying from $3\sqrt{d}n^{-1/4}$ to $10\sqrt{d}n^{-1/4}$.
	\item \textbf{Experiment 4}: $d\in\{100,200\}$, $r=2$. $n$ is varying from $100$ to $1000$ with $\lambda=3\sqrt{d}n^{-1/4}$.
	\item \textbf{Experiment 5}: $n\in\{100,200\}$, $r=2$. $d$ is varying from $100$ to $500$ with $\lambda=3\sqrt{d}n^{-1/4}$.
	\item \textbf{Experiment 6}: $n=10000$, $d=10$, $\lambda\in\{\sqrt{d}n^{-1/4},5\}$ (\textbf{Regime 3}). $r$ is varying from $2$ to $10$.
\end{itemize}
In \textbf{Experiment 1} \& \textbf{2} (\textbf{Regime 1} \& \textbf{2}), the error stays almost constant as $\lambda$ increases. Both cases fall into the \textit{strong} signal phase and an optimal rate of $O((dr/n)^{1/2})$ can be attained, suggested by Theorem \ref{thm:spectral}. While in \textbf{Experiment 3} (\textbf{Regime 3}) with the same range of $\lambda$, the phase transition effect is clearly demonstrated in the bottom panel of Figure \ref{fig:lambda_varying}: when $\lambda$ varies from $C_1d^{1/2}n^{-1/4}$ to $C_2r^{-1/2}$, the optimal rate $O(\lambda^{-1}(d/n)^{1/2})$ is linear in $\lambda^{-1}$; when $\lambda\ge C_2r^{-1/2}$, the optimal rate $O((dr/n)^{1/2})$ is again independent of $\lambda$.\\
In \textbf{Experiment 4} \& \textbf{5}, we screen the effect of varying $n$ and $d$, respectively in Figure \ref{fig:n_varying}. As expected, the error becomes smaller as $n$ grows (or $d$ decreases). The linearity between the error rate and $n^{-1/2}$ (or $d^{1/2}$) can be verified in the right panels, which is in accordance with Theorem \ref{thm:spectral}.\\
In \textbf{Experiment 6}, we let $r$ vary with other parameters fixed and focus on \textbf{Regime 3}, which is the most interesting case due to the phase transition effect in terms of rank $r$. As shown in Figure \ref{fig:r_varying}, the error rate $O(\lambda^{-1}(d/n)^{1/2})$ is constant in $r$ with $\lambda\in (C_1d^{1/2}n^{-1/4},C_2r^{-1/2})$ and when $\lambda\ge C_2r^{-1/2}$, the error rate increases with $r$.
\begin{figure}
	\centering
	\includegraphics[width=0.49\linewidth]{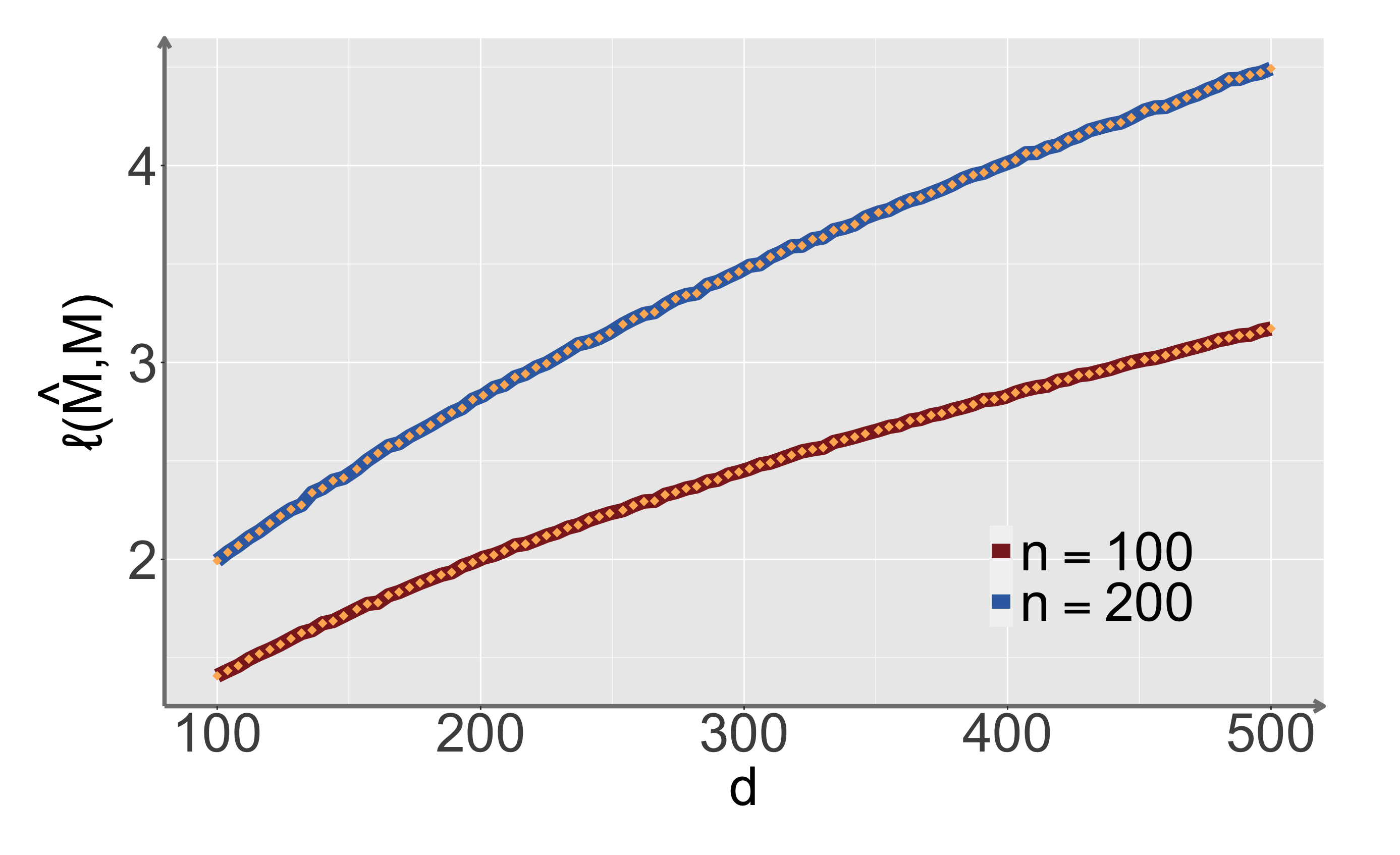}
	\includegraphics[width=0.49\linewidth]{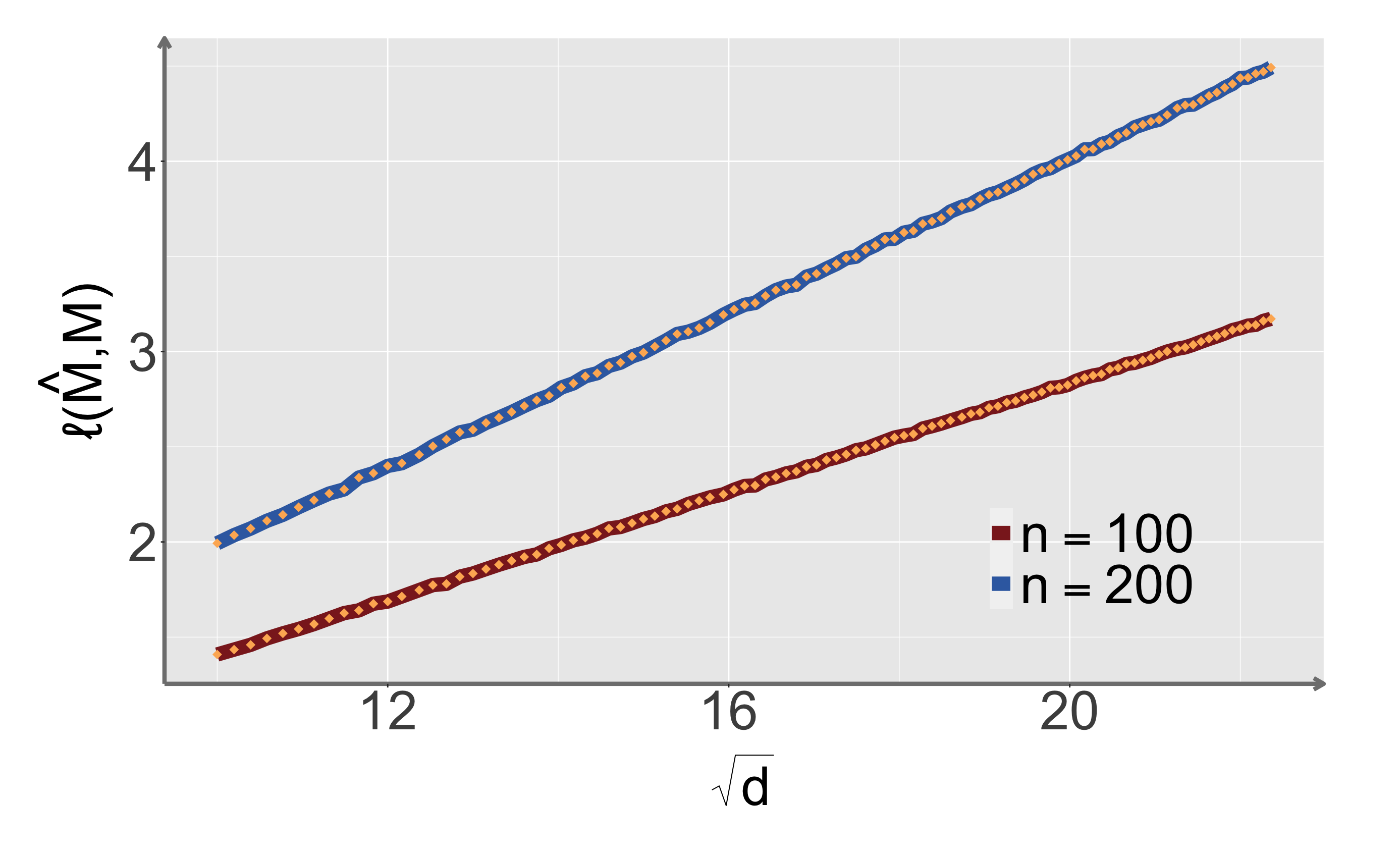}
	\caption{Experiments with $d$ varying. Left panel: $\ell(\hat \bM,\bM)$ against $d$; Right panel: $\ell(\hat \bM,\bM)$ against $d^{1/2}$. }
	\label{fig:d_varying}
\end{figure}

\section{Real data experiment}\label{sec:realdata}
\begin{figure}
	\centering
	\includegraphics[width=0.49\linewidth]{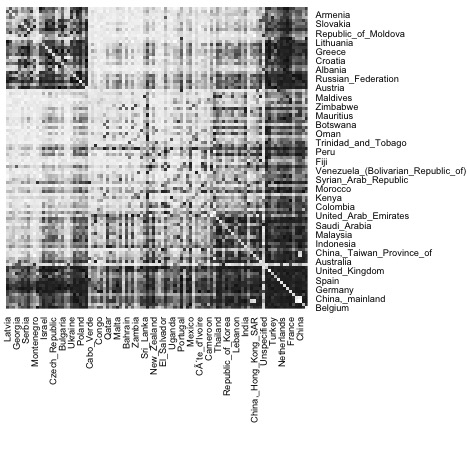}
	\includegraphics[width=0.49\linewidth]{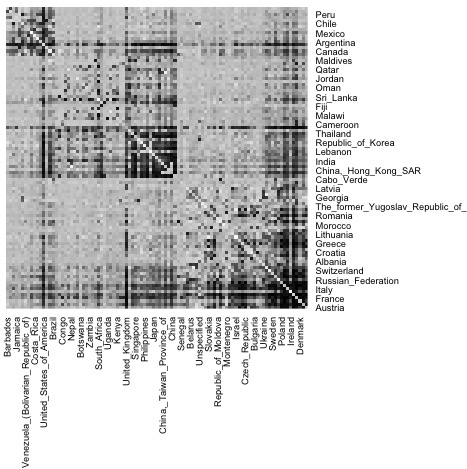}
	\caption{Heatmaps for $\hat \bM_1$ and $\hat \bM_2$.}
	\label{fig:M_1M_2}
\end{figure}

We present an application of our algorithm on a real-world dataset, which is a collection of multiple layers of worldwide food trading networks (\cite{de2015structural}), recording the trade flows of 30 food products between 99 countries. We pre-process the data the same as in \cite{jing2021community} and end up with a 3-rd order binary tensor $\bcalX$ of dimension $99 \times  99 \times  30$. Each layer of this tensor $[\bcalX]_{\cdot\cdot i}=\bX_i$ represents the adjacency matrix of one specific type of food product $i$, and nodes are different countries/regions which are common across all layers. As shown in \cite{jing2021community}, the layers could be clustered into two groups, one of which mainly consists of raw or unprocessed food and another is made of processed food. We adopt this clustering result as ground truth and assume all layers are generated independently according to two expected adjacency matrices $\bM_1,\bM_2\in\mathbb{R}^{99\times 99}$. Note that though throughout the paper the noise matrix $\bZ_i$ is assumed to be Gaussian , we believe the spectral aggregation can be applied to more general setting (for instance, observations with sub-gaussian noise). Our goal is to recover $\bM_1$ and $\bM_2$.  To make it adapted to our framework (as mentioned in Section \ref{sec:intro}), we first construct centered observations $\widetilde\bX_i=\bX_i-\bar\bX$, where $\bar\bX=n^{-1}\sum_{i=1}^n\bX_i$ is the sample average of adjacency matrices over all layers. Here, $\bar\bX$ serves as an estimate of $(\bM_1+\bM_2)/2$. Then we apply the spectral aggregation algorithm with rank $r=10$ to $\{\widetilde \bX_i\}_{i=1}^n$ to get $\hat {\bM}$. It turns out that the final result is not sensitive to choice of rank $r$. Finally we can construct $\hat\bM_1=\bar\bX+\hat {\bM}$ and $\hat\bM_2=\bar\bX-\hat {\bM}$. To appropriately visualize our result, we rearrange the order of columns and rows of $\hat \bM_1$ and $\hat \bM_2$ in the same way as in \cite{jing2021community}, which is based on the community labels estimated by tensor method therein, in order to have a glance of community structures. In Figure \ref{fig:M_1M_2}, the mean matrix in the left panel demonstrates a strong trend of global trading, while the other one shows the dominance of regional trading. These findings coincides with results in  \cite{jing2021community}, whereas we are estimating the difference of two center matrices instead of clustering all observations. Note that the results in \cite{jing2021community} require layer clustering before producing $\hat \bM_1, \hat\bM_2$ but our method does not.

\section{Discussion}\label{sec:discuss}
Our main focus in this paper is on the optimal estimation and computational limits for the two-component low-rank Gaussian mixtures. It is of great interest to investigate the minimax optimal estimation when the number of components is greater than two. Unfortunately, our spectral aggregation method is inapplicable and we cannot immediately see an easy generalization of the maximum likelihood estimator to the multi-component case. There are several possibilities. For instance, unlike the two-component case, it might be necessary to, at least partially, recover the latent labels before estimating the underlying low-rank components.  Indeed, the linear regression low-rank mixture model \citep{chen2021learning} was treated by this way. However, it is well recognized that consistent clustering often requires a much stronger condition on the signal strength. See, for instance, \cite{loffler2019optimality, wu2019randomly} and references therein. For the two-component symmetric case as in model (\ref{eq:LrMM}), consistent clustering requires a signal strength at least\footnote{To see this, one can simply assume the singular vectors $\bU$ and $\bV$ are available before hand.} in the order of $\Omega(1)$ when $r$ is a constant, which can be much more stringent than the condition required by the spectral aggregation method in Regime 3. It therefore indicates another possibility: there might exist some method that can reliably estimate the multiple low-rank components without the prerequisite of meaningful clustering. We leave this for future works. 

\newpage

\bibliographystyle{plainnat}
\bibliography{references} % see references.bib for bibliography management

\newpage

\appendix

\section{Proofs for main results}
\subsection{Proof of Theorem \ref{thm:mle}}
For technical reasons discussed in Section \ref{sec:MLE}, we split our proof into two cases, corresponding to the first and second statement in Theorem \ref{thm:mle}.
\subsubsection*{Case 1: $dr\log (nd)<n$}
In this regime, the standard tool to establish the convergence rate of MLE is applicable. To this end, we need to introduce the following notations. Define 
$$\bar{\mathcal{P}}_{d_1,d_2}(r,\lambda):=\left\{\frac{p_{\bM}+p_{\bM^\prime}}{2}:\bM^\prime \in \mathcal{M}_{d_1,d_2}(r,\lambda)\right\},\quad \bar{\mathcal{P}}^{1/2}_{d_1,d_2}(r,\lambda):=\left\{p^{\frac{1}{2}}:p\in \bar{\mathcal{P}}_{d_1,d_2}(r,\lambda)\right\}$$
and for any  small $\delta>0$, define a Hellinger ball centered at $p_{\bM}$ with radius $\delta$ by 
$$\bar{\mathcal{P}}^{1/2}_{d_1,d_2}(r,\lambda,\delta):=\left\{\bar p^{\frac{1}{2}}\in  \bar{\mathcal{P}}^{1/2}_d(\lambda): d_{\textsf{H}}\left(\bar p,p_{\bM}\right )\le \delta\right\}$$
We refer to $H_B(\epsilon, \bar{\mathcal{P}}^{1/2}_{d_1,d_2}(r,\lambda,\delta),L_2(\mu))$ as the $\epsilon$-bracketing entropy of $\bar{\mathcal{P}}^{1/2}_{d_1,d_2}(r,\lambda,\delta)$ under $L_2(\mu)$ metric with Lebesgue measure $\mu$ and view $\mathcal{J}_B(\delta,\bar{\mathcal{P}}^{1/2}_{d_1,d_2}(r,\lambda,\delta),L_2(\mu ))$ as the entropy integral of $\bar{\mathcal{P}}^{1/2}_{d_1,d_2}(r,\lambda,\delta)$, which is defined as
$$\mathcal{J}_B(\delta,\bar{\mathcal{P}}^{1/2}_{d_1,d_2}(r,\lambda,\delta),L_2(\mu )):=\int _{\delta^2/2^{13} }^\delta H_B^{1/2}(\epsilon, \bar{\mathcal{P}}^{1/2}_{d_1,d_2}(r,\lambda,\delta),L_2(\mu))d\epsilon \vee \delta $$
Now we state Theorem 7.4 in \cite{geer2000empirical} (adapted to our notation), which establishes the rate of convergence of MLE. 
\begin{lemma}[\cite{geer2000empirical}]\label{lem:mle_conv}
Take $\Psi(\delta)\ge \mathcal{J}_B(\delta,\bar{\mathcal{P}}^{1/2}_{d_1,d_2}(r,\lambda,\delta),L_2(\mu ))$ in such a way that $\Psi(\delta)/\delta^2$ is a non-increasing function of $\delta$. Then for a universal constant $c$, and for 
$$\sqrt n\delta_n^2\ge c\Psi(\delta_n)$$
we have for all $\delta\ge \delta_n$
$$\Prob(d_{\textsf{H}}(p_{\hat \bM_{\MLE}}, p_{\bM})>\delta)\le c\exp\left(-\frac{n\delta^2}{c^2}\right)$$
\end{lemma}
\noindent A combination of Lemma \ref{lem:hellinger} and Lemma \ref{lem:mle_conv} implies that the convergence rate of $\hat \bM_{\MLE}$ would entail an upper bound on the $\epsilon$-bracketing entropy $H_B(\epsilon, \bar{\mathcal{P}}^{1/2}_{d_1,d_2}(r,\lambda,\delta),L_2(\mu))$. Notice that for any $\delta>0$,
\begin{align}\label{lem-rel-bracket}
H_B(\epsilon, \bar{\mathcal{P}}^{1/2}_{d_1,d_2}(r,\lambda,\delta),L_2(\mu))&\overset{(a)}{\le} H_B(\epsilon, \bar{\mathcal{P}}^{1/2}_{d_1,d_2}(r,\lambda),L_2(\mu))\overset{(b)}{=}H_B({\epsilon}/{\sqrt{2}}, \bar{\mathcal{P}}_{d_1,d_2}(r,\lambda),d_{\textsf{H}})\nonumber \\ 
&\overset{(c)}{\le} H_B({\epsilon}, {\mathcal{P}}_{d_1,d_2}(r,\lambda),d_{\textsf{H}})
\end{align}
where (a) is due to $\bar{\mathcal{P}}^{1/2}_{d_1,d_2}(r,\lambda,\delta)\subset\bar{\mathcal{P}}^{1/2}_{d_1,d_2}(r,\lambda)$, (b) follows from the definition of Hellinger distance $d_{\textsf{H}}$ and (c) is due to the following fact (cf. Lemma 4.2 in \cite{geer2000empirical}, \cite{ho2016convergence}): for any $\bar p_1=\frac{1}{2}(p_{\bM_1}+p_{\bM})\in \bar{\mathcal{P}}_{d_1,d_2}(r,\lambda),\bar p_2=\frac{1}{2}(p_{\bM_2}+p_{\bM})\in \bar{\mathcal{P}}_{d_1,d_2}(r,\lambda)$
$$d_{\textsf{H}}^2(\bar p_1,\bar p_2)\le \frac{1}{2}d_{\textsf{H}}^2(p_{\bM_1},p_{\bM_2})$$
In view of \eqref{lem-rel-bracket}, it suffices to bound $H_B(\epsilon,\mathcal{P}_{d_1,d_2}(r,\lambda),d_{\textsf{H}})$. The following lemma characterize the size of bracketing entropy of $\mathcal{P}_{d_1,d_2}(r,\lambda)$.
\begin{lemma}\label{lem:bracketing-entropy} Assume $d_1\asymp d_2\asymp d$ then we have
	$$H_B(\epsilon,\mathcal{P}_{d_1,d_2}(r,\lambda),d_{\textsf{H}})\lesssim dr\log\left(\frac{d}{\epsilon}\right)$$
\end{lemma}
\noindent Using relation \eqref{lem-rel-bracket} and Lemma \ref{lem:bracketing-entropy} we can arrive at
$$\mathcal{J}_B(\delta,\bar{\mathcal{P}}^{1/2}_{d_1,d_2}(r,\lambda,\delta),L_2(\mu ))\lesssim \int _{\delta^2/2^{13} }^\delta \sqrt{dr\log\left(\frac{d}{\epsilon}\right)}d\epsilon \vee \delta \lesssim \delta \sqrt{dr\log \left(\frac{d}{\delta}\right)}$$
Now we can take $\Psi(\delta)=C\delta \sqrt{dr\log \left(\frac{d}{\delta}\right)}$ for some absolute constant $C>0$ and $\delta=\delta_n=\sqrt{\frac{dr}{n}\log (nd)}$, then we have $\Psi(\delta)/\delta^2=C\frac{1}{\delta} \sqrt{dr\log \left(\frac{d}{\delta}\right)}$ is a non-increasing function of $\delta $ and that 
$$\sqrt{n}\delta_n^2=\frac{dr}{\sqrt{n}}\log(nd)\ge c\frac{dr}{\sqrt n}\sqrt{\log (nd)} \sqrt{\log\left(\frac{\sqrt{nd}}{\sqrt{r\log (nd)}}\right)}=c\Psi(\delta_n)$$
By Lemma \ref{lem:mle_conv}, with probability at least $1-\exp(-cd\log^2(nd))$ we have
$$d_{\textsf{H}}(p_{\hat \bM_\MLE}, p_{\bM})\le C\sqrt{\frac{dr\log (nd)}{n}}$$
It suffices to use Lemma \ref{lem:hellinger} to connect the density estimation and parameter estimation. Notice that $\fro{\hat \bM_\MLE}+\fro{\bM}\asymp \lambda\sqrt{r}$. By Lemma \ref{lem:hellinger}, if $\lambda\sqrt{r}\lesssim 1$, with probability at least $1-\exp(-cd\log^2(nd))$:
$$\ell(\hat \bM_{\MLE},\bM) \lesssim (\lambda\sqrt{r})^{-1}\cdot d_{\textsf{H}}(p_{\hat \bM_\MLE}, p_{\bM})\le \frac{1}{\lambda}\sqrt{\frac{d\log (nd)}{n}}$$
If $\lambda\sqrt{r}\gtrsim 1$, note that in this case ($dr<n\log(nd)$), we have with probability at least $1-\exp(-cd\log^2(nd))$:
$$\min\{1,\ell(\hat \bM_\MLE,\bM)\} \lesssim  d_{\textsf{H}}(p_{\hat \bM_{\MLE}}, p_{\bM})\le \sqrt{\frac{dr\log (nd)}{n}}<1 $$
implying that  $\ell(\hat \bM_{\MLE},\bM) \lesssim \sqrt{{dr\log (nd)}/{n}}$. 
Combining two pieces we conclude that with probability at least $1-\exp(-cd\log^2(nd))$: 
$$\ell(\hat \bM_{\MLE},\bM)\le C\left(\sqrt{\frac{dr\log (nd)}{n}}\vee \frac{1}{\lambda}\sqrt{\frac{d\log (nd)}{n}}\right)$$
We can further have a bound in expectation:
$$\E \ell(\hat \bM_{\MLE},\bM)\le C\left(\sqrt{\frac{dr\log (nd)}{n}}\vee \frac{1}{\lambda}\sqrt{\frac{d\log (nd)}{n}}\right)$$
provided that $\lambda \le \exp(cd\log^2(nd))$.
\subsubsection*{Case 2: $dr\log(nd)\ge n$}
In this regime, our ultimate goal is to have $\fro{\hat \bM_{\MLE}-\bM}\lesssim \sqrt{dr\log(nd)/n}$ with high probability and in expectation and hence we can assume  $\fro{\hat \bM_{\MLE}-\bM}\ge c_0\sqrt{dr\log(nd)/n}$ for some absolute constant $c_0>0$ (otherwise we have the desired result). Without loss of generality, we assume $\fro{\hat\bM_{\MLE}-\bM}\le \fro{\hat\bM_{\MLE}+\bM}$. Unlike Case 1, we resort to KL divergence instead of Hellinger distance to establish the convergence rate. Let $P_\bM$ denote the distribution of \eqref{eq:LrMM} and recall the definition of KL divergence, for any $\bM,\bM^\prime\in\mathcal{M}_{d_1,d_2}(r,\lambda)$ we have 
\begin{align*}
	D_{\textsf{KL}}\left(p_{\bM} \big \| p_{ \bM^\prime}\right)=\int \left(\log \frac{p_{\bM}}{p_{\bM^\prime }}\right)dP_\bM
\end{align*}
Note that for fixed $\bM$ and $\bM^\prime$, we simply have $D_{\textsf{KL}}\left(p_{\bM} \big \| p_{ \bM^\prime}\right)=\E\log \left({p_{\bM}}(\bX)/{p_{\bM^\prime }(\bX)}\right)$ for $\bX\sim p_\bM$. On the other hand, by the definition of the maximum likelihood estimator $\hat \bM_{\MLE}$, we have
$$\frac{1}{n}\sum_{i=1}^n\log \frac{p_{\bM}(\bX_i)}{p_{\hat \bM_{\MLE}}(\bX_i)}\le 0$$
Therefore, we can have that
\begin{align}\label{ineq:mle_KL}
D_{\textsf{KL}}\left(p_{\bM} \big \| p_{ \hat \bM_{\MLE}}\right)\le -\frac{1}{n}\sum_{i=1}^n\log \frac{p_{\bM}(\bX_i)}{p_{\hat \bM_{\MLE}}(\bX_i)}+D_{\textsf{KL}}\left(p_{\bM} \big \| p_{ \hat \bM_{\MLE}}\right)	
\end{align}
Now we give an upper bound of RHS of \eqref{ineq:mle_KL}. To this end, we consider a ball in $\mathcal{M}_{d_1,d_2}(r,\lambda)$ with radius $\delta$, i.e., $\mathcal{M}(\delta):=\{\bM^\prime\in \mathcal{M}_{d_1,d_2}(r,\lambda): \fro{\bM^\prime-\bM}\le \delta\}$. Our aim is to bound the following quantity:
$$\theta_n(\delta):=\sup_{\bM^\prime\in \mathcal{M}(\delta)}\left|\frac{1}{n}\sum_{i=1}^n\log \frac{p_{\bM}(\bX_i)}{p_{\bM^\prime}(\bX_i)}-D_{\textsf{KL}}\left(p_{\bM} \big \| p_{\bM^\prime}\right)\right|	$$
Observe that
\begin{align*}
\log \frac{p_{\bM}(\bX)}{p_{\bM^\prime}(\bX)}&=\log \left(\frac{e^{-\frac{1}{2}\fro{\bX-\bM}^2}+e^{-\frac{1}{2}\fro{\bX+\bM}^2}}{e^{-\frac{1}{2}\fro{\bX-\bM^\prime }^2}+e^{-\frac{1}{2}\fro{\bX+\bM^\prime}^2}}\right)=\frac{1}{2}\fro{\bM^\prime}^2-\frac{1}{2}\fro{\bM}^2+\log \left(\frac{e^{\langle \bX,\bM \rangle}+e^{-\langle \bX,\bM \rangle}}{e^{\langle \bX,\bM^\prime \rangle}+e^{-\langle \bX,\bM^\prime \rangle}}\right)
\end{align*}
By log-sum-exp inequality, we have
\begin{align*}
|\langle \bX,\bM\rangle|-|\langle \bX,\bM^\prime\rangle|-\log 2\le \log \left(\frac{e^{\langle \bX,\bM \rangle}+e^{-\langle \bX,\bM \rangle}}{e^{\langle \bX,\bM^\prime \rangle}+e^{-\langle \bX,\bM^\prime \rangle}}\right)\le \log2+|\langle \bX,\bM\rangle|-|\langle \bX,\bM^\prime\rangle|
\end{align*}
Hence we have
\begin{align*}
\frac{1}{n}\sum_{i=1}^n\log \frac{p_{\bM}(\bX_i)}{p_{\bM^\prime}(\bX_i)}-D_{\textsf{KL}}\left(p_{\bM} \big \| p_{\bM^\prime}\right)\le 2\log2+\frac{1}{n}\sum_{i=1}^n\left[|\langle \bX_i,\bM\rangle|-|\langle \bX_i,\bM^\prime\rangle|\right]-\E\left[|\langle \bX,\bM\rangle|-|\langle \bX,\bM^\prime\rangle|\right]
\end{align*}
which implies $\theta_n(\delta)\le 2\log2+\tilde \theta_n(\delta)$, where
\begin{align*}
\tilde \theta_n(\delta):=\sup_{\bM^\prime\in \mathcal{M}(\delta)}\left|\frac{1}{n}\sum_{i=1}^n\left[|\langle \bX_i,\bM\rangle|-|\langle \bX_i,\bM^\prime\rangle|\right]-\E\left[|\langle \bX,\bM\rangle|-|\langle \bX,\bM^\prime\rangle|\right]\right|
\end{align*}
To get a high probability bound for $\tilde \theta_n(\delta)$, we first upper bound its expectation. By symmetrization (see, e.g., in \cite[Lemma 2.3.1]{van1996weak}), we have 
\begin{align*}
	\E \tilde \theta_n(\delta)\le 2 \E\left(\sup_{\bM^\prime\in \mathcal{M}(\delta)}\left|\frac{1}{n}\sum_{i=1}^n\varepsilon_i\left[|\langle \bX_i,\bM\rangle|-|\langle \bX_i,\bM^\prime\rangle|\right]\right|\right)
\end{align*}
where $\{\varepsilon_i\}_{i=1}^n$ are independent Rademacher random variables, which is independent of $\{\bX_i\}_{i=1}^n$. Denote $\phi_i(\bM^\prime)=|\langle \bX_i,\bM\rangle|-|\langle \bX_i,\bM^\prime\rangle|$, for any $\bM_1,\bM_2\in \mathcal{M}(\delta)$ we have 
$$|\phi_i(\bM_1)-\phi_i(\bM_2)|\le |\langle \bX_i, \bM_1-\bM_2\rangle|=|\langle \bX_i, \bM_1-\bM\rangle-\langle \bX_i, \bM_2-\bM\rangle|$$
which means $\phi_i(\bM^\prime)$ is $1$-Lipschitz in $\langle \bX_i,\bM^\prime-\bM\rangle$. By comparison theorem (\cite[Theorem 4.12]{ledoux1991probability}), we deduce that 
\begin{align*}
	\E\left(\sup_{\bM^\prime\in \mathcal{M}(\delta)}\left|\frac{1}{n}\sum_{i=1}^n\varepsilon_i\left[|\langle \bX_i,\bM\rangle|-|\langle \bX_i,\bM^\prime\rangle|\right]\right|\right)\le 	\E\left(\sup_{\bM^\prime\in \mathcal{M}(\delta)}\left|\frac{1}{n}\sum_{i=1}^n\varepsilon_i\langle \bX_i,\bM^\prime-\bM\rangle\right|\right)
\end{align*}
Hence we proceed as
\begin{align*}
	\E \tilde \theta_n(\delta)&\le 2\E\left(\sup_{\bM^\prime\in \mathcal{M}(\delta)}\left|\frac{1}{n}\sum_{i=1}^n\varepsilon_i\langle \bX_i,\bM^\prime-\bM\rangle\right|\right)=2\E\left(\sup_{\bM^\prime\in \mathcal{M}(\delta)}\left|\frac{1}{n}\sum_{i=1}^n\varepsilon_i\langle s_i\bM+\bZ_i,\bM^\prime-\bM\rangle\right|\right)\\
	&\le 2\E\left(\sup_{\bM^\prime\in \mathcal{M}(\delta)}\left|\left\langle \frac{1}{n}\sum_{i=1}^n\varepsilon_i\bM,\bM^\prime-\bM\right\rangle\right|\right)+2\E\left(\sup_{\bM^\prime\in \mathcal{M}(\delta)}\left|\left\langle \frac{1}{n}\sum_{i=1}^n\varepsilon_i\bZ_i,\bM^\prime-\bM\right\rangle\right|\right)\\
	&\overset{(a)}{\le} 2\delta\fro{\bM}\E\left|\frac{1}{n}\sum_{i=1}^n\varepsilon_i\right|+2\sqrt{2r}\delta\E\left\|\frac{1}{n}\sum_{i=1}^n\bZ_i\right\|\\
	&\overset{(b)}{\le} 2\delta\lambda\sqrt{\frac{r}{n}}+2\sqrt{2}\delta\sqrt{\frac{dr}{n}}\overset{(c)}{\lesssim } \delta \sqrt{\frac{dr}{n}}
\end{align*}
where in (a) we've used $\|\bM^\prime -\bM\|_*\le \text{rank}(\bM^\prime -\bM)\cdot \fro{\bM^\prime -\bM}\le \sqrt{2r}\fro{\bM^\prime -\bM}$, in (b) we have a simple bound for $\E\left|\frac{1}{n}\sum_{i=1}^n\varepsilon_i\right|\le 1/\sqrt{n}$ by Jensen's inequality, and (c) is due to the assumption $\lambda\lesssim  \sqrt{d}$. Define $\sigma^2:=\sup_{\bM^\prime\in \mathcal{M}(\delta)}\sum_{i=1}^n\E\left[|\langle \bX_i,\bM\rangle|-|\langle \bX_i,\bM^\prime\rangle|\right]^2$, notice that
\begin{align*}
	\left[|\langle \bX_i,\bM\rangle|-|\langle \bX_i,\bM^\prime\rangle|\right]^2&\le |\langle \bX_i,\bM-\bM^\prime\rangle|^2=|\langle s_i\bM+\bZ_i,\bM-\bM^\prime\rangle|^2\\
	&=\langle\bM,\bM-\bM^\prime\rangle^2+\langle\bZ_i,\bM-\bM^\prime\rangle^2+2s_i\langle\bM,\bM-\bM^\prime\rangle\langle\bZ_i,\bM-\bM^\prime\rangle
\end{align*}
Observe that $\langle\bZ_i,\bM-\bM^\prime\rangle\sim \mathcal{N}(0,\fro{\bM-\bM^\prime}^2)$, we have $\sigma^2\le n\fro{\bM}^2\delta^2+n\delta^2\lesssim n\lambda^2r\delta^2$, the last inequality is due to $\fro{\bM}\ge\lambda\sqrt{r}\gtrsim 1$ in this regime. Moreover, by Lemma 2.2.2 in \cite{van1996weak}, we have
\begin{align*}
	\left\|\max_i\sup_{\bM^\prime\in \mathcal{M}(\delta)}\left|\phi_i(\bM^\prime)-\E\phi_i(\bM^\prime)\right|\right\|_{\psi_2}\lesssim \sqrt{\log n}\max_i\left\|\sup_{\bM^\prime\in \mathcal{M}(\delta)}\left|\phi_i(\bM^\prime)-\E\phi_i(\bM^\prime)\right|\right\|_{\psi_2}
\end{align*}
It suffices to note that for each $i\in[n]$,
\begin{align*}
	\sup_{\bM^\prime\in \mathcal{M}(\delta)}\left|\phi_i(\bM^\prime)-\E\phi_i(\bM^\prime)\right|&\le \sup_{\bM^\prime\in \mathcal{M}(\delta)}\left|\langle\bX_i,\bM^\prime-\bM\rangle\right|+\sup_{\bM^\prime\in \mathcal{M}(\delta)}\E\left|\langle\bX_i,\bM^\prime-\bM\rangle\right|\\
	&\le \delta\fro{\bM}+\delta\sqrt{r}\op{\bZ_i}+\E\left(\delta\fro{\bM}+\delta\sqrt{r}\op{\bZ_i}\right)\\
	&\lesssim \delta\sqrt{dr}+\delta\sqrt{r}\op{\bZ_i}
\end{align*}
where in the last inequality we've used $\E\op{\bZ_i}\lesssim \sqrt{d}$ and $\lambda\lesssim \sqrt{d}$. By random matrix theory, we know $\E\op{\bZ_i}\asymp\sqrt{d}$ and $\op{\bZ_i}-\E\op{\bZ_i}$ is sub-gaussian, then  $\op{\op{\bZ_i}}_{\psi_2}\lesssim \sqrt{d}$. Hence
\begin{align*}
	\left\|\max_i\sup_{\bM^\prime\in \mathcal{M}(\delta)}\left|\phi_i(\bM^\prime)-\E\phi_i(\bM^\prime)\right|\right\|_{\psi_2}\lesssim \delta\sqrt{dr\log n}+\delta\sqrt{r\log n}\max_i\left\|\op{\bZ_i}\right\|_{\psi_2}\lesssim \delta\sqrt{dr\log n}
\end{align*}
Now we can invoke concentration inequality for suprema of empirical processes of unbounded functions \cite[Theorem 4]{adamczak2008tail}, we have for any $t\ge 0$:
\begin{align}\label{ineq:concentration_sup}
	\Prob\left(\tilde\theta_n(\delta)\geq C\delta\sqrt{\frac{dr\log n}{n}}+\delta\sqrt{\frac{dr}{n}}t\right)\le \exp\left(-ct^2\right)
\end{align}
Note that \eqref{ineq:concentration_sup} only holds for any given $\delta>0$. Consider any $\delta\in[\sqrt{dr/n},2\sqrt{dr}]$, let $\delta_j=2^{j}\sqrt{dr/n}$ for $j=0,1,\cdots,k^*+1$ with $k^*:=\lfloor\log_2(2\sqrt{n}) \rfloor$, then $\delta\in \bigcup_{j=1}^{k^*} [\delta_j,\delta_{j+1}]$. By construction, for any $\delta\in[\delta_j,\delta_{j+1}]$, we have $\delta \asymp \delta_j\asymp \delta_{j+1}$. Hence for fixed $j$, \eqref{ineq:concentration_sup} holds for any $\delta\in[\delta_j,\delta_{j+1}]$ up to change in constants $c>0$ and $C>0$. Take a union bound over all $j=0,1,\cdots,k^*+1$, we  have \eqref{ineq:concentration_sup} holds for any $\delta\in [\sqrt{dr/n},2\sqrt{dr}]$ and any $t\ge \log (k^*+2)\gtrsim \log\log n$. Combined with \eqref{ineq:mle_KL}, we conclude that
\begin{align}\label{ineq:kl-hp-bound}
D_{\textsf{KL}}\left(p_{\bM} \big \| p_{ \hat \bM_{\MLE}}\right)&\le \left|\frac{1}{n}\sum_{i=1}^n\log \frac{p_{\bM}(\bX_i)}{p_{\hat \bM_{\MLE}}(\bX_i)}-D_{\textsf{KL}}\left(p_{\bM} \big \| p_{ \hat \bM_{\MLE}}\right)\right|\le \tilde\theta_n(\fro{\hat \bM_{\MLE}-\bM})+2\log 2\nonumber \\
&\le  C\fro{\hat \bM_{\MLE}-\bM}\sqrt{\frac{dr\log (nd)}{n}}
\end{align}
where the last inequality holds with probability at least $1-(nd)^{-4}$, due to the facts that $\fro{\hat \bM_{\MLE}-\bM}\gtrsim \sqrt{dr\log (nd)/n}\ge 1$. Since $\fro{\hat \bM_{\MLE}-\bM}\gtrsim1$ and $\fro{\bM}\gtrsim1$ in this regime, it turns out that we can apply Lemma \ref{lem:KL} to get a lower bound of $D_{\textsf{KL}}\left(p_{\bM} \big \| p_{ \hat \bM_{\MLE}}\right)$, hence we have with probability at least $1-(nd)^{-4}$ that
\begin{align*}
	\fro{\hat \bM_{\MLE}-\bM}\le C\sqrt{\frac{dr\log (nd)}{n}}
\end{align*}
Finally, we can have a bound in expectation
\begin{align*}
	\E\fro{\hat \bM_{\MLE}-\bM}\le C\sqrt{\frac{dr\log (nd)}{n}}
\end{align*}
given that $\lambda\lesssim \sqrt{d}$.
\qed
\subsection{Proof of Theorem \ref{thm:spectral}}
In the proof, we consider the conditional model with sample splitting, i.e., $\bX_i^{(k)}\overset{d}{=}s_i^{(k)}\bM+\bZ_i^{(k)}$ for $k=1,2,3,4$ and $i\in[n_0]$. Let $\bU\bSigma\bV^\top $ denote the thin SVD of the signal matrix $\bM$ and recall that $d_1\asymp d_2\asymp d$.
\subsubsection*{Step 1:}
Denote $\bX_f=[\bX_1^{(1)},\cdots,\bX_{n_0}^{(1)}]\in \mathbb{R}^{d\times {n_0}d}$. A key observation is that $\hat \bu_1$ is also the leading eigenvector of $\frac{1}{n_0}\bX_f\bX_f^\top-d\bI_d$. Then we have
\begin{align}\label{step1-decomp}
	\frac{1}{n_0}\bX_f\bX_f^\top-d\bI_d&=\frac{1}{n_0}\sum_{i=1}^{n_0}\bX_i^{(1)}\bX_i^{(1)\top }-d\bI_d
%	&=\bM\bM^{\top}+\bM\left(\frac{1}{n_0}\sum_{i=1}^{n_0}s_i^{(1)}Z_i^{(1)\top }\right)+\left(\frac{1}{n_0}\sum_{i=1}^{n_0}s_i^{(1)}Z_i^{(1)}\right)\bM^{\top} +\frac{1}{n_0}\sum_{i=1}^{n_0}Z_i^{(1)}Z_i^{(1)\top}-dI_d\\
=\bM\bM^{\top }+\Delta
\end{align}
where 
$$	\Delta:=\bM\left(\frac{1}{n_0}\sum_{i=1}^{n_0}s_i^{(1)}\bZ_i^{(1)\top }\right)+\left(\frac{1}{n_0}\sum_{i=1}^{n_0}s_i^{(1)}\bZ_i^{(1)}\right)\bM^{\top} +\frac{1}{n_0}\sum_{i=1}^{n_0}\bZ_i^{(1)}\bZ_i^{(1)\top}-d\bI_{d_1}$$
Note that $\sum_{i=1}^{n_0}s_i^{(1)}\bZ_i^{(1)}$ ($\sum_{i=1}^{n_0}s_i^{(1)}\bZ_i^{(1)\top }$) is a $d_1\times d_2$ ($d_2\times d_1$) matrix of independent centered Gaussian entries with variance $n_0$, then by random matrix theory (e.g. \cite{vershynin2010introduction}) with probability at least $1-\exp(-cd)$, we have:
$$\left\|\bM\left(\frac{1}{n_0}\sum_{i=1}^{n_0}s_i^{(1)}\bZ_i^{(1)\top }\right)\right\|\lesssim \lambda_1\sqrt{\frac{d}{n}},\quad \left\|\left(\frac{1}{n_0}\sum_{i=1}^{n_0}s_i^{(1)}\bZ_i^{(1)}\right)\bM^{\top}\right\|\lesssim \lambda_1\sqrt{\frac{d}{n}}$$
Furthermore, since 
\begin{align*}
	\frac{1}{n_0}\sum_{i=1}^{n_0}\bZ_i^{(1)}\bZ_i^{(1)\top}-d\bI_{d_1}=d\left(\frac{1}{n_0d}\sum_{i=1}^{n_0}\sum_{j=1}^d[\bZ_{i}^{(1)}]_{:j}[\bZ_{i}^{(1)}]^\top_{:j}-\bI_{d_1}\right)
\end{align*}
where $[\bZ_{i}^{(1)}]_{:j}$ is the $j$-th column of $\bZ_{i}^{(1)}$. By concentration of sample covariance operator (e.g. \cite{koltchinskii2017concentration}), we have with probability at least $1-\exp(-cd)$:
\begin{align*}
	\left\|d\left(\frac{1}{n_0d}\sum_{i=1}^{n_0}\sum_{j=1}^d[\bZ_{i}^{(1)}]_{:j}[\bZ_{i}^{(1)}]^\top_{:j}-\bI_{d_1}\right)\right\|\lesssim\frac{d}{\sqrt{n}}
\end{align*}
By \eqref{step1-decomp} and eigenvalue perturbation theory (e.g. Corollary 8.1.6 in \cite{golub13}), we have\footnote{With slight abuse of notation, we use $\lambda_j(\cdot)$ to denote the $j$-th largest eigenvalue of a given matrix, while $\lambda_j$'s themselves are singular values of $\bM$.}
\begin{align*}
	\lambda_1^2-\|\Delta\|\le \lambda_1\left(\frac{1}{n_0}\bX_f\bX_f^\top-d\bI_{d_1}\right)\le \lambda_1^2+\|\Delta\|
\end{align*}
%\begin{align*}
%	\lambda_1\left(\frac{1}{n_0}X_fX_f^\top-dI_d\right)\ge\lambda_1\left(\bM\bM^{\top}\right)- \|\Delta\|=\lambda_1^2-\|\Delta\|
%\end{align*}
Therefore, we obtain 
\begin{align}\label{step1-lowerbound_np}
	\lambda_1^2-2\|\Delta\|\le \hat \bu_1^\top \bM\bM^{\top }\hat \bu_1\le \lambda_1^2+2\| \Delta \|
\end{align}
Hence with probability at least $1-\exp(-cd)$ we have
\begin{align}\label{step1-lowerbound}
	 \hat \bu_1^\top \bM\bM^{\top}\hat \bu_1\ge \lambda_1^2-2\| \Delta \|\ge \lambda_1^2-C\left(\lambda_1\sqrt{\frac{d}{n}}+\frac{d}{\sqrt{n}}\right)\gtrsim \lambda_1^2
\end{align}
where the last inequality holds provided that $\lambda^2\geq C_0\frac{d}{\sqrt{n}}$ for some large absolute constant $C_0>0$. 
\subsubsection*{Step 2:}
Observe that $\hat \bv_1$ is the leading eigenvector of $\frac{1}{n_0}\sum_{i=1}^{n_0}\bX_i^{(2)\top}\hat \bu_1\hat \bu_1^\top \bX_i^{(2)}-\bI_{d_2}$ and we have the following decomposition:
\begin{align}\label{step2-decomp}
	\frac{1}{n_0}\sum_{i=1}^{n_0}\bX_i^{(2)\top}\hat \bu_1\hat \bu_1^\top \bX_i^{(2)}-\bI_{d_2}=\bM^{\top}\hat \bu_1\hat \bu_1^\top \bM+\Delta^\prime
\end{align}
where 
$$\Delta^\prime:=\left(\frac{1}{n_0}\sum_{i=1}^{n_0}s_i^{(2)}\bZ_i^{(2)\top }\hat \bu_1\right)\hat \bu_1^\top \bM+\bM^{\top }\hat \bu_1\left(\frac{1}{n_0}\sum_{i=1}^{n_0}\hat \bu_1^\top s_i^{(2)}\bZ_i^{(2)}\right)+\frac{1}{n_0}\sum_{i=1}^{n_0}\bZ_i^{(2)\top }\hat \bu_1\hat \bu_1^\top \bZ_i^{(2)}-\bI_{d_2}$$
Due to the independence of $\hat \bu_1$ and $\{\bZ_i^{(2)}\}_{i=1}^{n_0}$, we conclude that $\frac{1}{n_0}\sum_{i=1}^{n_0}s_i^{(2)}\bZ_i^{(2)\top }\hat \bu_1\sim \calN(0,\frac{1}{n_0}\bI_{d_2})$, hence with probability at least $1-\exp(-cd)$:
\begin{align*}
	\left\|\left(\frac{1}{n_0}\sum_{i=1}^{n_0}s_i^{(2)}\bZ_i^{(2)\top }\hat \bu_1\right)\hat \bu_1^\top \bM\right\|\lesssim \lambda_1\sqrt{\frac{d}{n}}, \quad \left\|\bM^{\top }\hat \bu_1\left(\frac{1}{n_0}\sum_{i=1}^{n_0}\hat \bu_1^\top s_i^{(2)}\bZ_i^{(2)}\right)\right\|\lesssim \lambda_1\sqrt{\frac{d}{n}}
\end{align*}
Notice that $\frac{1}{n_0}\sum_{i=1}^{n_0}\bZ_i^{(2)\top }\hat \bu_1\hat \bu_1^\top \bZ_i^{(2)}\overset{d}{=}\frac{1}{n_0}\sum_{i=1}^{n_0}\bz_i\bz_i^\top$, where $\bz_i\sim \calN(0,\bI_{d_2})$ and $\bz_i$'s are independent. Again, by concentration of sample covariance operator, with probability at least $1-\exp(-cd)$:
\begin{align*}
	\left\|\frac{1}{n_0}\sum_{i=1}^{n_0}\bZ_i^{(2)\top }\hat \bu_1\hat \bu_1^\top \bZ_i^{(2)}-\bI_{d_2}\right\|\lesssim \sqrt{\frac{d}{n}}\vee \frac{d}{n}
\end{align*}
Therefore, \eqref{step2-decomp} and eigenvalue perturbation theory imply that
\begin{align*}
	\lambda_1\left(\bM^{\top}\hat \bu_1\hat \bu_1^\top \bM\right)-\|\Delta^\prime\|\le \lambda_1\left(\frac{1}{n_0}\sum_{i=1}^{n_0}\bX_i^{(2)\top}\hat \bu_1\hat \bu_1^\top \bX_i^{(2)}-\bI_{d_2}\right)\le \lambda_1\left(\bM^{\top}\hat \bu_1\hat \bu_1^\top \bM\right)+\|\Delta^\prime\|
\end{align*}
%Denote $\hat{\lambda_1}^2=\lambda_1 \left(\frac{1}{n_0}\sum_{i=1}^{n_0} \bX_i^{(2)}\hat \bu_1\hat \bu_1^\top \bX_i^{(2)}-I_d\right)$
%, by \eqref{step1-lowerbound_np} we have
%\begin{align*}
%	\lambda_1^2-2\|\Delta\|-\|\Delta^\prime\|\le \hat \lambda_1^2\le \lambda_1^2+2\|\Delta\|+\|\Delta^\prime\|
%\end{align*}
Combined with the decomposition \eqref{step2-decomp}, we can arrive at
\begin{align}\label{step2-lowerbound_np}
	\lambda_1^2-2\left(\|\Delta\| +\|\Delta^\prime\|\right)\le  \hat \bv_1^\top \bM^{\top}\hat \bu_1\hat \bu_1^\top \bM\hat \bv_1\le \lambda_1^2+2\left(\|\Delta\| +\|\Delta^\prime\|\right)
\end{align}
Thus, we get
\begin{align}\label{step2-lowerbound}
|\hat \bu_1^\top \bM\hat \bv_1|\asymp\lambda_1
\end{align}
with probability at least $1-\exp(-cd)$, provided that $\lambda^2\ge C_0\frac{d}{\sqrt{n}}$.
\subsubsection*{Step 3:}
We proceed our analysis by conditioning on the event $\{\eqref{step2-lowerbound} \text{~holds}\}$. Observe that
\begin{align}\label{step3-decomp}
\frac{1}{n_0}\sum_{i=1}^{n_0}(\hat \bu_1^\top \bX_i^{(3)}\hat \bv_1)\bX_i^{(3)}-\hat \bu_1\hat \bv_1^\top &=\frac{1}{n_0}\sum_{i=1}^{n_0}(\hat \bu_1^\top (s_i^{(3)}\bM+\bZ_i^{(3)})\hat \bv_1)(s_i^{(3)}\bM+\bZ_i^{(3)})\nonumber\\
&=:(\hat \bu_1^\top \bM\hat \bv_1)\bM+\Upsilon
\end{align}
where
$$\Upsilon:=\hat \bu_1^\top \bM\hat \bv_1\left( \frac{1}{n_0}\sum_{i=1}^{n_0}s_i^{(3)}\bZ_i^{(3)}\right)+\bM \left(\frac{1}{n_0}\sum_{i=1}^{n_0}s_i^{(3)}(\hat \bu_1^\top \bZ_i^{(3)}\hat \bv_1)\right)+\frac{1}{n_0}\sum_{i=1}^{n_0}(\hat \bu_1^\top \bZ_i^{(3)}\hat \bv_1)\bZ_i^{(3)}-\hat \bu_1\hat \bv_1^\top $$
Now we give an upper bound for $\|\Upsilon\|$. By random matrix theory we know with probability at least $1-\exp(-cd)$:
$$\left\|\hat \bu_1^\top \bM\hat \bv_1\left( \frac{1}{n_0}\sum_{i=1}^{n_0}s_i^{(3)}\bZ_i^{(3)}\right)\right\|\lesssim \lambda_1\sqrt{\frac{d}{n}}$$
Next, notice that $\frac{1}{n_0}\sum_{i=1}^{n_0}s_i^{(3)}(\hat \bu_1^\top \bZ_i^{(3)}\hat \bv_1)\sim \calN(0,\frac{1}{n_0})$, we have with probability at least $1-\exp(-cd)$:
$$\left\|\bM\left(\frac{1}{n_0}\sum_{i=1}^{n_0}s_i^{(3)}(\hat \bu_1^\top \bZ_i^{(3)}\hat \bv_1)\right)\right\|\lesssim  \lambda_1\sqrt{\frac{d}{n}}$$
It remains to bound $\frac{1}{n_0}\sum_{i=1}^{n_0}(\hat \bu_1^\top \bZ_i^{(3)}\hat \bv_1)\bZ_i^{(3)}-\hat \bu_1\hat \bv_1^\top$. Notice the following decomposition:
\begin{align*}
\frac{1}{n_0}\sum_{i=1}^{n_0}(\hat \bu_1^\top \bZ_i^{(3)}\hat \bv_1)\bZ_i^{(3)}=\frac{1}{n_0}\sum_{i=1}^{n_0}(\hat \bu_1^\top \bZ_i^{(3)}\hat \bv_1)\left[\mathcal{P}_{\hat \bu_1}\bZ_i^{(3)} \mathcal{P}_{\hat \bv_1}+\mathcal{P}^\perp_{\hat \bu_1}\bZ_i^{(3)} \mathcal{P}_{\hat \bv_1}+\mathcal{P}^\perp_{\hat \bu_1}\bZ_i^{(3)} \mathcal{P}^\perp_{\hat \bv_1}+\mathcal{P}_{\hat \bu_1}\bZ_i^{(3)}\mathcal{P}^\perp_{\hat \bv_1}\right]
\end{align*}
where $\mathcal{P}_\bu$ is the projection matrix onto the column space of $\bu$ and $\mathcal{P}_\bu^\perp$ is the projection matrix onto orthogonal complement of the column space of $\bu$. Since $\sum_{i=1}^{n_0}(\hat \bu_1^\top \bZ_i^{(3)}\hat \bv_1)^2\sim \chi^2_{n_0}$, by concentration for chi-square random variable with $n_0$ degrees of freedom (see \cite{laurent2000adaptive}), we have with probability at least {$1-\exp(-c\sqrt{d(d\wedge n)})$}:
\begin{equation*}
\left\|\frac{1}{n_0}\sum_{i=1}^{n_0}(\hat \bu_1^\top \bZ_i^{(3)}\hat \bv_1)\mathcal{P}_{\hat \bu_1}\bZ_i^{(3)} \mathcal{P}_{\hat \bv_1}-\hat \bu_1\hat \bv_1^\top\right\|=\left\|\hat \bu_1\hat \bv_1^\top\left( \frac{1}{n_0}\sum_{i=1}^{n_0}(\hat \bu_1^\top \bZ_i^{(3)}\hat \bv_1)^2-1\right)\right\|\lesssim {\sqrt{\frac{d}{n}}}
\end{equation*}
By property of Gaussian matrices we have $(\hat \bu_1^\top \bZ_i^{(3)}\hat \bv_1)\mathcal{P}^\perp_{\hat \bu_1}\bZ_i^{(3)} \mathcal{P}_{\hat \bv_1}\stackrel{d}{=}g_i\bZ_i^{(3)}$, where $g_i\overset{i.i.d}{\sim}\calN(0,1)$ and $\{g_i\}_{i=1}^{n_0}$ is independent of $\{\bZ_i^{(3)}\}_{i=1}^{n_0}$. Hence we have
\begin{equation}\label{step3-eq2}
\Prob\left(\left\|\frac{1}{n_0}\sum_{i=1}^{n_0}(\hat \bu_1^\top \bZ_i^{(3)}\hat \bv_1)\mathcal{P}^\perp_{\hat \bu_1}\bZ_i^{(3)} \mathcal{P}_{\hat \bv_1}\right\|\leq \frac{\sqrt{d}}{n}\sqrt{\sum_{i=1}^{n_0}g_i^2} \Bigg| \{g_i\}_{i=1}^{n_0}\right)\ge  1-\exp(-cd)	
\end{equation}
In addition, by concentration for chi-square random variable we have $\sqrt{\sum_{i=1}^{n_0}g_i^2}\lesssim \sqrt{n}$ with probability at least $1-\exp(-cn)$. Combined with \eqref{step3-eq2}, we arrive at
\[\left\|\frac{1}{n_0}\sum_{i=1}^{n_0}(\hat \bu_1^\top \bZ_i^{(3)}\hat \bv_1)\mathcal{P}^\perp_{\hat \bu_1}\bZ_i^{(3)} \mathcal{P}_{\hat \bv_1}\right\|\lesssim \sqrt{\frac{d}{n}}\]
with probability at least $1-\exp(-c(d\wedge n))$. Similar arguments can be applied to $(\hat \bu_1^\top \bZ_i^{(3)}\hat \bv_1)\mathcal{P}^\perp_{\hat \bu_1}\bZ_i^{(3)} \mathcal{P}^\perp_{\hat \bv_1}$ and $(\hat \bu_1^\top \bZ_i^{(3)}\hat \bv_1)\mathcal{P}_{\hat \bu_1}\bZ_i^{(3)}\mathcal{P}^\perp_{\hat \bv_1}$. Collecting four parts we can bound the last term of $\Upsilon$ as
\begin{align}\label{step3-2ndorderbound}
\left\|\frac{1}{n_0}\sum_{i=1}^{n_0}(\hat \bu_1^\top \bZ_i^{(3)}\hat \bv_1)\bZ_i^{(3)}-\hat \bu_1\hat \bv_1^\top \right\|\lesssim \sqrt{\frac{d}{n}}
\end{align}
with probability at least $1-\exp(-c(d\wedge n))$. Hence we have the following bound for $\|\Upsilon\|$ with probability at least  {$1-\exp(-c(d\wedge n))$}:
\begin{align}\label{step3-upsilonbound}
	\|\Upsilon\|\lesssim \lambda_1\sqrt{\frac{d}{n}}+\sqrt{\frac{d}{n}}
\end{align}
%\textcolor{purple}
%{
%Now by Wedin's $\sin\Theta$ theorem, we have
%\begin{equation}\label{step3-uverror}
%\max\{\|\sin\Theta(\tilde \bU,U)\|,\|\sin\Theta(\tilde \bV,V)\|\}\lesssim \frac{\|\Upsilon\|}{\sigma_{r}\left((\hat \bu_1^\top \bM\hat \bv_1)\bM\right)}
%\end{equation}
%Notice that by \eqref{step2-lowerbound}, with probability at least $1-\exp(-cd)$:
%\begin{align*}
%	\sigma_{r}\left((\hat \bu_1^\top \bM\hat \bv_1)\bM\right)= |\hat \bu_1^\top \bM\hat \bv_1|\sigma_{r}(\bM)\gtrsim \lambda_1\lambda
%\end{align*}
%Together with \eqref{step3-upsilonbound} we arrive at
%\begin{align}\label{step3-uverror}
%\max\{\|\sin\Theta(\tilde \bU,U)\|,\|\sin\Theta(\tilde \bV,V)\|\}\lesssim \frac{\lambda_1\sqrt{\frac{d}{n}}+\sqrt{\frac{d}{n}}}{\lambda_1\lambda}\lesssim \frac{\sqrt{d/n}}{\lambda}+\frac{\sqrt{d/n}}{\lambda_1\lambda}
%\end{align}
%}
{
For any $j\in[r]$, denote $\tilde\lambda_j$ is the $j$-th largest singular value of $\frac{1}{n_0}\sum_{i=1}^{n_0}(\hat \bu_1^\top \bX_i^{(3)}\tilde \bV_1)\bX_i^{(3)}-\hat \bu_1\hat \bv_1^\top$ and $\tilde \bu_j,\tilde \bv_j$ the corresponding left and right singular vectors.  By \eqref{step3-decomp} and perturbation theory for singular values we have
\begin{align}\label{step3-lambdatildebound}
	\sigma_j\left((\hat \bu_1^\top \bM\hat \bv_1) \bM\right)-\|\Upsilon\| \le \tilde\lambda_j\le\sigma_j\left((\hat \bu_1^\top \bM\hat \bv_1) \bM\right)+\|\Upsilon\|
\end{align}
By definition of singular value and singular vectors, we have that
\begin{align*}
\tilde\lambda_j=\tilde \bu_j^\top\left(\frac{1}{n_0}\sum_{i=1}^{n_0}(\hat \bu_1^\top \bX_i^{(3)}\hat \bv_1)\bX_i^{(3)}-\hat \bu_1\hat \bv_1^\top\right)\tilde \bv_j=(\hat \bu_1^\top \bM\hat \bv_1)(\tilde \bu_j^\top \bM\tilde \bv_j)+\tilde \bu_j^\top\Upsilon\tilde \bv_j
\end{align*} 
which implies 
\begin{align}\label{step3-bound}
\sigma_j\left((\hat \bu_1^\top \bM\hat \bv_1) \bM\right)-2\|\Upsilon\|\le (\hat \bu_1^\top \bM\hat \bv_1)(\tilde \bu_j^\top \bM\tilde \bv_j)\le \sigma_j\left((\hat \bu_1^\top \bM\hat \bv_1) \bM\right)+2\|\Upsilon\|
\end{align}
Using \eqref{step2-lowerbound}, it follows that with probability at least  {$1-\exp(-c(d\wedge n))$} such that for all $j\in[r]$:
\begin{align}\label{step3-absbound}
\left|\left|\tilde \bu_j^\top \bM\tilde \bv_j\right|- \lambda_j\right |\lesssim \frac{\|\Upsilon\|}{\lambda_1}=o(\lambda)
\end{align}
given that $\lambda^2\gtrsim \frac{d}{\sqrt{n}}$. Notice that this implies with overwhelming probability, we have
\begin{enumerate}
	\item $|\tilde \bu_j^\top \bM\tilde \bv_j|\asymp \lambda_j$
	\item $\tilde \bu_j^\top \bM\tilde \bv_j$ share the same sign with $\hat \bu_1^\top \bM\hat \bv_1$
\end{enumerate} 
As a consequence, $|\sum_{j=1}^r\tilde \bu_j^\top \bM\tilde \bv_j|\asymp \sum_{j=1}^r\lambda_j$ with probability at least  $1-\exp(-c(d\wedge n))$. These facts will be used in the following derivations.
\subsubsection*{Step 4:}
We proceed by conditioning on the event $\{\eqref{step3-bound},\eqref{step3-absbound} \text{~holds}\}$. Consider the rank-$r$ approximation of $\frac{1}{n_0}\sum_{i=1}^{n_0}\left(\sum_{j=1}^r\tilde \bu_j^\top \bX_i^{(4)}\tilde \bv_j\right)\bX_i^{(4)}-\sum_{j=1}^r\tilde \bu_j\tilde \bv_j^\top$, which admits the following decomposition:
\begin{align}\label{step4-decomp}
\frac{1}{n_0}\sum_{i=1}^{n_0}\left(\sum_{j=1}^r\tilde \bu_j^\top \bX_i^{(4)}\tilde \bv_j\right)\bX_i^{(4)}-\sum_{j=1}^r\tilde \bu_j\tilde \bv_j^\top&=\frac{1}{n_0}\sum_{i=1}^{n_0}\left(\sum_{j=1}^r(\tilde \bu_j^\top (s_i^{(4)}\bM+\bZ_i^{(4)})\tilde \bv_j\right)(s_i^{(4)}\bM+\bZ_i^{(4)})\nonumber\\
&=:\sum_{j=1}^r\left(\tilde \bu_j^\top \bM\tilde \bv_j\right)\bM+\Upsilon^\prime
\end{align}
where
\begin{align*}
\Upsilon^\prime&:=\sum_{j=1}^r\left(\tilde \bu_j^\top \bM\tilde \bv_j\right)\left( \frac{1}{n_0}\sum_{i=1}^{n_0}s_i^{(4)}\bZ_i^{(4)}\right)+\bM \left(\sum_{j=1}^r\tilde \bu_j^\top\left(\frac{1}{n_0}\sum_{i=1}^{n_0}s_i^{(4)}\bZ_i^{(4)}\right)\tilde \bv_j\right)\\
&+\frac{1}{n_0}\sum_{i=1}^{n_0}\left(\sum_{j=1}^r\tilde \bu_j^\top \bZ_i^{(4)}\tilde \bv_j\right)\bZ_i^{(4)}-\sum_{j=1}^r\tilde \bu_j\tilde \bv_j^\top 	
\end{align*}
Similar to that in step 3, we need to upper bound $\|\Upsilon^\prime\|$, the spectral norm of the perturbation term. The following bound is clear, which holds with probability at least $1-\exp(-cd)$:
\begin{align*}
\left\|\sum_{j=1}^r\left(\tilde \bu_j^\top \bM\tilde \bv_j\right)\left( \frac{1}{n_0}\sum_{i=1}^{n_0}s_i^{(4)}\bZ_i^{(4)}\right)\right\|\lesssim \sum_{j=1}^r\left(\tilde \bu_j^\top \bM\tilde \bv_j\right)\sqrt{\frac{d}{n}}\lesssim \left(\sum_{j=1}^r\lambda_j\right)\sqrt{\frac{d}{n}}
\end{align*}
Due to the rotation invariance of Gaussian and the orthogonality of $\tilde \bu_j$'s and $\tilde \bv_j$'s, we have with probability at least $1-\exp(-cd/r)$:
\begin{align*}
\left\|\bM \left(\sum_{j=1}^r\tilde \bu_j^\top\left(\frac{1}{n_0}\sum_{i=1}^{n_0}s_i^{(4)}\bZ_i^{(4)}\right)\tilde \bv_j\right)\right\|\lesssim \lambda_1\sqrt{\frac{d}{n}}
\end{align*}
The following decomposition is similar to that in step 3:
\begin{align}\label{step4-decomp2}
\frac{1}{n_0}\sum_{i=1}^{n_0}\left(\sum_{j=1}^r\tilde \bu_j^\top \bZ_i^{(4)}\tilde \bv_j\right)\bZ_i^{(4)}=\frac{1}{n_0}\sum_{i=1}^{n_0}\left(\sum_{j=1}^r\tilde \bu_j^\top \bZ_i^{(4)}\tilde \bv_j\right)\left[\mathcal{P}_{\tilde \bU}\bZ_i^{(4)} \mathcal{P}_{\tilde \bV}+\mathcal{P}^\perp_{\tilde \bU}\bZ_i^{(4)} \mathcal{P}_{\tilde \bV}+\mathcal{P}^\perp_{\tilde \bU}\bZ_i^{(4)} \mathcal{P}^\perp_{\tilde \bV}+\mathcal{P}_{\tilde \bU}\bZ_i^{(4)}\mathcal{P}^\perp_{\tilde \bV}\right]
\end{align}
By the property of Gaussian matrices, we have
\begin{align*}
\left\|\frac{1}{n_0}\sum_{i=1}^{n_0}\left (\sum_{j=1}^r\tilde \bu_j^\top \bZ_i^{(4)}\tilde \bv_j\right)\mathcal{P}_{\tilde \bU}\bZ_i^{(4)} \mathcal{P}_{\tilde \bV}\right\|{=}\left\|\frac{1}{n_0}\sum_{i=1}^{n_0}\text{Tr}\left (\tilde \bU^\top \bZ_i^{(4)}\tilde \bV\right){\tilde \bU^\top }\bZ_i^{(4)}\tilde \bV	\right\|\overset{d}{=}\left\|\frac{1}{n_0}\sum_{i=1}^{n_0}\text{Tr}\left (\bZ_{r,i}\right)\bZ_{r,i}\right\|
\end{align*}
where $\{\bZ_{r,i}\}_{i=1}^{n_0}$ are independent matrices of dimension $r\times r$ with i.i.d standard normal entries. Hence 
\begin{align*}
\left\|\frac{1}{n_0}\sum_{i=1}^{n_0}\left (\sum_{j=1}^r\tilde \bu_j^\top \bZ_i^{(4)}\tilde \bv_j\right)\mathcal{P}_{\tilde \bU}\bZ_i^{(4)} \mathcal{P}_{\tilde \bV}-\sum_{j=1}^r\tilde \bu_j\tilde \bv_j^\top\right\|\overset{d}{=}\left\|\frac{1}{n_0}\sum_{i=1}^{n_0}\text{Tr}\left (\bZ_{r,i}\right)\bZ_{r,i}-\bI_r\right\|
\end{align*}
The following lemma gives the concentration inequality of the above term.
\begin{lemma}\label{lemma-concentration}
	Let $\bZ,\bZ_1,\cdots,\bZ_n$ be $r\times r$ independent matrices with i.i.d standard normal entries. Then there exists a constant $C>0$ such that, for all $t>0$,  with probability  at least $1-e^{-t}$:
	$$\left\|\frac{1}{n}\sum_{i=1}^{n}\text{Tr}\left(\bZ_i\right)\bZ_i-\bI_r\right\| \le C\left(r\sqrt{\frac{t+\log(2r)}{n}}+r\frac{t+\log(2r)}{n}\right)$$
\end{lemma}
\noindent By Lemma \ref{lemma-concentration}, if $d<nr$, we take $t=d/r$, then we have $\left\|\frac{1}{n_0}\sum_{i=1}^{n_0}\text{Tr}\left(\bZ_{i,r}\right)\bZ_{i,r}-\bI_r\right\| \lesssim \sqrt{\frac{dr}{n}}$ with probability at least $1-\exp(-d/r)$, provided that $d\gtrsim r\log r$. If $d>nr$, we can take $t=\sqrt{nd/r}$, then we have $\left\|\frac{1}{n_0}\sum_{i=1}^{n_0}\text{Tr}\left(\bZ_{i,r}\right)\bZ_{i,r}-\bI_r\right\| \lesssim \sqrt{\frac{dr}{n}}$ with probability at least $1-\exp(-n)$, provided that $nd\gtrsim r^2\log^2r$. In summary, we have with probability at least $1-\exp(-c(d/r\wedge n))$:
\begin{align*}
\left\|\frac{1}{n_0}\sum_{i=1}^{n_0}\text{Tr}\left(\bZ_{i,r}\right)\bZ_{i,r}-\bI_r\right\| \lesssim \sqrt{\frac{dr}{n}}
\end{align*}
In addition, we have $\left(\sum_{j=1}^r\tilde \bu_j^\top \bZ_i^{(4)}\tilde \bv_j\right)\mathcal{P}^\perp_{\tilde \bU}\bZ_i^{(4)} \mathcal{P}_{\tilde \bV}\overset{d}{=}\sqrt{r}g_{i}\bZ_{i}^{(4)}$, where $g_{i}\overset{i.i.d}{\sim} N(0,1)$ and $\{g_{i}\}_{i=1}^{n_0}$ is independent of $\{\bZ_i^{(4)}\}_{i=1}^{n_0}$. Hence we have
\begin{align*}
\Prob\left(\left\|\frac{1}{n_0}\sum_{i=1}^{n_0}\left(\sum_{j=1}^r\tilde \bu_j^\top \bZ_i^{(4)}\tilde \bv_j\right)\mathcal{P}^\perp_{\tilde \bU}\bZ_i^{(4)} \mathcal{P}_{\tilde \bV}\right\| \ge \frac{\sqrt{dr}}{n}\sqrt{\sum_{i=1}^{n_0}g_i^2}\Bigg|\{g_i\}_{i=1}^{n_0}\right)\le \exp(-cd)
\end{align*}
Note that by concentration for chi-square random variable $\sqrt{\sum_{i=1}^{n_0}g_i^2}\lesssim \sqrt{n}$ with probability at least $1-\exp(-cn)$. Then we can conclude that with probability at least $1-\exp(-c(d\wedge n))$:
\[\left\|\frac{1}{n_0}\sum_{i=1}^{n_0}\left(\sum_{j=1}^r\tilde \bu_j^\top \bZ_i^{(4)}\tilde \bv_j\right)\mathcal{P}^\perp_{\tilde \bU}\bZ_i^{(4)} \mathcal{P}_{\tilde \bV}\right\|\lesssim \sqrt{\frac{dr}{n}}\] 
The bounds for $\frac{1}{n_0}\sum_{i=1}^{n_0}\left(\sum_{j=1}^r\tilde \bu_j^\top \bZ_i^{(4)}\tilde \bv_j\right)\mathcal{P}^\perp_{\tilde \bU}\bZ_i^{(4)} \mathcal{P}^\perp_{\tilde \bV}$ and $\frac{1}{n_0}\sum_{i=1}^{n_0}\left(\sum_{j=1}^r\tilde \bu_j^\top \bZ_i^{(4)}\tilde \bv_j\right)\mathcal{P}_{\tilde \bU}\bZ_i^{(4)}\mathcal{P}^\perp_{\tilde \bV}$ can be obtained similarly. We have with probability at least $1-\exp(-c(d/r\wedge n))$:
\begin{align}\label{step4-upsilonprimebound}
	\|\Upsilon^\prime\|\lesssim \left(\sum_{j=1}^r\lambda_j\right)\sqrt{\frac{d}{n}}+\sqrt{\frac{dr}{n}}
\end{align}
\subsubsection*{Step 5:}
Again, we continue on the event $\{\eqref{step4-upsilonprimebound} \text{~holds}\}$. In this step, we first construct $\hat\Lambda$, which is an estimator for the pre-factor $\Lambda^*:=|\sum_{j=1}^r\tilde \bu_j^\top \bM\tilde \bv_j|$ of the signal part in \eqref{step4-decomp}. Notice that
\begin{align}\label{step5-decomp}
	\frac{1}{n_0}&\left[\sum_{i=1}^{n_0}\left(\sum_{j=1}^r\tilde \bu_j^\top \bX_i^{(4)}\tilde \bv_j\right)^2-r\right]\\
	&=\left(\sum_{j=1}^r\tilde \bu_j^\top \bM\tilde \bv_j\right)^2+\frac{1}{n_0}\sum_{i=1}^{n_0}\left(\sum_{j=1}^r\tilde \bu_j^\top \bM\tilde \bv_j\right)\left(\sum_{j=1}^r\tilde \bu_j^\top \bZ_i^{(4)}\tilde \bv_j\right)+\frac{1}{n_0}\sum_{i=1}^{n_0}\left(\sum_{j=1}^r\tilde \bu_j^\top \bZ_i^{(4)}\tilde \bv_j\right)^2-r
\end{align}
The second term $\frac{1}{n_0}\sum_{i=1}^{n_0}\left(\sum_{j=1}^r\tilde \bu_j^\top \bM\tilde \bv_j\right)\left(\sum_{j=1}^r\tilde \bu_j^\top \bZ_i^{(4)}\tilde \bv_j\right)\overset{d}{=}\left|\sum_{j=1}^r\tilde \bu_j^\top \bM\tilde \bv_j\right|\sqrt{\frac{r}{n_0}}g$, with $g$ being standard normal. Therefore, with probability at least $1-\exp(-cd)$:
\begin{align*}
	\left|\frac{1}{n_0}\sum_{i=1}^{n_0}\left(\sum_{j=1}^r\tilde \bu_j^\top \bM\tilde \bv_j\right)\left(\sum_{j=1}^r\tilde \bu_j^\top \bZ_i^{(4)}\tilde \bv_j\right)\right|\lesssim \left(\sum_{j=1}^r\lambda_j\right)\sqrt{\frac{dr}{n}}
\end{align*}
The third term $\frac{1}{n_0}\sum_{i=1}^{n_0}\left(\sum_{j=1}^r\tilde \bu_j^\top \bZ_i^{(4)}\tilde \bv_j\right)^2-r\overset{d}{=}\frac{r}{n_0}(b-n_0)$, where $b{\sim }\chi_{n_0}^2$, hence we have with probability at least {$1-\exp(-c\sqrt{d(d\wedge n)})$}:
\begin{align*}
\left|\frac{1}{n_0}\sum_{i=1}^{n_0}\left(\sum_{j=1}^r\tilde \bu_j^\top \bZ_i^{(4)}\tilde \bv_j\right)^2-r\right|\lesssim r\sqrt{\frac{d}{n}}
\end{align*}
Therefore, we have with probability at least $1-\exp(-c\sqrt{d(d\wedge n)})$:
\begin{align*}
	\frac{1}{n_0}\left[\sum_{i=1}^{n_0}\left(\sum_{j=1}^r\tilde \bu_j^\top \bX_i^{(4)}\tilde \bv_j\right)^2-r\right]\gtrsim \left(\sum_{j=1}^r\lambda_j\right)^2-\left(\sum_{j=1}^r\lambda_j\right)\sqrt{\frac{dr}{n}}-r\sqrt{\frac{d}{n}}\gtrsim\frac{dr^2}{\sqrt{n}}
\end{align*}
where we've used the fact $\sum_{j=1}^r\lambda_j\gtrsim \frac{\sqrt{d}r}{n^{1/4}}$. Hence by definition of $\hat \Lambda$, with probability at least $1-\exp(-c\sqrt{d(d\vee n)})$:
\begin{align}\label{step5-lambda-max}
	\hat \Lambda^2=\max\left\{\frac{1}{n_0}\left[\sum_{i=1}^{n_0}\left(\sum_{j=1}^r\tilde \bu_j^\top \bX_i^{(4)}\tilde \bv_j\right)^2-r\right],\frac{dr^2}{\sqrt{n}}\right\}=\frac{1}{n_0}\left[\sum_{i=1}^{n_0}\left(\sum_{j=1}^r\tilde \bu_j^\top \bX_i^{(4)}\tilde \bv_j\right)^2-r\right]
\end{align}
The concentration inequalities of second and third term of \eqref{step5-decomp} also imply that with probability at least {$1-\exp(-c\sqrt{d(d\vee n)})$}:
\begin{align}\label{step5-lambdaerror}
	|\hat \Lambda-\Lambda^*|=\frac{|\hat \Lambda^2-\Lambda^{*2}|}{\hat \Lambda+\Lambda^*}\lesssim \frac{\left(\sum_{j=1}^r\lambda_j\right)\sqrt{\frac{dr}{n}}+r\sqrt{\frac{d}{n}}}{\sum_{j=1}^r\lambda_j}\le \sqrt{\frac{dr}{n}}+\frac{r}{\sum_{j=1}^r\lambda_j}\sqrt{\frac{d}{n}}
\end{align}
Since the RHS of \eqref{step5-lambdaerror} is of order $o\left(\sum_{j=1}^r\lambda_j\right)$, we have $\hat \Lambda\asymp\Lambda^*\asymp \sum_{j=1}^r\lambda_j$ with probability at least $1-\exp(-c\sqrt{d(d\vee n)})$. Next, denote $\check \bU$, $\check \bV$ the left and right leading $r$ singular vectors of$\frac{1}{n_0}\sum_{i=1}^{n_0}\left(\sum_{j=1}^r\tilde \bu_j^\top \bX_i^{(4)}\tilde \bv_j\right)\bX_i^{(4)}-\sum_{j=1}^r\tilde \bu_j\tilde \bv_j^\top$, then the best rank-$r$ approximation is given by
\begin{align*}
\check \bM&=\check \bU\check \bU^\top\left(\frac{1}{n_0}\sum_{i=1}^{n_0}\left(\sum_{j=1}^r\tilde \bu_j^\top \bX_i^{(4)}\tilde \bv_j\right)\bX_i^{(4)}-\sum_{j=1}^r\tilde \bu_j\tilde \bv_j^\top\right)\check \bV\check \bV^\top\\
&=\sum_{j=1}^r\left(\tilde \bu_j^\top \bM\tilde \bv_j\right)\check \bU\check \bU^\top \bM\check \bV\check \bV^\top+\check \bU\check \bU^\top \Upsilon^\prime\check \bV\check \bV^\top
\end{align*}
The error of low-rank approximation is characterized by the perturbation term, given by the following lemma.
\begin{lemma}\label{lemma-lowrank-asymmetric}
	Consider a rank-$r$ matrix $\bM\in\mathbb{R}^{d_1\times d_2}$ with its thin-SVD form $\bU\bSigma \bV^\top $, where $\bU\in\mathbb{O}_{d_1,r},\bV\in\mathbb{O}_{d_2,r}$ and $\bSigma=\text{diag}(\sigma_1,\cdots,\sigma_r)$, $\sigma_1\geq\sigma_2\ge\cdots \sigma_r>0$, let $E$ be a $d_1\times d_2$ perturbation matrix and $\hat \bM=\bM+\bE$. Denote $\hat \bM_r$ the best rank-$r$ approximation of $\hat \bM$. Suppose that $\sigma_r\ge 3\|\bE\|$, then there exists some absolute constant $C_0>0$ such that
	\begin{align*}
		\fro{\hat \bM_r-\bM}\le C_0\min\{\fro{\bE},\sqrt{r}\|\bE\|\}
	\end{align*}
\end{lemma}
\noindent By Lemma \ref{lemma-lowrank-asymmetric}, we have
\begin{align*}
\left\|\check \bM-\sum_{j=1}^r\left(\tilde \bu_j^\top \bM\tilde \bv_j\right)\bM\right\|_{\text{F}}\lesssim \sqrt{r}\|\Upsilon^\prime\|
\end{align*}
Recall that $\hat \bM=\check \bM/\hat \Lambda$. Denote $\eta^*=\text{sign}\left(\sum_{j=1}^r\tilde \bu_j^\top \bM\tilde \bv_j\right)$,  hence we have the following bound
\begin{align}\label{step5-error-decomp}
\fro{\hat \bM-\eta^*\bM}\le \left\| \hat\Lambda^{-1}\check \bM-\hat\Lambda^{-1}\sum_{j=1}^r\left(\tilde \bu_j^\top \bM\tilde \bv_j\right)\bM\right\|_{\text{F}}+\left\|\hat\Lambda^{-1}\sum_{j=1}^r\left(\tilde \bu_j^\top \bM\tilde \bv_j\right)\bM-\eta^*\bM\right\|_{\text{F}}
\end{align}
Using \eqref{step4-upsilonprimebound} and \eqref{step5-lambdaerror}, the first term can be bounded with probability at least $1-\exp(-c(d/r\wedge n))$:
\begin{align*}
\left\| \hat\Lambda^{-1}\check \bM-\hat\Lambda^{-1}\sum_{j=1}^r\left(\tilde \bu_j^\top \bM\tilde \bv_j\right)\bM\right\|_{\text{F}}\lesssim \frac{\sqrt{r}\|\Upsilon^\prime\|}{\hat\Lambda}\lesssim \sqrt{\frac{dr}{n}}+\frac{r}{\sum_{j=1}^r\lambda_j}\sqrt{\frac{d}{n}}
\end{align*}
Using \eqref{step5-lambdaerror}, the second term can be bounded with probability at least $1-\exp(-c\sqrt{d(d\wedge n)})$:
\begin{align*}
\left\|\hat\Lambda^{-1}\sum_{j=1}^r\left(\tilde \bu_j^\top \bM\tilde \bv_j\right)\bM-\eta^*\bM\right\|_{\text{F}}= \left|\frac{\Lambda^*-\hat\Lambda}{\hat\Lambda}\right|\|\bM\|_{\text{F}}\lesssim \sqrt{\frac{dr}{n}}+\frac{r}{\sum_{j=1}^r\lambda_j}\sqrt{\frac{d}{n}}
\end{align*}
where we've used the fact $\fro{\bM}\le \sum_{j=1}^r\lambda_j$. Take union bound over all events that we've conditioned on in previous steps, we conclude that with probability at least $1-\exp(-c(d/r\wedge n))$:
\begin{align}\label{step5-hpbound}
\min_{\eta\in\{\pm1\}}\fro{\hat \bM-\eta \bM}\lesssim \sqrt{\frac{dr}{n}}+\frac{r}{\sum_{j=1}^r\lambda_j}\sqrt{\frac{d}{n}}
\end{align}
\subsubsection*{Bound in expectation:}
Denote the event $Q:=\{\eqref{step5-hpbound} \text{~holds}\}$. Then we have the following bound in expectation:
\begin{align}\label{expectation-decomp}
\E\min_{\eta\in\{\pm 1\}}\fro{\hat \bM-\eta \bM}=\E\min_{\eta\in\{\pm 1\}}\fro{\hat \bM-\eta \bM}\mathbb{I}_Q+\E\min_{\eta\in\{\pm 1\}}\fro{\hat \bM-\eta \bM}\mathbb{I}_{Q^c}
\end{align}
Note that by Von Neumann's trace inequality, we have
\begin{align}\label{expectation-traceineq}
\Lambda^*&=\sum_{j=1}^r\tilde \bu_j^\top \bM\tilde \bv_j=\text{Tr}(\tilde \bU^\top \bM \tilde \bV)=\text{Tr}(\tilde \bU^\top \bU\bSigma \bV^\top  \tilde \bV)\le  \sum_{j=1}^r\sigma_j(\tilde \bU^\top \bU)\sigma_j(\bSigma \bV^\top  \tilde \bV)\le \sum_{j=1}^r\lambda_j
\end{align}
Since $\check \bM$ is a rank-$r$ projection of $\Lambda^* \bM+\Upsilon^\prime$ and $\hat \Lambda$ is lower bounded by $\frac{\sqrt{d}r}{n^{1/4}}$, we have the following upper bound using \eqref{expectation-traceineq}:
\begin{align}\label{expectation-mhat}
\fro{\hat \bM}=\fro{\hat\Lambda^{-1}\check \bM}\le \frac{\sqrt{r}}{\hat\Lambda}\left(\Lambda^*\op{\bM}+\op{\Upsilon^\prime}\right)\le \frac{n^{1/4}}{\sqrt{dr}}\left(\lambda_1\sum_{j=1}^r\lambda_j+\op{\Upsilon^\prime}\right)
\end{align}
Now we turn to bound $\E \op{\Upsilon^\prime}$, note that by definition we have 
\begin{align}\label{expectation-upsilon-decomp}
\E\op{\Upsilon^\prime}&\le \E\left\|\sum_{j=1}^r\left(\tilde \bu_j^\top \bM\tilde \bv_j\right)\left( \frac{1}{n_0}\sum_{i=1}^{n_0}s_i^{(4)}\bZ_i^{(4)}\right)\right\|+\E\left\|\bM \left(\sum_{j=1}^r\tilde \bu_j^\top\left(\frac{1}{n_0}\sum_{i=1}^{n_0}s_i^{(4)}\bZ_i^{(4)}\right)\tilde \bv_j\right)\right\|\nonumber \\
&+\E\left\|\frac{1}{n_0}\sum_{i=1}^{n_0}\left(\sum_{j=1}^r\tilde \bu_j^\top \bZ_i^{(4)}\tilde \bv_j\right)\bZ_i^{(4)}-\sum_{j=1}^r\tilde \bu_j\tilde \bv_j^\top\right\| 	
\end{align}
The first term of \eqref{expectation-upsilon-decomp} can be bounded as 
\begin{align*}
	 \E\left\|\sum_{j=1}^r\left(\tilde \bu_j^\top \bM\tilde \bv_j\right)\left( \frac{1}{n_0}\sum_{i=1}^{n_0}s_i^{(4)}\bZ_i^{(4)}\right)\right\|\lesssim  \left(\sum_{j=1}^r\lambda_j\right) \sqrt{\frac{d}{n}}
\end{align*}
The second term of \eqref{expectation-upsilon-decomp} can be bounded as 
\begin{align*}
	 \E\left\|\bM \left(\sum_{j=1}^r\tilde \bu_j^\top\left(\frac{1}{n_0}\sum_{i=1}^{n_0}s_i^{(4)}\bZ_i^{(4)}\right)\tilde \bv_j\right)\right\|\lesssim  \lambda_1 \sqrt{\frac{r}{n}}
\end{align*}
For the last term of \eqref{expectation-upsilon-decomp}, recall the decomposition \eqref{step4-decomp2}, we have
\begin{align*}
	\E\left\|\frac{1}{n_0}\sum_{i=1}^{n_0}\left(\sum_{j=1}^r\tilde \bu_j^\top \bZ_i^{(4)}\tilde \bv_j\right)\mathcal{P}^\perp_{\tilde \bU}\bZ_i^{(4)} \mathcal{P}_{\tilde \bV}\right\|&{=}\E\left[\E\left[\left\|\frac{\sqrt{r}}{n_0}\sum_{i=1}^{n_0}g_{i}Z_{i}^{(4)}\right\|\Bigg | \{g_i\}_{i=1}^n\right]\right]\lesssim \E\left[\frac{\sqrt{dr}}{n_0}\sqrt{\sum_{i=1}^{n_0}g_i^2}\right]\lesssim \sqrt{\frac{dr}{n}}
\end{align*}
Similar bounds hold for $\frac{1}{n_0}\sum_{i=1}^{n_0}\left(\sum_{j=1}^r\tilde \bu_j^\top \bZ_i^{(4)}\tilde \bv_j\right)\mathcal{P}^\perp_{\tilde \bU}\bZ_i^{(4)} \mathcal{P}^\perp_{\tilde \bV}$ and $\frac{1}{n_0}\sum_{i=1}^{n_0}\left(\sum_{j=1}^r\tilde \bu_j^\top \bZ_i^{(4)}\tilde \bv_j\right)\mathcal{P}_{\tilde \bU}\bZ_i^{(4)}\mathcal{P}^\perp_{\tilde \bV}$. It remains to find $\E \left\|\frac{1}{n_0}\sum_{i=1}^{n_0}\left (\sum_{j=1}^r\tilde \bu_j^\top \bZ_i^{(4)}\tilde \bv_j\right)\mathcal{P}_{\tilde \bU}\bZ_i^{(4)} \mathcal{P}_{\tilde \bV}-\sum_{j=1}^r\tilde \bu_j\tilde \bv_j^\top\right\|$. For simplicity, denote $\Gamma:=\left\|\frac{1}{n_0}\sum_{i=1}^{n_0}\left (\sum_{j=1}^r\tilde \bu_j^\top \bZ_i^{(4)}\tilde \bv_j\right)\mathcal{P}_{\tilde \bU}\bZ_i^{(4)} \mathcal{P}_{\tilde \bV}-\sum_{j=1}^r\tilde \bu_j\tilde \bv_j^\top\right\|$, then using Lemma \ref{lemma-concentration} we can get
\begin{align*}
\E\Gamma&=\int_{0}^\infty\Prob\left(\Gamma\ge t\right)dt=\int_{0}^{2r\sqrt{\frac{\log(2r)}{n_0}}}\Prob\left(\Gamma\ge t\right)dt+\int_{2r\sqrt{\frac{\log(2r)}{n_0}}}^\infty \Prob\left(\Gamma\ge t\right)dt\\
&\le 2r\sqrt{\frac{\log(2r)}{n_0}}+\int_{2r\sqrt{\frac{\log(2r)}{n_0}}}^\infty \frac{r}{2n_0}\left(\sqrt{\frac{n_0}{u+\log (2r)}}+2\right)\Prob\left(\Gamma\ge r\sqrt{\frac{u+\log(2r)}{n_0}}+r\frac{u+\log(2r)}{n_0}\right)du\\
&\lesssim r\sqrt{\frac{\log r}{n}}+\frac{r}{n}\int_{2r\sqrt{\frac{\log(2r)}{n_0}}}^\infty\left(\sqrt{\frac{n_0}{u+\log (2r)}}+2\right)\exp(-u)du\lesssim  r\sqrt{\frac{\log r}{n}}
\end{align*} 
Hence we can conclude that
\begin{align}
\E\op{\Upsilon^\prime}\le 	\left(\sum_{j=1}^r\lambda_j\right) \sqrt{\frac{d}{n}}+\sqrt{\frac{dr}{n}}
\end{align}
provided that $d\gtrsim r\log r$. By \eqref{expectation-mhat}, we have
\begin{align*}
	\E\fro{\hat \bM}\lesssim  \frac{n^{1/4}}{\sqrt{dr}}\left(\lambda_1\sum_{j=1}^r\lambda_j+\left(\sum_{j=1}^r\lambda_j\right) \sqrt{\frac{d}{n}}+\sqrt{\frac{dr}{n}}\right)\lesssim \frac{n^{1/4}}{\sqrt{dr}}\left(\lambda_1\sum_{j=1}^r\lambda_j\right)
\end{align*}
Hence
\begin{align*}
	\E\min_{\eta\in\{\pm 1\}}\fro{\hat \bM-\eta \bM}\le \E\fro{\hat \bM}+\E\fro{ \bM}\lesssim  \frac{\lambda_1n^{1/4}}{\sqrt{d}}\left(\frac{1}{\sqrt{r}}\sum_{j=1}^r\lambda_j\right)+\sqrt{\sum_{j=1}^r\lambda_j^2}
\end{align*}
Then \eqref{expectation-decomp} implies that 
\begin{align*}
	\E\min_{\eta\in\{\pm 1\}}\fro{\hat \bM-\eta \bM}&\le \sqrt{\frac{dr}{n}}+\frac{r}{\sum_{j=1}^r\lambda_j}\sqrt{\frac{d}{n}}+\left[\frac{\lambda_1n^{1/4}}{\sqrt{d}}\left(\frac{1}{\sqrt{r}}\sum_{j=1}^r\lambda_j\right)+\sqrt{\sum_{j=1}^r\lambda_j^2}\right]\exp(-c(r^{-1}d\wedge n))\\
	&\lesssim \sqrt{\frac{dr}{n}}+\frac{r}{\sum_{j=1}^r\lambda_j}\sqrt{\frac{d}{n}}
\end{align*}
provided that $\lambda_1\asymp\lambda_r\asymp \lambda$ and $\lambda\lesssim \exp(c(r^{-1}{d}\wedge n)-\log n)$.
{\subsection{Proof of Theorem \ref{thm:lower_bound}}
The main idea is to construct a set of sufficiently dissimilar hypotheses to apply Fano's method. To this end, we fix some $\bU_0\in \mathbb{O}_{d,r}$ and consider the ball centered at $\bU_0$ with radius of $\epsilon\in (0,\sqrt{2r}]$ under the chordal Frobenius-norm metric $\textsf{dist}(\bU_1,\bU_2):=\min_{\bO\in\mathbb{O}_r}\fro{\bU_1- \bU_2\bO}$ :
$$B_\epsilon(\bU_0):=\{\bU:\textsf{dist}(\bU,\bU_0)\le \epsilon\}$$
By Lemma 1 in \cite{cai2013sparse} and the equivalence between $\textsf{dist}(\cdot,\cdot)$ and $\|\sin\Theta(\cdot,\cdot)\|$, we have for any $\alpha\in(0,1)$, there exists $\{\bU_i^\prime\}_{i=1}^m$, a packing of $B_\epsilon(\bU_0)$ such that for some absolute constant $c_0>0$:
$$m\ge \left(\frac{c_0}{\alpha}\right)^{r(d-r)},\quad \min_{i< j}\textsf{dist}(\bU_i^\prime,\bU_j^\prime)\ge \alpha \epsilon$$
Denote $\bO_i=
\arg\min_{\bO\in\mathbb{O}_r}\fro{\bU_i^\prime- \bU_0\bO}$. Fix $\bSigma= \text{diag}(\lambda_1,\cdots,\lambda_r)$ with $\lambda_1= \cdots = \lambda_r=\lambda$ and $\bV\in\mathbb{O}_r$, we can construct $\bM_i= \bU_i^\prime \bO_i^\top\bSigma \bV^\top $ for $i=1,\cdots, m$. Notice that
\begin{align*}
\min_{\eta\in\{\pm1\}}\fro{\bM_i-\eta \bM_j}&=\min_{\eta\in\{\pm1\}}\fro{\bU_i^\prime \bO_i^\top\bSigma V^\top  -\eta \bU_j^\prime \bO_j^\top\bSigma \bV^\top }= \lambda \min_{\eta\in\{\pm1\}}\fro{\bU_i^\prime \bO_i^\top -\eta \bU_j^\prime \bO_j^\top}\\
&\ge \lambda \cdot\textsf{dist}(\bU_i^\prime,\bU_j^\prime)\ge \lambda \alpha\epsilon	
\end{align*}
Let $P_{\bM}$ denote the distribution of $\bX=s\bM+\bZ$ and let ${P}_j^{1:n}$ denote the distribution of $\{\bX_i^{(j)}=s_i \bM_j+\bZ_i, i=1,\cdots,n \}$, i.e, the $j$-th model parametrized by $\bM_j$ for $j=1,\cdots,m$. When $\fro{\bSigma}=\sqrt{r}\lambda \ge 1$, since $s$ has a Rademacher prior, using the log-sum inequality (see, e.g., \cite{do2003fast}) we have
\begin{align*}
\text{D}_{\textsf{KL}}({P}_j^{1:n}||{P}_k^{1:n})&\le \sum_{i=1}^n\frac{1}{2}\fro{\bM_j-\bM_k}^2= \frac{1}{2}n\fro{\bU_j^\prime \bO_j^\top\bSigma \bV^\top  - \bU_k^\prime \bO_k^\top\bSigma \bV^\top }^2=\frac{1}{2}n\lambda^2\fro{\bU_j^\prime \bO_j^\top - \bU_k^\prime \bO_k^\top }^2\\
&\le n\lambda^2\left(\textsf{dist}^2(\bU_j^\prime,\bU_0)+\textsf{dist}^2(\bU_k^\prime,\bU_0)\right)\le 	2n\lambda^2\epsilon^2
\end{align*}
When $\fro{\bSigma}=\sqrt{r}\lambda \le 1$, by Lemma 27 in \cite{wu2019randomly}, there exists a universal constant $C>0$, such that for any $\bU,\tilde \bU\in \mathbb{O}_{d,r}$:
\begin{align*}
\text{D}_\textsf{KL}(P_{\bU\bSigma \bV^\top}||P_{\tilde \bU\bSigma \bV^\top})&\le C\min_{\eta\in\{\pm 1\}}\fro{\text{vec}(\fro{\bSigma}^{-1}\bU\bSigma \bV^\top)-\eta\text{vec}(\fro{\bSigma}^{-1}\tilde \bU\bSigma \bV^\top)}^2\fro{\bSigma}^4\\
&\le C\fro{\bSigma}^4\frac{\|\bSigma\|^2}{\fro{\bSigma}^2}\min_{\eta\in\{\pm 1\}}\fro{\bU-\eta\tilde \bU}^2=Cr\lambda^4\min_{\eta\in\{\pm 1\}}\fro{\bU-\eta\tilde \bU}^2	
\end{align*}
which implies that
\begin{align*}
\text{D}_{\textsf{KL}}({P}_j^{1:n}||{P}_k^{1:n})&=n\text{D}_\textsf{KL}(P_{\bU_j^\prime \bO_j^\top\bSigma \bV^\top}||P_{\bU_k^\prime \bO_k^\top\bSigma \bV^\top})\le Cnr\lambda^4\min_{\eta\in\{\pm 1\}}\fro{\bU_j^\prime \bO_j^\top-\eta \bU_k^\prime \bO_k^\top}^2\\
&\le Cnr\lambda^4\left(\fro{\bU_j^\prime \bO_j^\top-\bU_0}^2+\min_{\eta\in\{\pm 1\}}\fro{\bU_0-\eta \bU_k^\prime \bO_k^\top}^2\right)\\
&=Cnr\lambda^4\left(\textsf{dist}^2(\bU_j^\prime,\bU_0)+\textsf{dist}^2(\bU_k^\prime,\bU_0)\right)\le Cnr\lambda^4\epsilon^2
\end{align*}
Hence we have 
$$\text{D}_{\textsf{KL}}({P}_j^{1:n}||{P}_k^{1:n})\le Cn\lambda^2(r\lambda^2\wedge 1)\epsilon^2$$
By Fano's lower bound on minimax risk (see, e.g., Proposition 15.12 in \cite{wainwright2019high}), we have
$$\inf_{\hat \bM}\sup_{\bM\in \mathcal{M}_{d_1,d_2}(r,\lambda)} \E \min_{\eta\in\{\pm 1\}}\fro{\hat \bM-\eta \bM}\ge \lambda\alpha\epsilon\left(1-\frac{Cn\lambda^2(r\lambda^2\wedge 1)\epsilon^2+\log 2}{r(d-r)\log(c_0/\alpha)}\right)$$
By choosing $\epsilon=\sqrt{\frac{r(d-r)}{C_0n\lambda^2(r\lambda^2\wedge 1)}}\wedge \sqrt{2r}$ for some large absolute constant $C_0>0$ and $\alpha=(c_0\wedge 1)/8$, we can guarantee that  $\left(1-\frac{Cn\lambda^2(r\lambda^2\wedge 1)\epsilon^2+\log 2}{r(d-r)\log(c_0/\alpha)}\right)\ge \frac{1}{2}$. Hence
\begin{align*}
\inf_{\hat \bM}\sup_{\bM\in \mathcal{M}_{d_1,d_2}(r,\lambda)} \E \min_{\eta\in\{\pm 1\}}\fro{\hat \bM-\eta \bM}\gtrsim \lambda\left(\sqrt{\frac{dr/n}{\lambda^2(r\lambda^2\wedge 1)}}\wedge \sqrt{r}\right)\gtrsim \left(\frac{1}{\lambda}\sqrt{\frac{d}{n}}+\sqrt{\frac{dr}{n}}\right)\wedge \lambda\sqrt{r}
\end{align*}\qed
}
\subsection{Proof of Theorem \ref{thm:comp}}
Denote the prior distribution for $(\bM,\bs)$ defined in \eqref{eq:testH} as $\Pi$, where $\bs=(s_1,\cdots,s_n)$ is the latent label vector. Let $(\bM^{(1)},\bs^{(1)}),(\bM^{(2)},\bs^{(2)})$ be two independent copies from prior distribution $\Pi$. By Theorem 2.6 in \cite{kunisky2019notes}, we have the following formula for $\|L_n^{\le D}\|$ under the additive Gaussian noise model:
\[\|L_n^{\le D}\|^2={\E_\Pi}\sum_{k=1}^D \frac{1}{k!}\langle \bs^{(1)},\bs^{(2)}\rangle^k\langle \bM^{(1)},\bM^{(2)}\rangle^k=1+\E_\Pi\sum_{k=1}^{\lfloor D/2\rfloor} \frac{1}{(2k)!}\langle \bs^{(1)},\bs^{(2)}\rangle^{2k}\langle \bM^{(1)},\bM^{(2)}\rangle^{2k}\]
The last inequality is due to the fact that $\langle \bs^{(1)},\bs^{(2)}\rangle$ in distribution equals to the sum of $n$ i.i.d. Rademacher random variables, denoted by $\sum_{i=1}^n U_i$, and hence $\E\langle \bs^{(1)},\bs^{(2)}\rangle^k=0$ for odd $k$. Hence we have
\begin{equation}\label{proofs-complower-eq1}
\E \langle \bs^{(1)},\bs^{(2)}\rangle^{2k}=\E \left(\sum_{i=1}^n U_i\right)^{2k}=\E \sum_{2k_1+\cdots+2k_n=2k} U_1^{2k_1}\cdots U_n^{2k_n}={n+k-1\choose k}
\end{equation}
Moreover, $\langle \bM^{(1)},\bM^{(2)} \rangle=\lambda^2 \langle \bu^{(1)},\bu^{(2)}\rangle\langle \bv^{(1)},\bv^{(2)}\rangle=\frac{\lambda^2}{d^2}\left(\sum_{i=1}^d U_i^{(1)}\right) \left(\sum_{i=1}^d U_i^{(2)}\right)$, where for $j=1,2$ $\{U_i^{(j)}\}_{i=1}^d$ are two independent copies of  $d$ i.i.d. Rademacher random variables. Since the even moment of standard normal is lower bounded by $1$, denote $d$ i.i.d. standard normal random variables by $\{g_i\}_{i=1}^d$ and then we have the following simple bound for the combination number:
\begin{align*}
\E\left(\sum_{i=1}^d U_i^{(1)}\right)^{2k}&=\E \sum_{2k_1+\cdots+2k_d=2k} (U_1^{(1)})^{2k_1}\cdots (U_d^{(1)})^{2k_d}\le \E \sum_{2k_1+\cdots+2k_d=2k} g_1^{2k_1}\cdots g_d^{2k_d}\\
&=E \left(\sum_{i=1}^d g_i\right)^{2k}= d^k(2k-1)!!
\end{align*}
Hence we have
\begin{equation}\label{proofs-complower-eq2}
\E \langle \bM^{(1)},\bM^{(2)} \rangle^{2k}=\frac{\lambda^{4k}}{d^{4k}}\E\left(\sum_{i=1}^d U_i^{(1)}\right)^{2k} \E\left(\sum_{i=1}^d U_i^{(2)}\right)^{2k}\le \frac{\lambda^{4k}}{d^{2k}}((2k-1)!!)^2
\end{equation}
Combining \eqref{proofs-complower-eq1} and \eqref{proofs-complower-eq2}, we arrive at
\begin{equation*}
\|L_n^{\le D}\|^2\le 1+\sum_{k=1}^{\lfloor D/2\rfloor}\frac{((2k-1)!!)^2}{(2k)!}{n+k-1\choose k}\frac{\lambda^{4k}}{d^{2k}}\le 1+\sum_{k=1}^{\lfloor D/2\rfloor}{n+k-1\choose k}\frac{\lambda^{4k}}{d^{2k}}=:1+\sum_{k=1}^{\lfloor D/2\rfloor}T_k
\end{equation*}
Notice that
$$\frac{T_{k+1}}{T_k}=\frac{\lambda^4}{d^2}\frac{{n+k\choose k+1}}{{n+k-1\choose k}}=\frac{\lambda^4}{d^2}\frac{n+k}{k+1}=\frac{\lambda^4}{d^2}\left(\frac{n}{k+1}-1\right)\lesssim\frac{\lambda^4 n}{d^2}\le \frac{1}{2} $$
provided that $\lambda^2\lesssim \frac{d}{\sqrt{n}}$. Together with $T_1=\frac{\lambda^4 n}{d^2}$, we have
$$\|L_n^{\le D}\|^2\le 1+O(T_1)=1+O\left(\frac{\lambda^4 n}{d^2}\right)$$\qed
\section{Proofs for technical lemmas}
\subsection{Proof of Lemma \ref{lem:hellinger}}
Our result is an application of the following lemma.
\begin{lemma}[Theorem 1 in \cite{davies2021lower}]\label{lem:hellinger_origin}
	Define 
$$\mathcal{F}:=\left\{f_{\bmu_0,\bmu_1}=\frac{1}{2}\mathcal{N}(\bmu_0,\bSigma)+\frac{1}{2}\mathcal{N}(\bmu_1,\bSigma)|	\bmu_0,\bmu_1\in \mathbb{R}^d,\bSigma\in\mathbb{R}^{d\times d},\bSigma\succ0,\bSigma=\bSigma^\top \right\}$$
For $f_{\bmu_0,\bmu_1},f_{\bmu_0^\prime,\bmu_1^\prime}\in \mathcal{F}$, define sets $S_1=\{\bmu_1-\bmu_0,\bmu_1^\prime-\bmu_0^\prime\}$, $S_2=\{\bmu_0^\prime-\bmu_0,\bmu_1^\prime-\bmu_1\}$, $S_3=\{\bmu_0^\prime-\bmu_1,\bmu_1^\prime-\bmu_0\}$ and vectors $\bv_k=\argmin_{\bs\in S_k}\op{\bs}_2$ for $k=1,2,3$. Let $\lambda_{\bSigma,\mathcal{U}}:=\max_{\bu:\op{\bu}_2=1,\bu\in\mathcal{U}}\bu^\top \bSigma \bu$ with $\mathcal{U}$ being the span of the vectors $\bv_1,\bv_2,\bv_3$. If $\op{\bv_1}_2\ge \min(\op{\bv_2}_2,\op{\bv_3}_2)/2$ and $\sqrt{\lambda_{\bSigma,\mathcal{U}}}=\Omega(\op{\bv_1})$, then 
	$$\op{f_{\bmu_0,\bmu_1}-f_{\bmu_0^\prime,\bmu_1^\prime}}_{\textsf{TV}}=\Omega\left(\min\left(1,\frac{\op{\bv_1}_2\min(\op{\bv_2}_2,\op{\bv_3}_2)}{\lambda_{\bSigma,\mathcal{U}}}\right)\right)$$
	and otherwise, we have that
	$$\op{f_{\bmu_0,\bmu_1}-f_{\bmu_0^\prime,\bmu_1^\prime}}_{\textsf{TV}}=\Omega\left(\min\left(1,\frac{\min(\op{\bv_2}_2,\op{\bv_3}_2)}{\sqrt{\lambda_{\bSigma,\mathcal{U}}}}\right)\right)$$
\end{lemma}
\noindent In our setting, $\bmu_0=\text{vec}(\bM)$, $\bmu_1=-\text{vec}(\bM)$, $\bSigma=\bI_{d_1}\otimes\bI_{d_2}$ and $\mathcal{F}=\mathcal{M}_{d_1,d_2}(r,\lambda)$. For any $p_\bM,p_{\bM_1}\in \mathcal{F}=\mathcal{M}_{d_1,d_2}(r,\lambda)$, we have $\op{\bv_1}_2=2\max\{\fro{\bM},\fro{\bM_1}\}$, $\op{\bv_2}_2=\fro{\bM-\bM_1}$, $\op{\bv_3}_2=\fro{\bM+\bM_1}$ and $\lambda_{\bSigma,\mathcal{U}}=1$. Notice that $\op{\bv_1}_2=2\max\{\fro{\bM},\fro{\bM_1}\}\geq  \ell(\bM_1,\bM)/2=\min(\op{\bv_2}_2,\op{\bv_3}_2)/2$ always holds.  By Lemma \ref{lem:hellinger_origin}, if $\op{\bv_1}_2\asymp  \fro{\bM}+\fro{\bM_1}\lesssim 1$, then 
$$d_{\textsf{TV}}(p_{\bM},p_{\bM_1})\gtrsim(\fro{\bM}+\fro{\bM_1})\ell(\bM_1,\bM)$$
Otherwise, we have
$$d_{\textsf{TV}}(p_{\bM},p_{\bM_1})\gtrsim \min\{1,\ell(\bM_1,\bM)\}$$
The result immediately follows by noting that the total variation distance is bounded by Hellinger distance.
\subsection{Proof of Lemma \ref{lem:KL}}
By definition of KL divergence, we have
\begin{align*}
D_{\textsf{KL}}(p_{\bM}\|p_{\bM_1})=\E\log \frac{p_{\bM}(\bX)}{p_{\bM_1}(\bX)}&=\frac{1}{2}\fro{\bM_1}^2-\frac{1}{2}\fro{\bM}^2+\E\log \left(\frac{e^{\langle \bX,\bM \rangle}+e^{-\langle \bX,\bM \rangle}}{e^{\langle \bX,\bM_1\rangle}+e^{-\langle \bX,\bM_1\rangle}}\right)
\end{align*}
where $\bX\sim p_{\bM}$. By log-sum-exp inequality, we have 
\begin{align*}
\log \left(\frac{e^{\langle \bX,\bM \rangle}+e^{-\langle \bX,\bM \rangle}}{e^{\langle \bX,\bM_1 \rangle}+e^{-\langle \bX,\bM_1 \rangle}}\right)\geq |\langle \bX,\bM\rangle|-|\langle \bX,\bM_1\rangle|-\log 2
\end{align*}
It follows that
\begin{align*}
D_{\textsf{KL}}(p_{\bM}\|p_{\bM_1})\ge \frac{1}{2}\fro{\bM_1}^2-\frac{1}{2}\fro{\bM}^2+\E\left[|\langle \bX,\bM\rangle|-|\langle \bX,\bM_1\rangle|\right]-\log 2
\end{align*}
Recall that $\bX\overset{d}{=}s\bM+\bZ$, and for brevity denote $\E_{\bx}$ the expectation over $\bx$. Then we have that 
\begin{align*}
	\E|\langle \bX,\bM\rangle|&=\E_s\E_\bZ|\langle s\bM+\bZ,\bM\rangle|=\E_s\E_\bZ|s\fro{\bM}^2+\langle \bZ,\bM\rangle|\\
	&=\E_s\left[\fro{\bM}\sqrt{\frac{2}{\pi}}e^{-\frac{\fro{\bM}^2}{2}}+s\fro{\bM}^2(1-2\Phi\left(-s\fro{\bM}\right))\right]\\
	&=\fro{\bM}\sqrt{\frac{2}{\pi}}e^{-\frac{\fro{\bM}^2}{2}}+\fro{\bM}^2(\Phi\left(\fro{\bM}\right)-\Phi\left(-\fro{\bM}\right))
\end{align*}
where the third equality is due to $|s\fro{\bM}^2+\langle\bZ,\bM\rangle|\big|s\sim |\mathcal{N}(s\fro{\bM}^2,\fro{\bM}^2)|$ and $\Phi(\cdot)$ denotes the cumulative distribution function (cdf) of standard normal. Likewise, we obtain that
\begin{align*}
	\E|\langle \bX,\bM_1\rangle|=\fro{\bM_1}\sqrt{\frac{2}{\pi}}e^{-\frac{\langle\bM,\bM_1\rangle^2}{2\fro{\bM_1}^2}}+\langle\bM,\bM_1\rangle\left(\Phi\left(\frac{\langle\bM,\bM_1\rangle}{\fro{\bM_1}}\right)-\Phi\left(-\frac{\langle\bM,\bM_1\rangle}{\fro{\bM_1}}\right)\right)
\end{align*}
Thus we have that
\begin{align}\label{ineq:KL-lower}
D_{\textsf{KL}}(p_{\bM}\|p_{\bM_1})&\ge \frac{1}{2}\fro{\bM_1}^2-\frac{1}{2}\fro{\bM}^2+\fro{\bM}\sqrt{\frac{2}{\pi}}e^{-\frac{\fro{\bM}^2}{2}}+\fro{\bM}^2(\Phi\left(\fro{\bM}\right)-\Phi\left(-\fro{\bM}\right))\nonumber\\
&-\fro{\bM_1}\sqrt{\frac{2}{\pi}}e^{-\frac{\langle\bM,\bM_1\rangle^2}{2\fro{\bM_1}^2}}-\langle\bM,\bM_1\rangle\left(\Phi\left(\frac{\langle\bM,\bM_1\rangle}{\fro{\bM_1}}\right)-\Phi\left(-\frac{\langle\bM,\bM_1\rangle}{\fro{\bM_1}}\right)\right)-\log 2
\end{align}
Without loss of generality, we assume $\langle\bM,\bM_1\rangle>0$. Using the upper and lower bound for cdf of standard normal (see, e.g., \cite{abramowitz1948handbook}), we obtain
\begin{align}\label{ineq:phiM}
\Phi\left(\fro{\bM}\right)-\Phi\left(-\fro{\bM}\right)=1-2\Phi\left(-\fro{\bM}\right)\ge 1-2\sqrt{\frac{2}{\pi}}\frac{1}{\fro{\bM}+\sqrt{\fro{\bM}^2+8/\pi}}e^{-\frac{\fro{\bM}^2}{2}}
\end{align}
\begin{align}\label{ineq:phiMM1}
\Phi\left(\frac{\langle\bM,\bM_1\rangle}{\fro{\bM_1}}\right)-\Phi\left(-\frac{\langle\bM,\bM_1\rangle}{\fro{\bM_1}}\right)&\leq 1-2\sqrt{\frac{2}{\pi}}\frac{1}{\frac{\langle\bM,\bM_1\rangle}{\fro{\bM_1}}+\sqrt{\frac{\langle\bM,\bM_1\rangle^2}{\fro{\bM_1}^2}+4}}e^{-\frac{\langle\bM,\bM_1\rangle^2}{2\fro{\bM_1}^2}}
\end{align}
It follows from \eqref{ineq:KL-lower} \eqref{ineq:phiM} and \eqref{ineq:phiMM1} that
\begin{align*}
D_{\textsf{KL}}(p_{\bM}\|p_{\bM_1})&\ge \frac{1}{2}\fro{\bM_1}^2+\frac{1}{2}\fro{\bM}^2-2\sqrt{\frac{2}{\pi}}\frac{\fro{\bM}^2e^{-\frac{\fro{\bM}^2}{2}}}{\fro{\bM}+\sqrt{\fro{\bM}^2+8/\pi}}-\fro{\bM_1}\sqrt{\frac{2}{\pi}}e^{-\frac{\langle\bM,\bM_1\rangle^2}{2\fro{\bM_1}^2}}\nonumber\\
&-\langle\bM,\bM_1\rangle\left(1-2\sqrt{\frac{2}{\pi}}\frac{e^{-\frac{\langle\bM,\bM_1\rangle^2}{2\fro{\bM_1}^2}}}{\frac{\langle\bM,\bM_1\rangle}{\fro{\bM_1}}+\sqrt{\frac{\langle\bM,\bM_1\rangle^2}{\fro{\bM_1}^2}+4}}\right)-\log 2\\
&= \frac{1}{2}(1-\epsilon)\left(\fro{\bM_1}^2+\fro{\bM}^2-2\langle\bM,\bM_1\rangle\right)+\frac{1}{2}\epsilon\fro{\bM}^2-2\sqrt{\frac{2}{\pi}}\frac{\fro{\bM}^2e^{-\frac{\fro{\bM}^2}{2}}}{\fro{\bM}+\sqrt{\fro{\bM}^2+8/\pi}}\\
&-\langle\bM,\bM_1\rangle\left(\epsilon-2\sqrt{\frac{2}{\pi}}\frac{e^{-\frac{\langle\bM,\bM_1\rangle^2}{2\fro{\bM_1}^2}}}{\frac{\langle\bM,\bM_1\rangle}{\fro{\bM_1}}+\sqrt{\frac{\langle\bM,\bM_1\rangle^2}{\fro{\bM_1}^2}+4}}\right)-\log 2
\end{align*}
where $\epsilon:=\frac{1}{\fro{\bM_1}}\sqrt{\frac{2}{\pi}}e^{-\frac{\langle\bM,\bM_1\rangle^2}{2\fro{\bM_1}^2}}$. Observe that
\begin{align*}
	&\frac{\fro{\bM}^2}{2\fro{\bM_1}}\sqrt{\frac{2}{\pi}}e^{-\frac{\langle\bM,\bM_1\rangle^2}{2\fro{\bM_1}^2}}-\langle\bM,\bM_1\rangle\left(\frac{1}{\fro{\bM_1}}\sqrt{\frac{2}{\pi}}e^{-\frac{\langle\bM,\bM_1\rangle^2}{2\fro{\bM_1}^2}}-2\sqrt{\frac{2}{\pi}}\frac{e^{-\frac{\langle\bM,\bM_1\rangle^2}{2\fro{\bM_1}^2}}}{\frac{\langle\bM,\bM_1\rangle}{\fro{\bM_1}}+\sqrt{\frac{\langle\bM,\bM_1\rangle^2}{\fro{\bM_1}^2}+4}}\right)\\
	&\ge \sqrt{\frac{2}{\pi}}e^{-\frac{\langle\bM,\bM_1\rangle^2}{2\fro{\bM_1}^2}}\left[\frac{\fro{\bM}^2}{2\fro{\bM_1}}-\langle\bM,\bM_1\rangle\left(\frac{1}{\fro{\bM_1}}-\frac{1}{\fro{\bM}+1}\right)\right]
\end{align*}
Now we need to show
\begin{align}\label{ineq:need-to-show}
\frac{\fro{\bM}^2}{2\fro{\bM_1}}-\langle\bM,\bM_1\rangle\left(\frac{1}{\fro{\bM_1}}-\frac{1}{\fro{\bM}+1}\right)\ge 0
\end{align}
It suffices to show
$$\frac{1}{2}\fro{\bM}\ge \fro{\bM_1}\frac{\fro{\bM}+1-\fro{\bM_1}}{\fro{\bM}+1}$$
If $\fro{\bM}+1\le \fro{\bM_1}$, then the inequality is trivial.
If $\fro{\bM}+1=K\fro{\bM_1}$ for some constant $K>1$, then 
$$ \fro{\bM_1}\frac{\fro{\bM}+1-\fro{\bM_1}}{\fro{\bM}+1}=\frac{K-1}{K}\fro{\bM_1}\le \frac{K}{2}\fro{\bM_1}-\frac{1}{2}=\frac{1}{2}\fro{\bM}$$
as long as 
$$\frac{K^2-2K+2}{2K}\fro{\bM_1}\ge \frac{1}{2}$$
Since $K^2-2K+2>0$, the inequality holds provided that $\fro{\bM}\ge \frac{2K-2}{K^2-2K+2}$. If $\fro{\bM_1}=o(\fro{\bM})$, then 
$$ \fro{\bM_1}\frac{\fro{\bM}+1-\fro{\bM_1}}{\fro{\bM}+1}\le \fro{\bM_1}\le \frac{1}{2}\fro{\bM}$$
Therefore, we conclude that \eqref{ineq:need-to-show} holds and hence we obtain that
\begin{align*}
D_{\textsf{KL}}(p_{\bM}\|p_{\bM_1})&\ge \frac{1}{2}(1-\epsilon)\left(\fro{\bM_1}^2+\fro{\bM}^2-2\langle\bM,\bM_1\rangle\right)-2\sqrt{\frac{2}{\pi}}\frac{\fro{\bM}^2e^{-\frac{\fro{\bM}^2}{2}}}{\fro{\bM}+\sqrt{\fro{\bM}^2+8/\pi}}-\log2\\
&\ge c_0\fro{\bM-\bM_1}^2
\end{align*}
provided that  $\fro{\bM}\ge C_0$ and $\fro{\bM-\bM_1}\ge C_1$. Due to symmetry, we can apply the same argument to $D_{\textsf{KL}}(p_{\bM}\|p_{-\bM_1})=D_{\textsf{KL}}(p_{\bM}\|p_{\bM_1})$ and the proof is completed.\qed
\subsection{Proof of Lemma \ref{lemma-concentration}}
The idea is to apply Berstein's type matrix inequality, e.g., Proposition 2 in \cite{koltchinskii2011nuclear}, and check the conditions therein are satisfied. To begin with, it's easy to verify that $\E(\text{Tr}(Z)Z-I_r)=0$. Then we check that $\|\text{Tr}(Z)Z-I_r\|$ is sub-exponential, which can be seen via the following derivation:
\begin{align*}
	\left\|\|\text{Tr}(\bZ)\bZ-\bI_r\|\right\|_{\psi_1}\le \left\|\|\text{Tr}(\bZ)\|\right\|_{\psi_2}\left\|\|\bZ\|\right\|_{\psi_2}+1\lesssim r
\end{align*}
where the second inequality follows from the fact that $\text{Tr}(\bZ)\sim N(0,r)$ and $\Prob(|\|\bZ\|-2\sqrt{r}|\ge t)\le 2\exp(-t^2/2)$. In addition, we need to bound $\left\|\frac{1}{n}\sum_{i=1}^n \E (\text{Tr}\left(\bZ_i\right)\bZ_i-\bI_r)(\text{Tr}\left(\bZ_i\right)\bZ_i-\bI_r)^\top \right\|^{1/2}$. Notice that
\begin{align*}
	\left\|\frac{1}{n}\sum_{i=1}^n \E (\text{Tr}\left(\bZ_i\right)\bZ_i-\bI_r)(\text{Tr}\left(\bZ_i\right)\bZ_i-\bI_r)^\top \right\|=\left\|\E \text{Tr}^2\left(\bZ\right)\bZ\bZ^\top-\bI_r\right\|
\end{align*}
The $l$-th diagonal entry of $\E \text{Tr}^2\left(\bZ\right)\bZ\bZ^\top$ can be computed as
\begin{align*}
	\E [\text{Tr}^2\left(\bZ\right)\bZ\bZ^\top]_{ll}=\E\left[\bZ_{ll}^4+\left(\sum_{j\ne l}\bZ_{jj}^2\right)\left(\sum_{j\ne l}\bZ_{lj}^2\right)+\bZ_{ll}^2\sum_{j\ne l}\bZ_{lj}^2+\bZ_{ll}^2\sum_{j\ne l}\bZ_{jj}^2\right]=r^2+2
\end{align*}
For $(l_1,l_2)$-th entry of $\E \text{Tr}^2\left(\bZ\right)\bZ\bZ^\top$ such that $l_1\ne l_2$, we have
\begin{align*}
	\E [\text{Tr}^2\left(\bZ\right)\bZ\bZ^\top]_{l_1l_2}=\E\left[\left(\sum_{j=1}^r\bZ_{jj}^2\right)\left(\sum_{j=1}^r \bZ_{l_1j}\bZ_{l_2j}\right)+2\left(\sum_{i<j}\bZ_{ii}\bZ_{jj}\right)\left(\sum_{j=1}^r\bZ_{l_1j}\bZ_{l_2 j}\right)\right]=0
\end{align*}
Hence we have
\begin{align*}
	\left\|\frac{1}{n}\sum_{i=1}^n \E (\text{Tr}\left(\bZ_i\right)\bZ_i-\bI_r)(\text{Tr}\left(\bZ_i\right)\bZ_i-\bI_r)^\top \right\|^{1/2}\lesssim r
\end{align*}
Applying matrix Berstein's inequality, we complete the proof.\qed
\subsection{Proof of Lemma \ref{lemma-lowrank-asymmetric}} 
We first prove a symmetric version of this lemma and then extend it to the desired non-symmetric version using standard dilation technique. Now we restate the symmetric version. 
\begin{lemma}\label{lemma-lowrank-symmetric}
	Consider a rank-$r$ matrix $\bM\in\mathbb{R}^{d\times d}$ with eigen-decomposition $\bU\Lambda \bU^\top $, where $\bU\in\mathbb{O}_{d_1,r}$ and $\bLambda=\text{diag}(\lambda_1,\cdots,\lambda_r)$, $|\lambda_1|\geq|\lambda_2|\ge\cdots |\lambda_r|>0$, let $\bE$ be a $d\times d$ symmetric 	perturbation matrix and $\hat \bM=\bM+\bE$. Denote $\hat \bM_r$ the best rank-$r$ approximation of $\hat \bM$. Suppose that $|\lambda_r|\ge (4+c_0)\|\bE\|$ for any constant $c_0>0$, then there exists some absolute constant $C_0>0$ such that
	\begin{align*}
		\fro{\hat \bM_r-\bM}\le C_0\sqrt{r}\|\bE\|
	\end{align*}
\end{lemma}
\noindent \textit{Proof of Lemma \ref{lemma-lowrank-symmetric}}
\\Denote $\hat \bU$ the leading $r$ (in absolute value) eigenvectors of $\hat \bM$. First, by Theorem 1 in \cite{xia2021normal}, we have the following identity holds:
\begin{align}\label{lemma-identity}
	\hat \bU\hat \bU^\top -\bU\bU^\top = \sum_{k\ge 1}\mathcal{S}_{\bM,k}(\bE)
\end{align}
where
$$\mathcal{S}_{\bM,k}(\bE)=\sum_{\bs:s_1+\cdots+s_{k+1}=k}(-1)^{1+\tau(\bs)}\mathfrak{P}^{-s_1}\bE\mathfrak{P}^{-s_2}\cdots \bE\mathfrak{P}^{-s_{k+1}}$$
Here $\bs=(s_1,\cdots,s_{k+1})$ contains non-negative indices, $\tau(\bs)=\sum_{j=1}^{k+1}\mathbb{I}(s_j>0)$ is the number of positive indices in $\bs$ and $\mathfrak{P}^{-1}=\bU\bLambda^{-1}\bU^\top$ and $\mathfrak{P}^0=\bU_\perp \bU_\perp^\top $. By definition of $\hat \bM_r$, utilizing \eqref{lemma-identity} we have
\begin{align}\label{lemma-decomp}
	\fro{\hat \bM_r-\bM}&=\fro{\hat \bU\hat \bU^\top (\bM+\bE)\hat \bU\hat \bU^\top -\bM}\nonumber\\
	&\le \fro{(\hat \bU\hat \bU^\top -\bU\bU^\top)\bM(\hat \bU\hat \bU^\top -\bU\bU^\top)}\\
	&+\fro{(\hat \bU\hat \bU^\top -\bU\bU^\top)\bM+\bM(\hat \bU\hat \bU^\top -\bU\bU^\top)+\hat \bU\hat \bU^\top \bE\hat \bU\hat \bU^\top}\nonumber\\
	&{\le} \fro{(\hat \bU\hat \bU^\top -\bU\bU^\top)\bM(\hat \bU\hat \bU^\top -\bU\bU^\top)}+\fro{\mathcal{S}_{\bM,1}(\bE)\bM+\bM\mathcal{S}_{\bM,1}(\bE)+P_\bU \bE P_{ \bU}}\nonumber\\
	&+\fro{\sum_{k\ge 2}\mathcal{S}_{\bM,k}(\bE)\bM+\bM\sum_{k\ge 2}\mathcal{S}_{\bM,k}(\bE)}+\fro{P_{\hat \bU} EP_{\hat \bU}-P_\bU \bE P_{ \bU}}
\end{align}
We are going to bound each term of \eqref{lemma-decomp}. Notice that for any $k\ge 1$
\begin{align*}
\|\mathcal{S}_{\bM,k}(\bE)\|\le \sum_{\bs:s_1+\cdots+s_{k+1}=k}\|\mathfrak{P}^{-s_1}\bE\mathfrak{P}^{-s_2}\cdots \bE\mathfrak{P}^{-s_{k+1}}\|\le {2k \choose k}\left(\frac{\|\bE\|}{|\lambda_r|}\right)^{k}\le \left(\frac{4\|\bE\|}{|\lambda_r|}\right)^{k}
\end{align*}
\begin{align*}
\op{\mathcal{S}_{\bM,k}(\bE)M}&=\op{\sum_{\bs:s_1+\cdots+s_{k+1}=k}(-1)^{1+\tau(\bs)}\mathfrak{P}^{-s_1}\bE\mathfrak{P}^{-s_2}\cdots \bE\mathfrak{P}^{-s_{k+1}}\bU\bLambda \bU^\top }\\
&\overset{(a)}{\le} \sum_{\bs:s_1+\cdots+s_{k+1}=k,s_{k+1}>0}\op{\mathfrak{P}^{-s_1}\bE\mathfrak{P}^{-s_2}\cdots \mathfrak{P}^{-s_k}\bE\bU\Lambda^{-s_{k+1}+1}}\\
&\le {2k \choose k}\op{\bE}\left(\frac{\|\bE\|}{|\lambda_r|}\right)^{k-1}\lesssim \op{\bE}\left(\frac{4\|\bE\|}{|\lambda_r|}\right)^{k-1}
\end{align*}
where in (a) we used the fact $\mathfrak{P}^{0}U=U_\perp U_\perp^\top U=0 $. Therefore, for the first term of \eqref{lemma-decomp} we have
\begin{align*}
\op{(\hat \bU\hat \bU^\top -\bU\bU^\top)\bM(\hat \bU\hat \bU^\top -\bU\bU^\top)}&=\op{\sum_{k_1,k_2\ge 1}\mathcal{S}_{\bM,k_1}(\bE)\bM\mathcal{S}_{\bM,k_2}(\bE)}\le \sum_{k_1,k_2\ge 1}\op{\mathcal{S}_{\bM,k_1}(\bE)\bM}\|\mathcal{S}_{\bM,k_2}(\bE)\|\\
&\le  \sum_{k_1\ge 1}\op{\mathcal{S}_{\bM,k_1}(\bE)\bM}\|\mathcal{S}_{\bM,1}(\bE)\|+\sum_{k_1\ge 1,k_2\ge 2}\op{\mathcal{S}_{\bM,k_1}(\bE)\bM}\|\mathcal{S}_{\bM,k_2}(\bE)\|\\
&\lesssim \frac{\op{\bE}^2}{|\lambda_r|}\sum_{k_1\ge 1}\left(\frac{4\|\bE\|}{|\lambda_r|}\right)^{k_1-1}+\op{\bE}\sum_{k_1\ge1,k_2\ge 2}\left(\frac{4\|\bE\|}{|\lambda_r|}\right)^{k_1+k_2-1}\\
&\lesssim \frac{\op{\bE}^2}{|\lambda_r|}\lesssim \op{\bE}
\end{align*}
Since $\text{rank}\left((\hat \bU\hat \bU^\top -\bU\bU^\top)\bM(\hat \bU\hat \bU^\top -\bU\bU^\top)\right)\le 2r$, we have
\begin{align*}
	\fro{(\hat \bU\hat \bU^\top -\bU\bU^\top)\bM(\hat \bU\hat \bU^\top -\bU\bU^\top)}\le \sqrt{2r}\op{(\hat \bU\hat \bU^\top -\bU\bU^\top)\bM(\hat \bU\hat \bU^\top -\bU\bU^\top)}\lesssim \sqrt{r}\op{\bE}
\end{align*}
The second term of \eqref{lemma-decomp} can bounded as
\begin{align*}
\fro{\mathcal{S}_{\bM,1}(\bE)\bM+\bM\mathcal{S}_{\bM,1}(\bE)+P_\bU \bE P_{ \bU}}&=\fro{\mathfrak{P}^{0}\bE\mathfrak{P}^{-1}\bU\bLambda \bU^\top +\bU\bLambda \bU^\top\mathfrak{P}^{-1}\bE\mathfrak{P}^{0} +\bU\bU^\top \bE \bU\bU^\top}\\
&=\fro{\bU_\perp \bU_\perp^\top \bE\bU\bU^\top +\bU\bU^\top \bE \bU_\perp \bU_\perp^\top +\bU\bU^\top \bE \bU\bU^\top}\\
&\lesssim \sqrt{r}\op{\bE}
\end{align*}
For the third term in \eqref{lemma-decomp}, we have
\begin{align*}
	\fro{\sum_{k\ge 2}\mathcal{S}_{\bM,k}(\bE)\bM+\bM\sum_{k\ge 2}\mathcal{S}_{\bM,k}(\bE)}\le 2\sum_{k\ge 2}\fro{\mathcal{S}_{\bM,k}(\bE)\bM}\lesssim \sqrt{r}\op{\bE}\sum_{k\ge 2}\left(\frac{4\|\bE\|}{|\lambda_r|}\right)^{k-1}\lesssim \sqrt{r}\op{\bE}
\end{align*}
It remains to bound the last term of \eqref{lemma-decomp}, which can be done as follows
\begin{align*}
	\fro{P_{\hat \bU} \bE P_{\hat \bU}-P_\bU \bE P_{\bU}}&=\fro{(\hat \bU\hat \bU^\top-\bU\bU^\top) \bE \hat \bU\hat \bU^\top +\bU\bU^\top \bE(\hat \bU\hat \bU^\top-\bU\bU^\top)}\le 2\fro{(\hat \bU\hat \bU^\top-\bU\bU^\top) \bE}\\
	&\le 2\sqrt{r}\op{\bE}\sum_{k\ge 1}\op{\mathcal{S}_{\bM,k}(\bE)}\lesssim \sqrt{r}\op{\bE}
\end{align*}
Collecting all pieces, by \eqref{lemma-decomp} we arrive at
\begin{align*}
	\fro{\hat \bM_r-\bM}\lesssim \sqrt{r}\op{\bE}
\end{align*}\qed
\\\textit{Proof of Lemma \ref{lemma-lowrank-asymmetric}}
\\Now we turn to the proof of Lemma \ref{lemma-lowrank-asymmetric}. Define
$$
\bM:=\begin{bmatrix}
    0     & \bM \\
    \bM^\top       & 0
\end{bmatrix},\quad 
\hat \bM=\begin{bmatrix}
    0     & \hat \bM\\
    \hat \bM^\top       & 0
\end{bmatrix},\quad
\hat \bM_r^*=\begin{bmatrix}
    0     & \hat \bM_r\\
    \hat \bM_r^\top       & 0
\end{bmatrix},\quad 
\bE^*=
\begin{bmatrix}
   0    & \bE \\
    \bE^\top     & 0
\end{bmatrix}$$
Also define
$$
\bTheta=\frac{1}{\sqrt{2}}\begin{bmatrix}
     \bU     & \bU \\
     \bV       & -\bV
\end{bmatrix},\quad 
\hat\bTheta=\frac{1}{\sqrt{2}}\begin{bmatrix}
    \hat \bU     & \hat \bU \\
    \hat \bV       & -\hat \bV
\end{bmatrix}
$$
Notice that $\bTheta$ and $\hat \bTheta$ are the eigenvectors of $\bM$ and $\hat \bM$, respectively. By construction we have $|\lambda_{2r}(\bM)|=\sigma_r$ and $\op{\bE^*}=\op{\bE}$. Then applying Lemma \ref{lemma-lowrank-symmetric} we have
\begin{align*}
	\fro{\hat \bM_r-\bM}=\frac{1}{\sqrt{2}}\left\|
	\begin{bmatrix}
    0     & \hat \bM_r-\bM \\
    \hat \bM_r^\top -\bM^\top       & 0
\end{bmatrix}\right\|_{\text{F}}=\frac{1}{\sqrt{2}}\fro{\hat \bM_r-\bM}\lesssim \sqrt{r}\|\bE\|
\end{align*}
\qed
\subsection{Proof of Lemma \ref{lem:bracketing-entropy}}
Consider a $\epsilon$-net for $\mathcal{M}_{d_1,d_2}(r):=\{\bM\in \mathbb{R}^{d_1\times d_2}: \text{rank}(\bM)=r\}$ endowed with metric $\fro{\cdot }$, denoted by $\mathcal{N}_\eps(\mathcal{M}_{d_1,d_2}(r))=\{\bM_1,\bM_2,\cdots,\bM_N\}$, we have its cardinality $|\mathcal{N}_\eps(\mathcal{M}_{d_1,d_2}(r))|=N\le 
\left(\frac{5}{\epsilon}\right)^{(d_1+d_2)r}$ (see, e.g., \cite{zhang2018tensor}). Then for any $i\in [N]$, we can have a ball centered at $\bM_i$ with radius $\epsilon$, that is, $\mathcal{B}_\epsilon(\bM_i):=\{\bM\in \mathcal{M}_{d_1,d_2}(r):\fro{\bM-\bM_i}\le \epsilon\}$. Hence $\mathcal{M}_{d_1,d_2}(r)\subseteq \cup_{i=1}^n\mathcal{B}_\epsilon(\bM_i)$. Now for any $p_\bM\in \mathcal{P}_{d_1,d_2}(r,\lambda)$ with $p_\bM(X)=(2\pi)^{-d_1d_2/2}\exp\left(-\frac{1}{2}\fro{\bX-\bM}^2\right)$, there exists $j\in[N]$ such that $\bM\in \mathcal{B}_\epsilon(\bM_j)$, we consider the following functions with $\delta=\epsilon/\sqrt{d_1d_2}$:
$$ \begin{cases} 
      l_{\bM_j}(\bX)= \exp\left(-\frac{1}{2}(1+\frac{1}{\delta})\epsilon^2\right)(2\pi)^{-d_1d_2/2}\exp\left(-\frac{\fro{\bX-\bM_j}^2}{2(1+\delta)^{-1}}\right) \\
      u_{\bM_j}(\bX)= \exp\left(\frac{\epsilon^2}{2\delta}\right)(2\pi)^{-d_1d_2/2}\exp\left(-\frac{\fro{\bX-\bM_j}^2}{2(1+\delta)}\right)
   \end{cases}
$$
We first check the bracket $[l_{\bM_j},u_{\bM_j}]$ contains $p_\bM$, which follow from the following observation:
$$\fro{\bX-\bM}^2=\fro{\bX-\bM_j+\bM_j-\bM}^2\le (1+\delta)\fro{\bX-\bM_j}^2+(1+\delta^{-1})\epsilon^2$$
$$\fro{\bX-\bM}^2=\fro{\bX-\bM_j+\bM_j-\bM}^2\geq (1+\delta)^{-1}\fro{\bX-\bM_j}^2-\delta^{-1}\epsilon^2$$
where the inequality follows from the inequality $(a+b)^2\le (1+\delta)a^2+(1+\delta^{-1})b^2$ for any $\delta>0$ and the fact that $\bM\in \mathcal{B}_\epsilon(\bM_j)$. Hence we have $l_{\bM_j}(\bX)\le p_\bM(X) \le u_{\bM_j}(\bX)$. It remains to calculate $d_{\textsf{H}}(l_{\bM_j},u_{\bM_j})$. Note that by definition of Hellinger distance, we have
\begin{align*}
d_{\textsf{H}}^2(l_{\bM_j},u_{\bM_j})&=\exp\left(-\frac{\delta+1}{2\delta}\epsilon^2\right)+\exp\left(\frac{\epsilon^2}{2\delta}\right)-2\exp\left(-\frac{\delta+1}{4\delta}\epsilon^2\right)\exp\left(\frac{\epsilon^2}{4\delta}\right)\left(\frac{2}{1+\delta+(1+\delta)^{-1}}\right)^{d_1d_2/2}\\
&\le 2\cosh\left(\frac{\epsilon^2}{2\delta }\right)-2\exp\left(-\frac{\epsilon^2}{4}\right)[\cosh\left(\ln (1+\delta ) \right)]^{-d_1d_2/2}\\
&\le 2\left(1+\frac{\epsilon^4}{4\delta^2}\right)-2\left(1-\frac{\epsilon^2}{4}\right)\left(1-\frac{\delta^2 d_1d_2}{4}\right)\\
&\le \frac{\epsilon^4}{2\delta^2}+\frac{\epsilon^2}{2}+\frac{\epsilon^2\delta^2d_1d_2}{8}+\frac{\delta^2d_1d_2}{2}\lesssim \epsilon^2d_1d_2
\end{align*}
Hence we can take $\epsilon=\epsilon^\prime /\sqrt{d_1d_2}$, then $d_{\textsf{H}}(l_{M_j},u_{M_j})\le \epsilon^{\prime}$. Since $d_1\asymp d_2\asymp d$, the cardinality of brackets becomes 
$$\log N\le (d_1+d_2)r\log \left(\frac{5\sqrt{d_1d_2}}{\epsilon^\prime}\right)\lesssim dr\log \left(\frac{d}{\epsilon^\prime}\right)$$
\qed

\end{document}